\pdfoutput=1
\documentclass{article}

\usepackage[utf8]{inputenc}
\usepackage[T1]{fontenc}
\usepackage[margin=1in]{geometry}
\usepackage{amsmath,amssymb,amsthm}
\usepackage{booktabs}
\usepackage{array}
\usepackage{multirow}
\usepackage{caption}
\usepackage{subcaption}
\usepackage{enumitem}
\usepackage{graphicx}
\usepackage{microtype}
\usepackage[numbers,sort&compress]{natbib}
\usepackage{xcolor}
\usepackage{hyperref}
\usepackage{algorithm}
\usepackage{algpseudocode}
\usepackage{nicefrac}

\hypersetup{
  colorlinks=true,
  linkcolor=blue!45!black,
  citecolor=blue!45!black,
  urlcolor=blue!55!black,
  pdftitle={Robin-Neumann Coupling of PINN and FEM Solvers: A
    Steklov-Poincare View, with Application to Fluid-Structure
    Interaction with Contact},
  pdfauthor={Mikel Landajuela}
}

\newtheorem{theorem}{Theorem}

\newtheorem{corollary}[theorem]{Corollary}
\theoremstyle{definition}

\theoremstyle{remark}
\newtheorem{remark}{Remark}

\newcommand{\R}{\mathbb{R}}

\newcommand{\tagF}{\mathrm{F}}    
\newcommand{\tagN}{\mathrm{P}}    
\newcommand{\OmegaF}{\Omega_{\tagF}}
\newcommand{\OmegaN}{\Omega_{\tagN}}
\newcommand{\Omegaf}{\Omega_{\!f}}

\newcommand{\SigmaI}{\Sigma}
\newcommand{\SF}{S_{\tagF}}
\newcommand{\SN}{S_{\tagN}}
\newcommand{\bSF}{\mathbf{S}_{\tagF}}
\newcommand{\bSN}{\mathbf{S}_{\tagN}}
\newcommand{\sigmaF}{\sigma_{\tagF}}
\newcommand{\sigmaN}{\sigma_{\tagN}}
\newcommand{\chiF}{\chi_{\tagF}}
\newcommand{\chiN}{\chi_{\tagN}}
\newcommand{\uF}{u_{\tagF}}
\newcommand{\uN}{u_{\tagN}}
\newcommand{\nF}{n_{\tagF}}
\newcommand{\nN}{n_{\tagN}}
\newcommand{\TDN}{T_{\mathrm{DN}}}
\newcommand{\TRN}{T_{\mathrm{RN}}}
\newcommand{\hSN}{\widehat{S}_{\tagN}}      
\newcommand{\EN}{E_{\tagN}}                  
\newcommand{\epsN}{\varepsilon_{\tagN}}      
\newcommand{\That}{\widehat{T}_{\mathrm{RN}}}
\newcommand{\rhohat}{\widehat{\rho}}         
\newcommand{\Miface}{M_{\SigmaI}}
\newcommand{\gint}{g_{\mathrm{int}}}
\newcommand{\gpw}{g^{\mathrm{pw}}}
\newcommand{\norm}[1]{\left\lVert #1 \right\rVert}
\newcommand{\abs}[1]{\left\lvert #1 \right\rvert}
\newcommand{\alphastar}{\alpha^{\star}}
\newcommand{\dt}{\Delta t}
\newcommand{\dn}{\partial_n}
\newcommand{\dnF}{\partial_{\nF}}
\newcommand{\dnN}{\partial_{\nN}}
\newcommand{\NN}{\mathcal{N}}
\newcommand{\Ll}{\mathcal{L}}
\newcommand{\Tinner}{T_{\!\mathrm{in}}}
\newcommand{\Kouter}{K_{\mathrm{out}}}
\newcommand{\figmaybe}[3][\linewidth]{\includegraphics[width=#1]{#2}}

\title{Robin--Neumann Coupling of PINN and FEM Solvers:\\
  A Steklov--Poincar\'e View, with Application to Fluid--Structure
  Interaction with Contact}

\author{%
  Mikel Landajuela \\
  Lawrence Livermore National Laboratory\\
  \texttt{landajuelala1@llnl.gov}
}
\date{}

\begin{document}
\maketitle

\begin{abstract}
Physics-informed neural networks (PINNs) are meshless and carry moving
geometry and topology change through resampling of collocation points;
the finite-element method (FEM) is the workhorse for boundary-fitted
discretisations. Coupling the two across a shared interface promises
the best of both, yet the few existing PINN--FEM schemes are validated
only empirically. We put the coupling on a domain-decomposition
footing. Viewing each solver as a Steklov--Poincar\'e (trace-to-flux)
operator, we transfer the classical Dirichlet--Neumann (DN) divergence
diagnosis and its Robin--Neumann (RN) cure, including a closed-form,
sweep-free interface impedance $\alphastar$, and prove a PINN-specific
contraction theorem: a trained network realises only a perturbed
Steklov operator with a per-step training residual, and RN still
contracts, with no shared-eigenbasis hypothesis, to a floor set by the
achieved training loss. Because a PINN has no stiffness matrix, we
introduce a Fourier-mode interface probe that recovers the network's
resolvable Steklov eigenvalues to within $0.5\%$ (to $10^{-5}$ on the
lowest mode), an order of
magnitude cheaper and orders of magnitude more accurate than
nodal-basis or random-trace estimators at matched budget, and that
doubles as a diagnostic of the network's spectral cap. The theory
predicts measured PINN--FEM contraction rates to within $7\%$ on 1D and
2D Poisson couplings, and a two-slab analogue of the large-added-mass
regime (the interface spectral instability of cardiovascular FSI,
without its inertia) shows RN's per-mode impedance matching winning
decisively where tuned scalar relaxation saturates. We then demonstrate
the framework on a Stokes/rigid-disc problem with Alart--Curnier
contact: the meshless PINN fluid absorbs the topology change at contact
by collocation exclusion alone, no remeshing and no cut cells, and the
static-equilibrium contact reaction matches the submerged weight to
$0.4\%$ under mesh refinement. We quantify the remaining
limitations: the warm-started PINN drifts off the Stokes manifold over
long horizons (tracked by a training-independent residual monitor), and
matched FEM--FEM benchmarks attribute the pre-impact squeeze-film
signatures to PINN under-resolution, leaving the fully resolved dynamic
regime to future work.
\end{abstract}

\section{Introduction}
\label{sec:intro}

\subsection{Why hybrid PINN--FEM coupling, and why FSI with contact}
\label{sec:intro:why}

Physics-informed neural networks~\citep{RaissiPerdikarisKarniadakis2019,
KarniadakisEtAl2021,LuMengMaoKarniadakis2021} have become a popular
mesh-free alternative for forward and inverse PDE problems: a PINN
replaces meshing and quadrature with random sampling and automatic
differentiation, so domain geometry, and changes thereof, enter only
through the collocation set. The finite-element method, by contrast, is
the workhorse of science and engineering
simulation~\citep{HughesCottrellBazilevs2005}; it is unmatched on
body-fitted geometries with strong gradients, and benefits from decades
of solver infrastructure and rigorous convergence theory. A solver that
is part FEM and part PINN is attractive whenever one part of a problem
suits FEM (a thin solid with sharp stress concentrations) and another
suits a PINN (a fluid region whose mesh is expensive to maintain), and
more generally to embed learned constitutive or operator components
inside an established FEM pipeline~\citep{MeethalEtAl2023,
MituschFunkeKuchta2021,TartakovskyEtAl2020}.

\begin{figure}[t]
\centering
\figmaybe[0.95\linewidth]{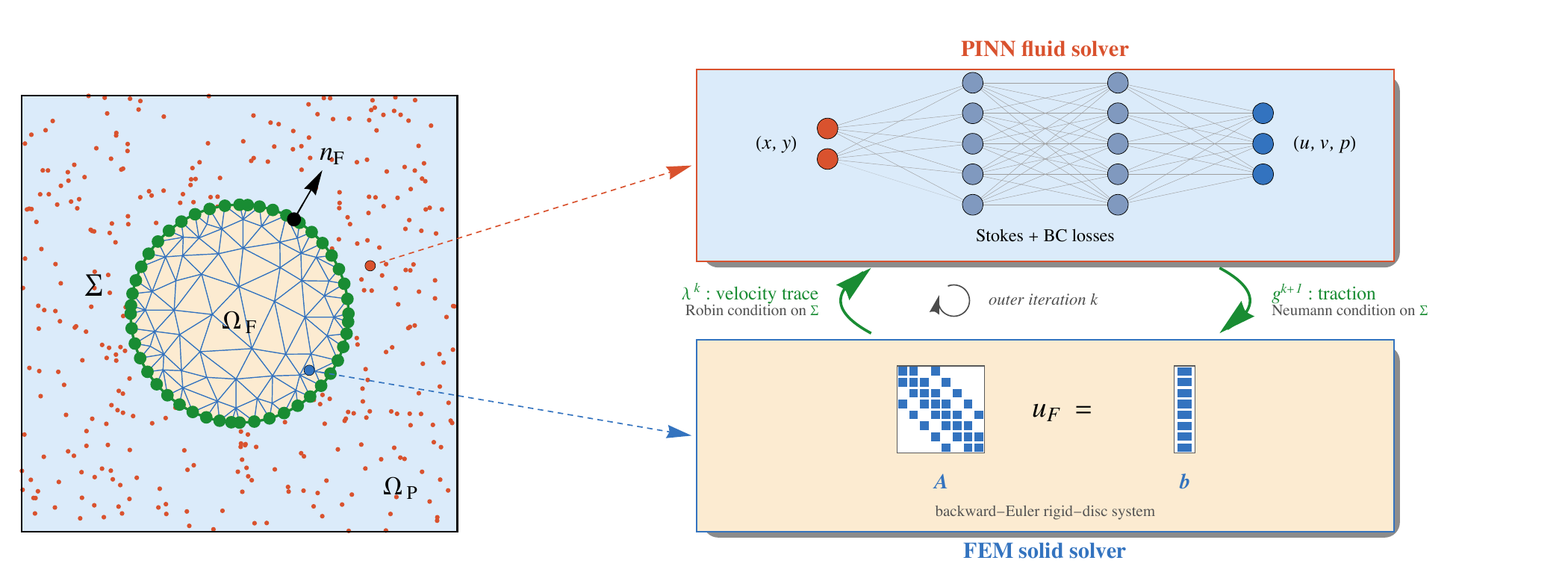}{FSI coupling block diagram.}
\caption{Partitioned PINN--FEM coupling, illustrated on the FSI
application of Section~\ref{sec:fsi}. Left: the mesh-free PINN fluid
$\OmegaN$ (collocation cloud) surrounds the FEM solid $\OmegaF$ (rigid
disc, body-fitted mesh); the two solvers meet only at the shared
interface $\SigmaI$. Right: at outer iteration $k$ the solid sends a
velocity trace $\lambda^k$, imposed on the PINN as a Robin condition,
and the fluid returns a traction $g^{k+1}$, a Neumann load on the solid,
closing a fixed-point loop. When this loop contracts, how fast, and with which
transmission conditions is the subject of
Sections~\ref{sec:theory}--\ref{sec:fsi}.}
\label{fig:fsi-coupling}
\end{figure}

The application that motivates this paper is \textbf{fluid--structure
interaction with contact}. In FSI, an immersed solid (a heart-valve
leaflet, a falling particle) deforms and moves through a fluid; the two
are coupled through stress and velocity matching at their shared
interface. When the solid hits a wall or another solid, classical
body-fitted FEM grids \emph{collapse}: the fluid mesh degenerates as the
gap closes to zero, forcing aggressive remeshing, ALE pseudo-elasticity
tricks, or a cut-cell
formulation~\citep{HronTurek2006,BurmanFreiFernandezGerosa2022}, the last
echoing the unfitted-mesh and Nitsche-XFEM treatments that immerse the
structure to sidestep body-fitted
remeshing~\citep{alauzet2016nitsche,fernandez2015splitting,fernandez2020splitting}.
A PINN
fluid has no such problem: its discretisation is a point cloud, and gap
closure is implemented simply by removing the contact slab from the
collocation set. \emph{Topology change is carried by the sampler.} This
is the concrete advantage we set out to exploit.
Figure~\ref{fig:fsi-coupling} previews the resulting solver: a mesh-free
PINN fluid and a body-fitted FEM solid exchange a velocity trace and an
interface traction around a partitioned fixed-point loop, the iteration
whose convergence this paper analyses.

\subsection{The state of PINN--FEM coupling}
\label{sec:intro:state}

The number of published methods that actually couple a PINN solver to an
FEM solver across a shared interface is surprisingly small.
\citet{SnyderTezaurWentland2023} introduced a Schwarz-style alternating
method for PINN--FEM coupling; the hybrid FEM--NN framework of
\citet{MituschFunkeKuchta2021} replaces a constitutive component of an
FEM model with a neural network; and \citet{MeethalEtAl2023} embed a
PINN inside an FEM iteration. These works establish that the coupling is
feasible and useful, but treat the partitioned iteration largely
empirically: none analyses its convergence rate, its divergence modes,
or the operator-spectral conditions under which a given transmission
scheme contracts. A larger, faster-growing line of work imports
classical domain-decomposition ideas \emph{into} the neural solver
itself: extended and conservative PINNs
\citep{JagtapKarniadakis2020,JagtapKharazmiKarniadakis2020}, finite-basis
PINNs \citep{MoseleyMarkhamNissenMeyer2023}, and their explicit Schwarz
reading \citep{DoleanHeinleinMishraMoseley2024}. But this remains
PINN--PINN: the Steklov--Poincar\'e substructuring view and the
optimal-Robin contraction theory standard for FEM--FEM and partitioned
FSI have not been carried over to the PINN--FEM case. Finally, our PINN
Steklov estimator and the spectral cap it exposes connect to the
well-documented spectral bias of neural solvers
\citep{RahamanEtAl2019,WangTengPerdikaris2021,WangYuPerdikaris2022,
KrishnapriyanEtAl2021}; we test a Fourier-feature
remedy~\citep{TancikEtAl2020} directly and find it does not lift the cap
at our budgets. This paper supplies the missing analysis by carrying the
classical domain-decomposition framework~\citep{QuarteroniValli1991,
QuarteroniValli1999,ToselliWidlund2005}, the optimised Robin
transmission conditions of optimized Schwarz
methods~\citep{GanderMagoulesNataf2002,Gander2006}, from which the
closed-form impedance of Section~\ref{sec:theory} descends, and the
optimal-Robin lessons of
partitioned FSI~\citep{CausinGerbeauNobile2005,BadiaNobileVergara2008,
BurmanFernandez2009,fernandez2013fully,fernandez2015fullydecoupled,
fernandez2016robin,landajuela2016coupling,Landajuela2017} over to the
PINN--FEM setting, where one subdomain
solver is a trained network rather than a second discretisation of the
same kind. We also benchmark the closest prior coupling directly: our
implementation of the overlapping Schwarz scheme
of~\citet{SnyderTezaurWentland2023}, in its weak-Dirichlet PINN--FOM
form, runs head-to-head against DN and RN at matched architecture and
inner budget (Appendix~\ref{app:inner-pinn:e6}). The comparison is
informative in both directions: with a generous overlap the Schwarz
iteration reaches moderate tolerance in fewer outer iterations, while
RN contracts monotonically on every seed, reaches a lower floor, and
needs no overlap, a requirement, not a preference, at the sharp
fluid--solid interface of the FSI application.

\subsection{Contributions}
\label{sec:intro:contributions}

\paragraph{(C1) Steklov--Poincar\'e framework for PINN--FEM coupling,
with a PINN-specific contraction theorem.}
We give a self-contained derivation of the DN and RN partitioned
iterations for the case where one subdomain is solved by FEM and the
other by a PINN (Section~\ref{sec:theory}). Three results are classical
domain-decomposition theory transferred to this setting: the DN
divergence diagnosis, the existence of a contractive band of Robin
impedances with the closed-form sweep-free default
$\alphastar = \sqrt{\sigmaF^{\min}\,\sigmaN^{\max}}$, and energy
boundedness for a co-training variant. The PINN-specific core is
Theorem~\ref{thm:pinn-floor}: a trained network realises only a
\emph{perturbed} Steklov operator with a per-step training residual, and
the RN iteration still contracts, with no shared-eigenbasis hypothesis,
to a floor set by the achieved training loss
(Remark~\ref{rem:loss-floor}), the network's spectral cap itself acting
as the contraction mechanism (Corollary~\ref{cor:cap}).

\paragraph{(C2) A sampling-based estimator of the PINN Steklov
operator.} The FEM operator $\SF$ is a Schur complement and essentially
free; the PINN operator $\SN$ is a trained-network response and much
harder to access. We propose a Fourier-mode interface probe
(Algorithm~\ref{alg:sn-probe}) that recovers continuum eigenvalues to
within $0.5\%$ on the modes the network resolves ($10^{-5}$ on the
lowest) and exposes its
spectral cap on the modes it cannot, and at matched training budget is an
order of magnitude cheaper and orders of magnitude more accurate than
nodal-basis or random-trace estimators (Table~\ref{tab:sn-estimators}).

\paragraph{(C3) The added-mass regime is where RN wins.} On 1D and 2D
Poisson coupling the operator-spectral theory predicts the empirical
PINN--FEM rates to within $7\%$ (Section~\ref{sec:exp:poisson}). The
practical payoff is sharpest in the regime FSI cares about: with a large
added-mass effect (a light structure in a heavy incompressible fluid),
the per-mode DN amplification is strongly mode-dependent, a tuned scalar
relaxation \emph{saturates}, and RN's per-mode impedance matching
contracts ever faster as the added mass grows
(Section~\ref{sec:exp:addedmass}; the two-slab analogue reproduces the
interface \emph{spectral} instability of added mass, not its inertia).

\paragraph{(C4) Application to Stokes--rigid-disc FSI with contact.} We
demonstrate the framework on the geometry that motivated the work
(Section~\ref{sec:fsi}): a rigid disc (FEM) inside a viscous fluid
(PINN, on the outer annulus). The meshless PINN fluid handles topology
change at contact by collocation exclusion (no remeshing, no cut
cells), and the contact reaction matches the submerged weight $\Pi$ to
within $0.4\%$ at static equilibrium under mesh refinement (a
force-balance check of the coupled solver's equilibrium bookkeeping;
the near-wall fluid field itself is under-resolved, which we quantify
separately). We delimit
the methodology's current scope (Section~\ref{sec:fsi:limits}): the
warm-started PINN drifts off the Stokes manifold over long time
horizons, and the pre-impact dynamic signatures are PINN
under-resolution artefacts that matched FEM--FEM benchmarks do not
reproduce.

\section{Steklov--Poincar\'e coupling of PINN and FEM}
\label{sec:theory}

\subsection{Configurations, transmission conditions, and roles}
\label{sec:problem}

\begin{figure}[t]
\centering
\begin{subfigure}[b]{0.32\linewidth}
  \centering
  \figmaybe{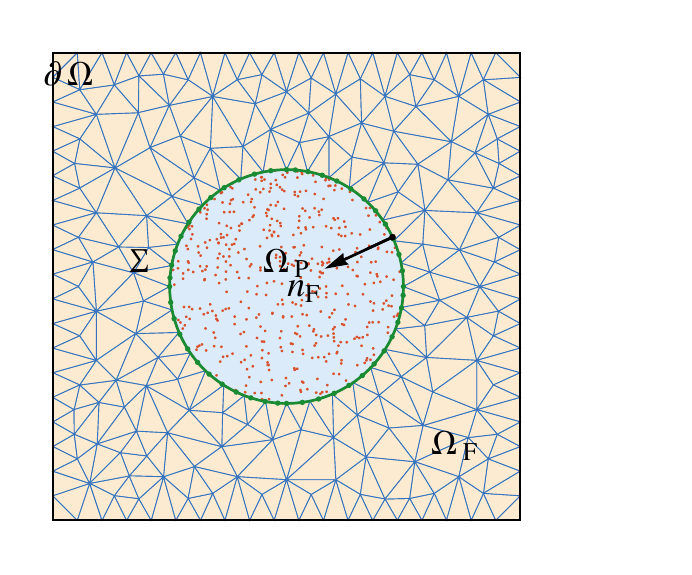}{Natural orientation.}
  \caption{\emph{Natural}: FEM on the outer annulus, PINN on the inner
  disc. Used for the E5/E7/E10 Poisson experiments.}
  \label{fig:domains:natural}
\end{subfigure}\hfill
\begin{subfigure}[b]{0.32\linewidth}
  \centering
  \figmaybe{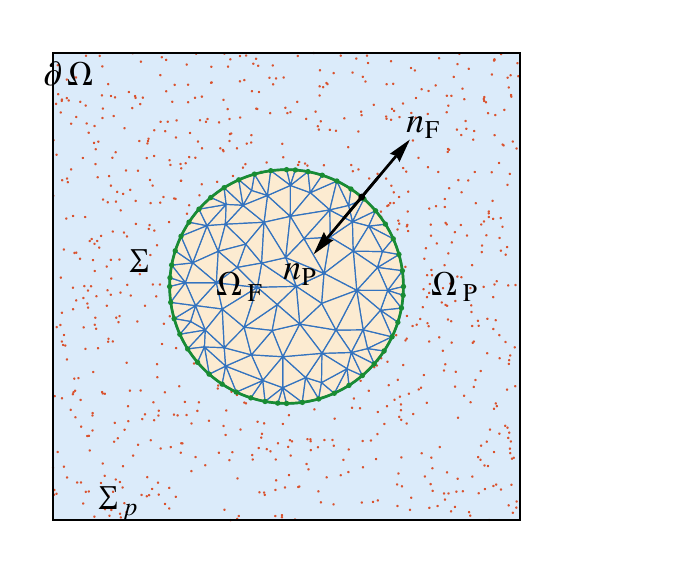}{FSI-aligned orientation.}
  \caption{\emph{FSI-aligned}: FEM on the inner disc (rigid solid), PINN
  on the outer annulus (viscous fluid).}
  \label{fig:domains:fsi-aligned}
\end{subfigure}\hfill
\begin{subfigure}[b]{0.32\linewidth}
  \centering
  \figmaybe{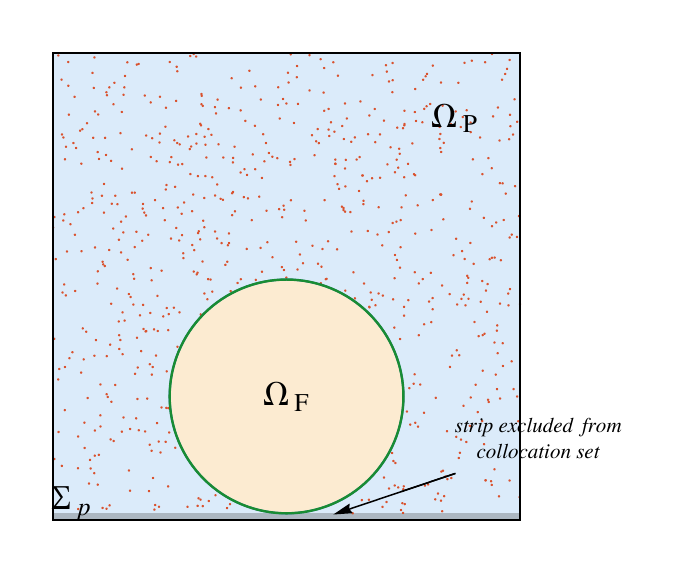}{Contact at bottom wall, collocation strip excluded.}
  \caption{\emph{Contact}: the disc has descended to $y_c = R +
  \varepsilon_g$ and the strip $\{y < \varepsilon_g\}$ is excluded from
  the collocation set.}
  \label{fig:domains:contact}
\end{subfigure}
\caption{Two role assignments for partitioned PINN--FEM coupling (a, b)
and the contact-engaged configuration that motivates the FSI
application (c). The PINN subdomain is mesh-free in every panel; the
shared interface $\SigmaI$ is green.}
\label{fig:domains}
\end{figure}

Let $\Omega \subset \R^d$ be a bounded Lipschitz domain decomposed
non-overlappingly into a \emph{FEM subdomain} $\OmegaF$ and a
\emph{PINN subdomain} $\OmegaN$, with interface
$\SigmaI = \overline{\OmegaF} \cap \overline{\OmegaN}$
(Figure~\ref{fig:domains}). As a scalar model problem we use the
Poisson equation
\begin{equation}
-\Delta u = f \quad \text{in } \Omega, \qquad u = 0 \quad \text{on } \partial \Omega ,
\label{eq:poisson}
\end{equation}
on a 1D split segment and a 2D disc-in-square, each in two orientations:
FEM-on-disc / PINN-on-annulus (the \emph{natural} orientation) and the
\emph{FSI-aligned} converse. The target application
(Section~\ref{sec:fsi}) is a rigid disc moving through a viscous fluid,
with the FEM subdomain on the disc and the PINN subdomain on the moving
fluid annulus.

A solution of~\eqref{eq:poisson} restricted to $\OmegaF$ and $\OmegaN$
satisfies, at $\SigmaI$, two transmission conditions:
$\uF = \uN$ (trace continuity) and $\dnF \uF = \dnF \uN$ (flux
continuity), where $\nF$ is the unit normal pointing out of $\OmegaF$
(so $\nN = -\nF$). A \emph{partitioned} solver enforces one condition
exactly on each side and iterates the other. \textbf{Dirichlet--Neumann
(DN):} FEM imposes $\uF|_\SigmaI = \lambda^k$, PINN imposes
$\dnN \uN|_\SigmaI = - \dnF \uF^k|_\SigmaI$, update
$\lambda^{k+1} = \uN^k|_\SigmaI$. \textbf{Robin--Neumann (RN):} FEM
imposes the Robin condition $\dnF \uF + \alpha\, \uF = \alpha\,\lambda^k
+ g^k$ for impedance $\alpha > 0$; PINN imposes the same Neumann
condition; update $(\lambda^{k+1}, g^{k+1})$ from the PINN trace and
flux. RN is the partitioned-FSI fix for the added-mass instability of
DN~\citep{CausinGerbeauNobile2005,BadiaNobileVergara2008,fernandez2015fullydecoupled,fernandez2016robin};
we show the
same fix works for PINN--FEM for the same spectral reason.

The DN/RN definitions attach the Dirichlet (or Robin) datum to $\OmegaF$
and the Neumann datum to $\OmegaN$, but the interface
equation~\eqref{eq:steklov-interface} is symmetric in the two sides, so
the assignment is a modelling choice. We use \emph{Dirichlet-on-FEM,
Neumann-on-PINN} for Poisson coupling (the cheap-FEM-Dirichlet choice),
and \emph{Robin-on-PINN, Neumann-on-FEM} for FSI, forced by the physics:
the kinematic condition $u_{\mathrm{fluid}}|_\SigmaI = \dot
d_{\mathrm{solid}}$ is a velocity trace on the fluid (PINN) side, while
the integrated traction $F_y = \int_{\SigmaI}(\sigma_f n)\cdot e_y\,ds$
is a Neumann forcing of the solid (FEM) ODE. The added-mass instability
lives on the fluid side in both formulations, so Robin-on-PINN remains
the contraction-restoring choice; all spectral results below apply after
the $\SF \leftrightarrow \SN$ relabelling, with $\alphastar$ invariant.
We confirm empirically that the RN contraction itself survives this
relabelling (Section~\ref{sec:exp:2d}). The full per-piece boundary
conditions for both orientations are collected in
Appendix~\ref{app:bcs}, the shared-interface-DOF and sign
conventions in Appendix~\ref{app:conventions}, and a symbol glossary in
Appendix~\ref{app:glossary}.

\subsection{The interface equation and the DN/RN maps}
\label{sec:theory:interface}

For each subdomain $\bullet \in \{F, N\}$ write the subdomain solution
at interface trace $\lambda$ as the sum of a homogeneous (zero-forcing)
and an inhomogeneous (forcing-only, zero-trace) part. The
\emph{Steklov--Poincar\'e (Dirichlet-to-Neumann) operator} is the
linear part of the trace-to-flux map
\citep{QuarteroniValli1991,QuarteroniValli1999},
\begin{equation}
S_\bullet \, \lambda := \dnF u_\bullet^{\mathrm{hom}}(\lambda)\big|_\SigmaI ,
\qquad \chi_\bullet := -\dnF u_\bullet^{\mathrm{inh}}\big|_\SigmaI ,
\label{eq:steklov-def}
\end{equation}
so the flux at trace $\lambda$ is $\dnF u_\bullet(\lambda) = S_\bullet
\lambda - \chi_\bullet$, with a sign flip on the PINN side
($\nN=-\nF$):
\begin{equation}
\dnF \uF(\lambda^\star) = \SF \lambda^\star - \chiF , \qquad
\dnF \uN(\lambda^\star) = -\SN \lambda^\star + \chiN .
\label{eq:fluxes}
\end{equation}
Imposing flux continuity gives the \emph{Steklov--Poincar\'e interface
equation} on $\SigmaI$ alone,
\begin{equation}
\boxed{ \; (\SF + \SN)\, \lambda^\star = \chiF + \chiN \quad \text{on } \SigmaI \; }
\label{eq:steklov-interface}
\end{equation}
Each subdomain response in~\eqref{eq:fluxes} is the \emph{exact} PDE
extension of its interface datum, so the bulk unknowns are condensed
into $\SF,\SN$, which is why~\eqref{eq:steklov-interface} lives on
$\SigmaI$ alone; the two schemes' update rules then reduce to a linear
recurrence on the trace error $e^k = \lambda^k - \lambda^\star$
(derived step by step in Appendix~\ref{app:dn-rn-derivation}):
\begin{equation}
e^{k+1} = -\, \SN^{-1} \SF \, e^k =: \TDN \, e^k ,
\qquad
e^{k+1} = (\alpha I - \SN)(\alpha I + \SF)^{-1} \, e^k =: \TRN(\alpha) \, e^k .
\label{eq:dn-rn-map}
\end{equation}
In practice each outer iteration is one cheap FEM solve plus one PINN
training run, communicating only through $\lambda$ (and, for RN, the
flux $g$); the two modes differ in exactly one step, the FEM
sub-problem's boundary condition, which we state together as
Algorithm~\ref{alg:dn-rn}. For the FSI assignment the roles of ``FEM
solve'' and ``PINN train'' swap, but the algorithm is otherwise
unchanged after the vector generalisation of
Section~\ref{sec:fsi:model}.

\begin{algorithm}[t]
\caption{\textsc{Partitioned PINN--FEM iteration} (Dirichlet/Robin-on-FEM,
Neumann-on-PINN). Mode \texttt{m}=\texttt{DN} gives Dirichlet--Neumann;
\texttt{m}=\texttt{RN} gives Robin--Neumann.}
\label{alg:dn-rn}
\begin{algorithmic}[1]
\Require Mode $\texttt{m} \in \{\texttt{DN}, \texttt{RN}\}$; initial trace
$\lambda^0$ (and, for \texttt{RN}, flux $g^0$ and impedance $\alpha>0$,
sweep-free default $\alphastar=\sqrt{\sigmaF^{\min}\sigmaN^{\max}}$ of
Theorem~\ref{thm:rn}); under-relaxation
$\omega\in(0,1]$ (default $1$ for Poisson, ${\approx}0.3$ for vector
Stokes); budget $\Kouter$, tolerance $\mathrm{tol}$.
\For{$k = 0, 1, \ldots, \Kouter - 1$}
  \State \textbf{FEM solve.} $-\Delta u=f$ on $\OmegaF$ with either
  $\uF^k|_\SigmaI=\lambda^k$ (\texttt{DN}) or the Robin condition
  $\dnF\uF+\alpha\uF=\alpha\lambda^k+g^k$ on $\SigmaI$ (\texttt{RN}).
  \State \textbf{Read FEM trace/flux.} $\lambda^k_F\gets\uF^k|_\SigmaI$;
  $g^k_F\gets \Miface^{-1}(K\uF^k-f)[\SigmaI]$ via~\eqref{eq:units-fix}.
  \State \textbf{PINN (Neumann) train.} minimise the PINN loss enforcing
  $-\Delta\uN=f$ on $\OmegaN$ and $\dnN\uN|_\SigmaI=-g^k_F$.
  \State \textbf{Read PINN trace.} $\lambda^k_N\gets\uN^k|_\SigmaI$.
  \State \textbf{Update.} \texttt{DN}: $\lambda^{k+1}\gets\lambda^k_N$;
  \texttt{RN}: $(\lambda,g)^{k+1}\gets(1-\omega)(\lambda,g)^k+\omega(\lambda^k_N,g^k_F)$.
  \If{$\norm{\lambda^k_N-\lambda^k_F}<\mathrm{tol}$} \Return \EndIf
\EndFor
\end{algorithmic}
\end{algorithm}

\subsection{Convergence theory: classical transfer and a PINN-specific
contraction theorem}
\label{sec:theory:theorems}

The classical results share one working hypothesis: that $\SF$ and
$\SN$ admit a \emph{common eigenbasis} on $\SigmaI$ with strictly
positive eigenvalues, collapsing the maps of~\eqref{eq:dn-rn-map} onto
independent scalar recurrences per interface mode. On the continuum
this holds whenever $\SigmaI$ carries a natural shared mode basis
\emph{and the remainder of each subdomain's boundary respects it}: a
single point in 1D; concentric circles in 2D, where both operators are
angular convolutions diagonalised by $\{\cos(k\theta),\sin(k\theta)\}$,
with $\sigmaN(k)=k/R$ for a disc and the closed
form~\eqref{eq:sn-annulus} for an annulus. Two distinct effects degrade
the hypothesis, and we keep them separate: a \emph{geometric} one (a
non-polar outer boundary, such as the square in our 2D experiments,
destroys the shared eigenbasis \emph{already in the continuum}; a clean
FEM--FEM benchmark on that geometry indeed diverges with no contractive
$\alpha$, Section~\ref{sec:exp:2d}) and a \emph{discretisation} one
(operators on different meshes co-diagonalise only approximately even
on polar geometry). The closed-form rates below are therefore the
continuum, polar-boundary limit; the genuinely-discrete case is treated
in Appendix~\ref{app:proofs}.

\begin{theorem}[DN divergence diagnosis;
compare~\citealp{QuarteroniValli1999}]
\label{thm:dn}
If $\SF$ and $\SN$ are simultaneously diagonalisable with positive
eigenvalues $\{\sigma_{\tagF,k}\},\{\sigma_{\tagN,k}\}$ on a common
basis, then $\rho(\TDN) = \max_k \sigma_{\tagF,k}/\sigma_{\tagN,k}$, and
DN diverges whenever this exceeds $1$.
\end{theorem}

\begin{theorem}[RN contractive band and the optimised-Robin impedance;
compare~\citealp{QuarteroniValli1999,Gander2006}]
\label{thm:rn}
Under the same assumption, for every $\alpha>0$,
$\rho(\TRN(\alpha)) = \max_k \abs{\alpha-\sigma_{\tagN,k}}/(\alpha+\sigma_{\tagF,k})$.
Consequently RN contracts if and only if
$\alpha > \alpha_{\min} := \tfrac12\max_k(\sigma_{\tagN,k}-\sigma_{\tagF,k})_+$:
unlike DN, RN admits a non-empty contractive band
$(\alpha_{\min},\infty)$ for \emph{every} spectrum. When the binding
modes are spectrally aligned ($\sigma_{\tagF,k}=\sigma_{\tagN,k}=\sigma_k$,
the symmetric case of optimized Schwarz methods, with the combined
spectrum spanning $[\sigmaF^{\min},\sigmaN^{\max}]$), the minimax over
$\alpha$ is attained by equioscillation at
$\alphastar=\sqrt{\sigmaF^{\min}\sigmaN^{\max}}$, with
$\rho(\TRN(\alphastar)) = (\sqrt{\kappa}-1)/(\sqrt{\kappa}+1)<1$,
$\kappa=\sigmaN^{\max}/\sigmaF^{\min}$. For mismatched spectra
$\alphastar$ remains a closed-form, sweep-free default inside the band,
but it need not be the iteration optimum
(Appendix~\ref{app:proofs:rn}).
\end{theorem}

\begin{theorem}[Energy boundedness for RN co-training: decrease up to a
bounded forcing term]
\label{thm:cotraining}
Let $E^k = \tfrac12 \norm{\nabla\uF^k}^2_{L^2(\OmegaF)} +
\tfrac12\Ll_{\tagN}(\theta^k) +
\tfrac{\alpha}{2}\norm{\uF^k-\uN^k}^2_{L^2(\SigmaI)}$, where
$\Ll_{\tagN}$ is the PINN residual--Dirichlet loss. For sufficiently
small PINN learning rate and any $\alpha$ in the contractive band of
Theorem~\ref{thm:rn} (in particular $\alphastar$ in the aligned case),
$E^{k+1}\le E^k$ up to a bounded forcing term; the energy is therefore
bounded, decreasing whenever the trace mismatch dominates the forcing
offset.
\end{theorem}

Full proofs, including a complete proof of
Theorem~\ref{thm:cotraining} and the discrete consequences when the
eigenbases of $\SF$ and $\SN$ differ, are in
Appendix~\ref{app:proofs}.

Theorems~\ref{thm:dn}--\ref{thm:cotraining} treat the network as a
\emph{black-box exact} operator $\SN$ acting in a shared eigenbasis: no
approximation error, no optimisation error, and no training appears, and
the shared-basis hypothesis fails already in the continuum for the square
outer boundary of our 2D runs (the first effect above). We therefore add
a genuinely PINN-specific statement. A trained network does not realise
$\SN$ but a \emph{perturbed} Steklov operator $\hSN=\SN+\EN$, where $\EN$
is the measured spectral cap (it zeroes the high interface modes the
network cannot represent, Remark~\ref{rem:spectral-cap}), together with a
per-step residual from finite training. The cap is \emph{data-linear} and
perturbs the iteration operator (it moves the rate); the residual is
additive and sets a floor. The result uses operator norms only---no
shared eigenbasis---so it covers the square-boundary geometry that
Theorems~\ref{thm:dn}--\ref{thm:rn} miss.

\begin{theorem}[Inexact-PINN RN contraction to an $O(\epsN)$ floor; no
shared eigenbasis required]
\label{thm:pinn-floor}
Let $\SF$ be the exact FEM Steklov operator (SPD, smallest eigenvalue
$\sigmaF^{\min}>0$), and suppose at each outer step the PINN realises a
perturbed Steklov operator $\hSN=\SN+\EN$ on the interface and returns a
trace with per-step additive residual $\tilde r^k$ obeying
$\norm{\tilde r^k}\le C_\alpha\,\epsN$, where
$C_\alpha=\norm{(\alpha I+\SF)^{-1}}=1/(\alpha+\sigmaF^{\min})$ and
$\epsN$ bounds the PINN approximation-plus-optimisation error (uniform in
$k$ because Algorithm~\ref{alg:dn-rn} retrains to a fixed loss tolerance
each step). Then the RN trace error satisfies the perturbed linear
iteration $e^{k+1}=\That(\alpha)\,e^k+\tilde r^k$ with
$\That(\alpha)=\TRN(\alpha)-\EN(\alpha I+\SF)^{-1}$, whose operator norm
obeys the triangle bound
\begin{equation}
  \rhohat(\alpha):=\norm{\That(\alpha)}
  \;\le\; \norm{\TRN(\alpha)} + \frac{\norm{\EN}}{\alpha+\sigmaF^{\min}}.
  \label{eq:rhohat}
\end{equation}
If $\rhohat(\alpha)<1$, then for all $k$
\begin{equation}
  \norm{e^k}\;\le\;\rhohat(\alpha)^k\,\norm{e^0}
  \;+\;\frac{C_\alpha\,\epsN}{1-\rhohat(\alpha)},
  \qquad
  \limsup_{k\to\infty}\norm{e^k}\;\le\;\frac{C_\alpha\,\epsN}{1-\rhohat(\alpha)}
  \;=\;O(\epsN).
  \label{eq:pinn-floor}
\end{equation}
As $\epsN\to0$ and $\EN\to0$ the floor vanishes and
$\rhohat\to\rho(\TRN)$, recovering Theorem~\ref{thm:rn}.
\end{theorem}

The condition $\rhohat(\alpha)<1$ is an assumption, not a consequence:
the cap inflates the operator norm of a non-normal $\That$, and
\eqref{eq:rhohat} can exceed $1$ even when $\rho(\TRN)<1$. It holds
precisely when the cap concentrates $\EN$ in the modes that
$(\alpha I+\SF)^{-1}$ damps, which is the empirically observed regime and
the content of the next corollary.

\begin{corollary}[The spectral cap can open a contractive band]
\label{cor:cap}
Split the interface space $V=V_{\mathrm{res}}\oplus V_{\mathrm{cap}}$ into
resolved (low) and capped (high) bands. If $\EN$ acts as the cap on
$V_{\mathrm{cap}}$ (so $\hSN\approx0$ there) and the FEM Steklov
eigenvalues grow with mode, then on $V_{\mathrm{cap}}$ the iteration
reduces to $\That\approx\alpha(\alpha I+\SF)^{-1}$, of norm
$\norm{\That|_{V_{\mathrm{cap}}}}\le\alpha/(\alpha+\sigmaF^{\mathrm{high}})<1$.
The capped band thus contracts even where the exact $\TRN$ mixes
incompatible bases and need not contract: the cap is the contraction
mechanism, not merely a diagnostic.
\end{corollary}

\begin{remark}[A loss-controlled floor]
\label{rem:loss-floor}
The residual $\epsN$ is not a free constant: it is controlled by the
\emph{achieved} PINN training loss. For the Neumann/Robin subproblem,
coercivity of the elliptic operator and a trace inequality give the
a-posteriori bound
$\epsN=\norm{\lambda_{\tagN}^k-\lambda_{\tagN}^{\mathrm{exact}}}_{H^{1/2}(\SigmaI)}
\le C_{\mathrm{tr}}C_{\mathrm{stab}}\,\Ll_{\tagN}(\theta^k)^{1/2}$
(Appendix~\ref{app:proofs:pinn-floor}), so the floor scales as
$O(\Ll_{\tagN}^{1/2})$: the plateau height is set by how well the PINN was
trained. A training-budget sweep on the 2D natural orientation confirms
this falsifiable prediction: the measured RN plateau tracks the achieved
$\Ll_{\tagN}^{1/2}$ with log--log slope $1.13$ ($R^2=0.89$ over $27$ runs
spanning a $400\times$ budget range, three seeds;
Appendix~\ref{app:E25}).
Bounding the optimisation part of $\Ll_{\tagN}$ further by an NTK
gradient-descent rate is possible but fragile and unnecessary here; we
leave it to future work.
\end{remark}

\subsection{Estimating \texorpdfstring{$\SN$}{SN} for a PINN}
\label{sec:theory:sn-estimator}

For a PINN there is no stiffness matrix to take a Schur complement of.
A \emph{basis-vector probe} (train with $\uN|_\SigmaI = e_j$, read the
flux at all nodes to recover column $j$ of $\bSN$) is expensive, costing
one training per column, and the delta-spike traces are hard for an MLP,
contaminating the matrix with $4$--$100\%$ asymmetry. The contamination
is diagnostic precisely \emph{because} the true operator is symmetric:
the continuous Steklov--Poincar\'e (Dirichlet-to-Neumann) operator of a
self-adjoint elliptic problem is itself self-adjoint and positive on the
interface, an immediate consequence of Green's identity, so a
well-resolved discrete $\bSN$ must be symmetric up to quadrature error
\citep{QuarteroniValli1999,ToselliWidlund2005}. Measured asymmetry of
tens of percent therefore signals that the nodal traces excite modes the
network cannot represent, not a genuine property of $\bSN$. We instead
use a
\emph{Fourier-mode probe} (Algorithm~\ref{alg:sn-probe}): on a circular
interface train PINNs with $\uN|_\SigmaI\in\{\cos(k\theta),
\sin(k\theta)\}$, read the radial flux, and project back onto the same
mode. For a disc the continuum eigenvalues are $\sigmaN(k)=k/R$,
recovered to within $0.5\%$ for the resolved modes
(Table~\ref{tab:sn-estimators}); for an annulus with
outer Dirichlet at radius $r$ and interface at radius $R$,
\begin{equation}
\sigmaN(k) = -\frac{k}{r}\,\frac{\rho^{2k}+1}{\rho^{2k}-1}, \quad \rho = R/r,
\qquad
\sigma_0 = -\frac{1}{r \, \ln \rho}
\label{eq:sn-annulus}
\end{equation}
is matched to a few percent in the resolvable band. Because the Steklov
operator is self-adjoint on a smooth interface, the projection is the
per-mode eigenvalue; where the PINN cannot resolve a high mode the
measured $\hat\sigmaN(k)$ collapses toward $0$ instead of following the
continuum $k/R$, so the probe doubles as a spectral-cap diagnostic
(Remark~\ref{rem:spectral-cap}). The probe is also far more
\emph{sample-efficient}, because each training targets a single
resolvable mode rather than a delta spike or a random trace: on the
disc at $n_{\mathrm{iface}}=32$, \emph{five} trainings pin $\sigmaN(1)$
to $10^{-5}$, whereas the basis-vector and random-trace estimators need
${\sim}n_{\mathrm{iface}}$ trainings to even form $\bSN$ and are
${\sim}10^3\times$ less accurate on the low modes
(Table~\ref{tab:sn-estimators}; implementation in
Appendix~\ref{app:sn-probe}).

\begin{table}[t]
\centering
\caption{Estimating the disc PINN Steklov operator
($n_{\mathrm{iface}}=32$, matched $\Tinner=1500$ Adam steps): relative
error of $\sigmaN(k)$ against the continuum $k/R$. The structured
Fourier probe is orders of magnitude more accurate with an order of
magnitude fewer trainings.}
\label{tab:sn-estimators}
\small
\begin{tabular}{@{}lrrrr@{}}
\toprule
estimator & PINN trainings & $k=1$ & $k=2$ & $k=3$ \\
\midrule
Fourier (proposed)    & $5$  & $1{\times}10^{-5}$   & $1.2{\times}10^{-3}$ & $4.5{\times}10^{-3}$ \\
basis-vector          & $33$ & $1.2{\times}10^{-2}$ & $4.9{\times}10^{-2}$ & $7.2{\times}10^{-1}$ \\
random-LSQ ($N{=}32$) & $33$ & $1.3{\times}10^{-2}$ & $3.7{\times}10^{-3}$ & $2.8{\times}10^{-1}$ \\
random-LSQ ($N{=}64$) & $65$ & $6.0{\times}10^{-3}$ & $9.9{\times}10^{-3}$ & $1.1{\times}10^{-1}$ \\
\bottomrule
\end{tabular}
\end{table}

One discrete \emph{units pitfall} deserves a flag: the FEM side outputs
an \emph{integrated} reaction (units of force), whereas the PINN's
Neumann loss expects a \emph{pointwise} flux; the conversion
$\gpw=\Miface^{-1}\gint$ through the interface mass
matrix~\eqref{eq:units-fix} is essential, and omitting it inverts both
the DN and RN headline conclusions in 2D
(Remark~\ref{rem:units}, Appendix~\ref{app:conventions}). After the
fix, PINN--FEM tracks FEM--FEM to within $7\%$ on all rates we
measured.

\begin{remark}[The spectral cap and why it matters]
\label{rem:spectral-cap}
The collapse $\hat\sigmaN(k)\!\to\!0$ above some cap $k_\star$ is not an
artefact of the probe but a manifestation of the well-documented
\emph{spectral bias} of neural networks: within a finite training
budget, gradient descent on an MLP cannot represent a high-wavenumber
interface trace
\citep{RahamanEtAl2019,WangYuPerdikaris2022,KrishnapriyanEtAl2021,
WangTengPerdikaris2021}; here the annulus PINN of the FSI-aligned
orientation reaches only $\abs{\sigmaN}_{\max}\!\approx\!6$ against the
continuum's unbounded spectrum~\eqref{eq:sn-annulus}
(Section~\ref{sec:exp:2d}, Figure~\ref{fig:e9-fourier-swap}a). The cap has two
opposite consequences for coupling, both of which we observe: when the
dangerous modes are the high ones, as in the eigenbasis-mismatched
FSI-aligned geometry, the cap is benign and even \emph{rescues} RN by
zeroing $\bSN$'s unresolved high modes (Section~\ref{sec:exp:2d},
Corollary~\ref{cor:cap}); when
accuracy in the high band is required, it is a hard accuracy ceiling.
Fourier-feature embeddings
\citep{TancikEtAl2020,WangTengPerdikaris2021}, the standard remedy, do
not lift it at fixed width and depth (Figure~\ref{fig:e11-ablation}),
and neither does any gradient-descent lever we tried---second-order
optimisation, adaptive collocation, or exactly divergence-free
conditioning (Appendix~\ref{app:E22}). Yet the cap is an
\emph{optimisation} effect, not a
hard representational wall: replacing the gradient-trained network with a
random-feature model fitted by a single least-squares solve removes the
frequency dependence entirely and recovers the interface operator to
finite-element accuracy (Appendix~\ref{app:E23}; Figure~\ref{fig:e9-cap-recovery}
demonstrates the same recovery, mode for mode, on the scalar annulus
probe of Figure~\ref{fig:e9-fourier-swap}). The practical levers,
then, are a non-gradient (least-squares) inner solve, or---staying with
gradient descent---a longer inner budget on a \emph{small} network, which
recedes the cap only to a plateau well short of the least-squares solve;
enlarging the network instead degrades it (the E24 scale-up of
Appendix~\ref{app:E22}).
\end{remark}

\section{Numerical experiments: Poisson coupling}
\label{sec:exp:poisson}

\paragraph{1D split-segment.} On $\Omega=(0,1)$ split at $x=L$ (P1 FEM
left, MLP PINN right, manufactured solution), the empirical DN rate
tracks the theory line $\sigmaF/\abs{\sigmaN}=(1-L)/L$ to within $2\%$
on $7$ of $11$ sweep points and recovers the divergence sign on all
$11$; at $L=0.3$ DN diverges within four iterations while RN contracts
to the inner-PINN noise floor, the coupling energy growing
${\sim}30\times$ under DN and descending to a plateau under RN, the
signature of Theorem~\ref{thm:cotraining}. The full sweep, traces, and
energy curves are in Appendix~\ref{app:1d}
(Figure~\ref{fig:1d}).

\subsection{2D Poisson: disc-in-square, both orientations}
\label{sec:exp:2d}

Mesh $n_{\mathrm{iface}}=32$, $h=0.06$, manufactured
$u^\star=\sin(\pi x)\sin(\pi y)$. In the \emph{natural} orientation (FEM
annulus, PINN disc; E5) DN diverges at empirical rate $1.442$, matching
the FEM--FEM analogue's $1.452$ to within $0.7\%$, and RN at the
spectrally-derived $\alphastar$ contracts the field error
$0.025\to3\times10^{-3}$ over ten iterations; the
operator-spectral predictions agree with the empirical rates to within
$7\%$ (details, the co-training variant, and the head-to-head with the
overlapping-Schwarz coupling of \citet{SnyderTezaurWentland2023} in
Appendix~\ref{app:inner-pinn}). The closed-form
$\alphastar$ pays a quantifiable premium over a numerical
$\alpha$-search, rate $0.78$ at $\alphastar=4.16$ vs.\ $0.668$ at the
searched $\alpha=10.31$ (Appendix~\ref{app:inner-pinn:e10}), the
mismatched-spectra effect anticipated in Theorem~\ref{thm:rn}. The
estimate of $\sigmaN^{\max}$ feeding the closed form matters less than
one might fear: $\alphastar=4.16$ above uses the raw basis-probe
spectrum ($\abs{\sigmaN}^{\max}=67.4$), while the cleaned estimate
($\abs{\sigmaN}^{\max}\approx8.05$) gives $\alphastar=1.44$
(Appendix~\ref{app:inner-pinn}); the bowl is broad enough that both
operate comparably (the E10 sweep measures rate $0.75$ at its
$\alpha=1.41$ grid point).

\begin{figure}[t]
\centering
\figmaybe[0.95\linewidth]{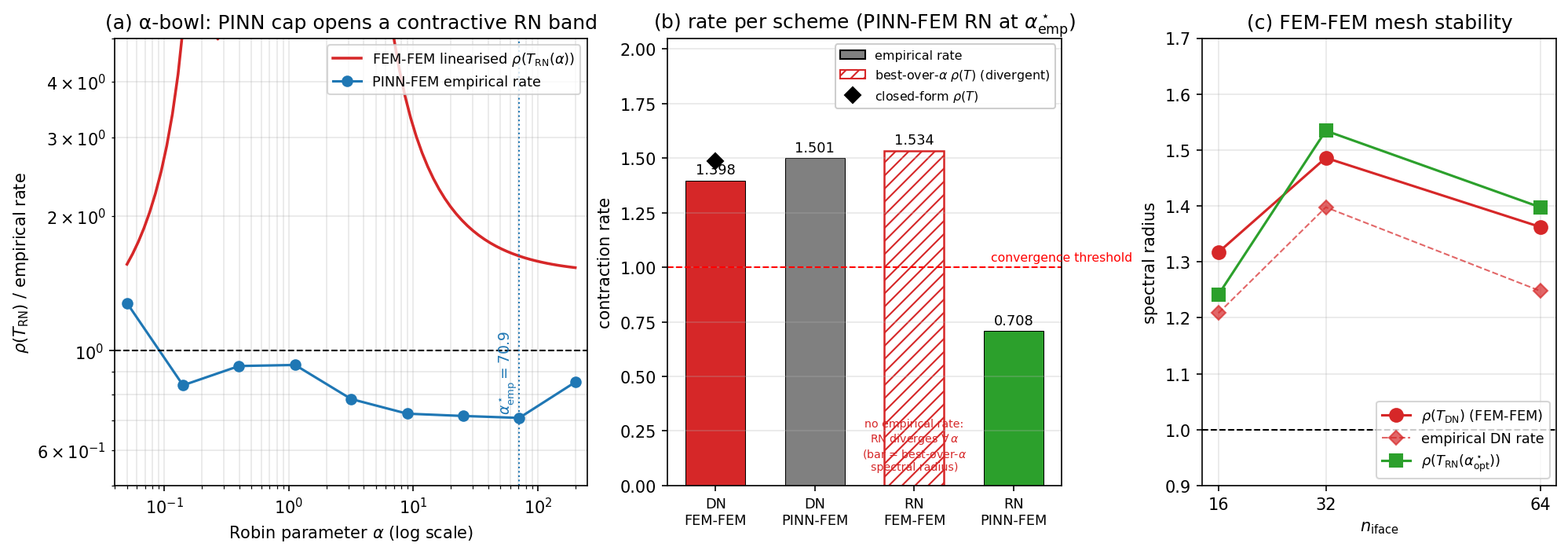}{E9 three-layer comparison: PINN-FEM has a contractive RN band (a) that FEM-FEM lacks; PINN cap costs $+7\%$ on DN, gains a contractive RN (b); FEM-FEM divergence is mesh-stable (c).}
\caption{Three-layer comparison for the FSI-aligned orientation. (a)
$\rho(\TRN(\alpha))$ for the FEM--FEM analogue (red) vs.\ empirical
PINN--FEM rates (blue): the PINN's spectral cap opens a contractive band
FEM--FEM lacks. (b) One contraction-rate bar per scheme at its best
operating point (solid: empirical; black diamond: closed-form
$\rho(\TDN)$): DN diverges on both sides ($+7\%$ PINN penalty); RN
contracts only on the PINN--FEM side. The RN FEM--FEM column has no
empirical bar \emph{by design}, RN diverging there for every $\alpha$;
its hatched bar is the best-over-$\alpha$ linearised spectral radius
${\approx}1.53>1$. (c) FEM--FEM at $n_{\mathrm{iface}}\in\{16,32,64\}$:
$\rho(\TDN)>1$ at every level, a continuum (geometric) property, not a
discretisation artefact.}
\label{fig:e9-three-layer}
\end{figure}

The \emph{FSI-aligned} orientation (FEM disc, PINN annulus; E9) is more
nuanced and becomes legible only once three layers are compared
(Figure~\ref{fig:e9-three-layer}). The continuum circular annulus has
$\abs{\sigmaF(k)/\sigmaN(k)}<1$ per mode, naively predicting DN
convergence. But the \emph{square} outer boundary destroys polar
separation: $\SF$ and $\SN$ no longer share an eigenbasis, $\TDN$ mixes
modes across incompatible bases, and a clean FEM--FEM benchmark on the
same geometry diverges ($\rho(\TDN)=1.486$, empirical $1.398$) with
\emph{no} contractive $\alpha\in[0.05,200]$, a geometric property
stable under mesh refinement ($\rho(\TDN)\in[1.32,1.49]$ at
$n_{\mathrm{iface}}\in\{16,32,64\}$). Replacing the FE annulus by the
PINN is then \emph{dual-edged}: the PINN's spectral cap (Fourier-probe
$\abs{\sigmaN}_{\max}\approx6$ vs.\ the clean FE $\approx225$) costs a
$+7\%$ DN penalty but \emph{rescues RN}, opening a broad contractive
band $\alpha\in[0.14,200]$ (empirical optimum $0.708$ at
$\alpha\approx71$). Mechanistically this is Corollary~\ref{cor:cap} in
action: the cap zeroes $\SN$'s unresolved high modes in the numerator
of $\TRN(\alpha)$ exactly where the eigenbasis mismatch is worst, and
the growing disc Steklov magnitudes damp them in the denominator.
Re-running the coupled iteration with the
\emph{FSI transmission roles} (Robin-on-PINN/Neumann-on-FEM) leaves the
contraction intact: on its own single-seed $\alpha$-grid (a separate
sweep over $[0.2,200]$, hence the slightly different band edge and
optimum from the E9 search just quoted) the swap experiment finds both
placements contractive with near-identical best rates ($0.72$ vs.\
$0.74$), confirming the role-invariance claimed in
Section~\ref{sec:problem} (Appendix~\ref{app:inner-pinn}). Across both
orientations the Steklov--Poincar\'e framework predicts every empirical
PINN--FEM rate we measured to within $7\%$.

\subsection{The large-added-mass regime: where RN wins}
\label{sec:exp:addedmass}

A partitioned-iteration practitioner will ask whether RN at $\alphastar$
beats simply \emph{under-relaxing} the divergent DN map. On the
disc-in-square the answer is no, because it is a \emph{low-added-mass}
problem: $\sigmaN(k)=k/R$ and $\sigmaF(k)$ both grow linearly, so the
per-mode amplification $r_k=\sigmaF(k)/\sigmaN(k)$ is essentially
constant in $k$, a single scalar relaxation is already near-optimal, and
RN sits at the PINN trace-error floor (the full relaxation-ladder
benchmark, over five seeds, is in Appendix~\ref{app:baselines}). The
regime that motivates partitioned FSI is the opposite. In cardiovascular
FSI~\citep{landajuela2018numerical} a light structure is immersed in a
heavy, nearly incompressible
fluid, and the \emph{added-mass} operator the fluid presents to the
interface dominates the \emph{low} modes and decays with $k$, making
$r_k$ strongly mode-dependent; DN then diverges with
$\rho(\TDN)=\max_k r_k=\mathrm{Ma}\gg1$, the added-mass
instability~\citep{CausinGerbeauNobile2005,BadiaNobileVergara2008,fernandez2013fully}.

\begin{figure}[t]
\centering
\begin{subfigure}[b]{0.49\linewidth}
  \centering
  \figmaybe{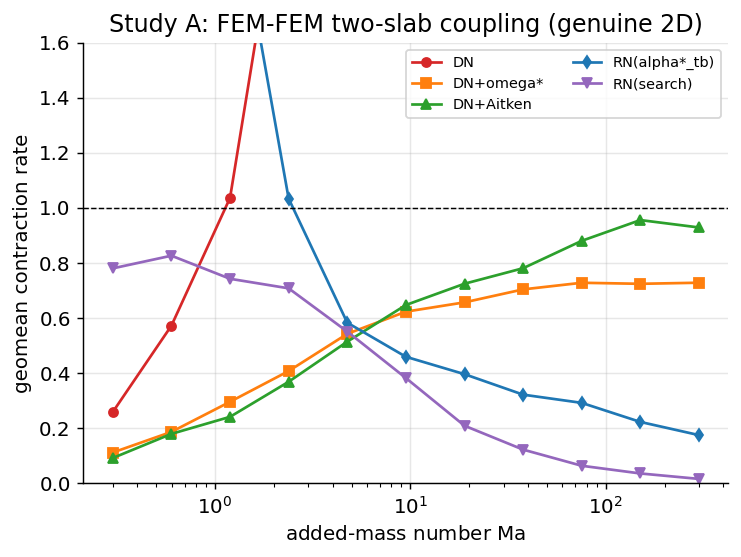}{Rate vs.\ Ma, FEM--FEM.}
  \caption{Geometric-mean contraction rate (defined in
  Section~\ref{sec:exp:addedmass}) vs.\ $\mathrm{Ma}$ in the 2D coupled
  iteration. Crossover at $\mathrm{Ma}\approx5$: relaxation saturates
  while RN keeps improving.}
\end{subfigure}\hfill
\begin{subfigure}[b]{0.49\linewidth}
  \centering
  \figmaybe{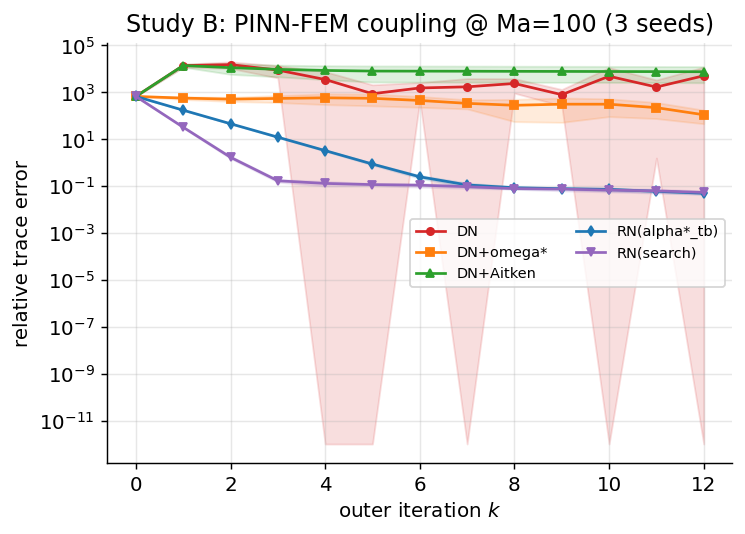}{PINN--FEM coupling, Ma=100.}
  \caption{A \emph{real PINN} on the solid slab, $\mathrm{Ma}=100$,
  three seeds: RN converges to the PINN floor while DN and Aitken
  diverge.}
\end{subfigure}
\caption{\textbf{The added-mass regime is where RN wins.} The two-slab
Poisson analogue emulates the interface spectral instability of added
mass (not its inertia). RN contracts ever faster as added mass grows
while scalar relaxation saturates; the effect persists from the operator
level through FEM--FEM to a trained-PINN subdomain.}
\label{fig:addedmass}
\end{figure}

A scalar Poisson problem has no inertia but reproduces this interface
\emph{spectral} instability exactly. We use a \emph{genuine 2D}
two-slab analogue (a thin high-diffusivity ``fluid'' slab stacked on a
thick ``solid'' slab, sharing a flat interface; the full geometry,
parameters, and the FEM--FEM versus trained-PINN study split are in
Appendix~\ref{app:addedmass-config} and Figure~\ref{fig:addedmass-config}),
with lateral interface modes $q_k=k\pi/W$ (the cosine harmonics
$\cos(q_k x)$ selected by the homogeneous-Neumann side walls of width $W$;
the boundary conditions that fix these wavenumbers are derived in
Appendix~\ref{app:addedmass-config}) and flat-slab eigenvalue
$\sigma(q)=\kappa\,q\coth(qL)$: a growing $\sigmaN(q)\approx\kappa_\tagN
q$ against a \emph{constant} $\sigmaF(q)\approx\kappa_\tagF/L_\tagF$ (the
lumped added-mass limit), so $r_k\approx\mathrm{Ma}/k$ with the single
dial $\mathrm{Ma}=\sigmaF(q_1)/\sigmaN(q_1)$. As $\mathrm{Ma}\to\infty$,
DN diverges ($\rho=\mathrm{Ma}$) and a scalar relaxation
\emph{saturates} at $(K-1)/(K+1)$ (it cannot straddle the mode spread),
whereas RN places the added mass $\alpha+\sigmaF$ in the
\emph{denominator} of its map and contracts ever faster.
Table~\ref{tab:addedmass} and Figure~\ref{fig:addedmass} bear this out
in a genuine 2D coupled iteration. Throughout, the reported
\emph{contraction rate} (the ``geomean contraction rate'' of
Figure~\ref{fig:addedmass}) is the geometric mean
$\bigl(\norm{e^{\Kouter}}/\norm{e^{0}}\bigr)^{1/\Kouter}$ of the
per-iteration interface-error reduction over the $\Kouter=12$ outer
iterations: ${<}1$ converges, ${>}1$ (labelled ``div.'') diverges. On
this metric,
below $\mathrm{Ma}\approx5$ tuned relaxation wins, but past the crossover
RN wins decisively. At
$\mathrm{Ma}=300$ scalar relaxation has saturated ($\omega^\ast=0.73$,
Aitken $0.93$) and DN has diverged, while RN contracts at $0.18$ from the
sweep-free $\alphastar$ and at $0.02$ from a searched optimum. The win
survives a \emph{trained-network} subdomain: with a real PINN on the
solid slab at $\mathrm{Ma}=100$ over three seeds, DN, Aitken and
$\omega^\ast$ all diverge or stall while RN converges to the PINN floor
($\sim5\times10^{-2}$) at a sweep-free rate $0.45\pm0.004$
(Figure~\ref{fig:addedmass}b). The spectral cap is harmless here because
it truncates only the high modes, which are the \emph{stable} ones in
this regime, so the resolvable band covers exactly the dangerous low
modes. RN's advantage is thus real and large precisely in the
FSI-relevant regime.

\begin{table}[t]
\centering
\caption{\textbf{RN wins in the large-added-mass regime} (genuine 2D
two-slab coupling, FEM--FEM realisation). Each entry is the
geometric-mean per-iteration contraction rate over $\Kouter=12$ outer
iterations, as defined in the text: ${<}1$ converges, ${>}1$ diverges
(marked ``div.''); the bold entry in each row is the best (smallest) rate
for that $\mathrm{Ma}$. As $\mathrm{Ma}$ grows, DN diverges and tuned
scalar relaxation ($\omega^\ast$, Aitken) saturates near
$(K{-}1)/(K{+}1)$, while RN contracts ever faster, the opposite ranking
to the low-added-mass disc-in-square (Appendix~\ref{app:baselines}).
$\alphastar$ is designed for the added-mass regime and is not meant for
$\mathrm{Ma}\lesssim2$.}
\label{tab:addedmass}
\small
\begin{tabular}{@{}lrrrrl@{}}
\toprule
$\mathrm{Ma}$ & DN ($\omega{=}1$) & DN$+\omega^\ast$ & DN$+$Aitken
 & RN($\alphastar$) & RN($\alpha$-search) \\
\midrule
$0.3$  & $0.26$  & $\mathbf{0.11}$ & $\mathbf{0.09}$ & $5.8$\,(div.) & $0.78$ \\
$9.5$  & $7.7$\,(div.) & $0.62$ & $0.65$ & $0.46$ & $\mathbf{0.38}$ \\
$37.8$ & $23$\,(div.)  & $0.70$ & $0.78$ & $0.32$ & $\mathbf{0.12}$ \\
$300$  & $235$\,(div.) & $0.73$ & $0.93$ & $0.18$ & $\mathbf{0.02}$ \\
\bottomrule
\end{tabular}
\end{table}

\section{Application to fluid--structure interaction with contact}
\label{sec:fsi}

We now move to the configuration that motivated the work: a rigid disc
(FEM) inside a viscous fluid (PINN, on the moving outer annulus). The
focus of this section is the meshless handling of \emph{topology change}
at contact.

\subsection{Model, discretisation, and vector RN}
\label{sec:fsi:model}

The fluid is governed by steady Stokes,
\begin{equation}
- \nabla \cdot \sigma_f(u, p) = -\rho_f g \, e_y \,\text{ in } \Omegaf(t), \qquad
\nabla \cdot u = 0 \,\text{ in } \Omegaf(t),
\label{eq:stokes}
\end{equation}
$\sigma_f(u,p)=-pI+2\mu_f\varepsilon(u)$, with no-slip on top/left/right
walls and Navier-slip on the bottom wall $\Sigma_p$, and $u=\dot
y_c(t)e_y$ on the moving disc surface $\SigmaI(t)$. The solid is a single
ODE on the centre height $y_c(t)$,
\begin{equation}
m_s \ddot y_c = -m_s g + F_y(t) + F_c(t), \qquad
F_y = \int_{\SigmaI(t)} (\sigma_f n) \cdot e_y \, ds ,
\label{eq:solid}
\end{equation}
discretised by backward Euler, with the normal contact reaction $F_c$
computed by an Alart--Curnier augmented-Lagrangian
projection~\citep{AlartCurnier1991}; a closed-form one-step
over-relaxation $\omega_{\mathrm{AL}}=m_s/(\gamma_c\dt^2)$
(Appendix~\ref{app:contact}) lands the disc at $y_c=R+\varepsilon_g$ to
machine precision (exact for the linearised backward-Euler position
map). The PINN fluid is a $(u,v,p)$ MLP (width 32, depth 4,
Tanh). \emph{At each time step the collocation set is regenerated on the
current fluid domain} $\Omegaf(t)=\Omega_+\setminus(B((0,y_c(t)),R)\cup\{y<\varepsilon_g\})$;
the topology change at contact is therefore carried by the sampler,
with no remeshing, no cut cell, and no level set
(Figure~\ref{fig:domains}c; the three near-contact configurations are
detailed in Figure~\ref{fig:contact-topology}). For Stokes the interface trace is
the vector velocity $u|_\SigmaI\in\R^2$ and the conjugate datum is the
traction $(\sigma_f n)|_\SigmaI\in\R^2$; the \emph{vector RN} iteration
is Algorithm~\ref{alg:dn-rn} with $\lambda,g$ promoted to vector fields
and $\SF,\SN$ the $(2n_{\mathrm{iface}})\times(2n_{\mathrm{iface}})$
vector Steklov operators. All spectral statements carry over verbatim.
The complete time-stepped solver, assembling the backward-Euler solid
step, the vector-RN inner coupling, the Alart--Curnier augmented-Lagrangian
contact projection, and the per-step collocation regeneration, is stated as
Algorithm~\ref{alg:contact-fsi} in Appendix~\ref{app:contact}; the network
architecture and the solver block diagram are in
Appendix~\ref{app:fsi-figs}.

\subsection{Static-Stokes validation and the vector-RN recipe}
\label{sec:fsi:recipe}

On a fixed-geometry manufactured Stokes problem the vector RN iteration
needs four practical ingredients to be transient-stable
(Appendix~\ref{app:E11-hyper}): the unit-correct pointwise Neumann
transfer~\eqref{eq:units-fix}; warm-started PINN inner solves (carrying
weights and Adam moments across outer iterations); joint
under-relaxation on $(\lambda,g)$ with $\omega\in[0.25,0.30]$; and an
early stop at the first sustained residual increase, since the
stabilised iteration is transient-stable rather than asymptotically
stable. With $(\alpha,\omega,\Kouter)=(30,0.30,17)$ the iteration
contracts monotonically $3.87\to0.158$ on residual and reaches a
relative trace error $0.136$ versus the manufactured solution. The PINN
vector Steklov operator, estimated by the vector Fourier probe, is
capped at $k^\star\approx4$ ($\abs{\sigmaN}^{\max}_{k\ge1}=3.06$ at
$k=2$ vs.\ the FE reference's $119.7$); the cap, not the coupling
scheme, is the binding accuracy limit (Appendix~\ref{app:E11-hyper}).

\subsection{Free fall and contact by collocation exclusion}
\label{sec:fsi:contact}

\paragraph{Free fall (E12).} A disc released from rest at $y_c(0)=0.5$
($g=1$, $\mu_f=\rho_f=1$, $\rho_s=2$) falls under drag balance: the
vertical force $F_y$ plateaus at the buoyancy level
$\rho_f\pi R^2 g=0.196$ through $t\approx0.25$, and a Newton-balance fit
on $t\in[0.04,0.25]$ gives an effective drag coefficient $c\in[0.19,
0.32]$, consistent with a confined, sub-terminal disc
(Appendix~\ref{app:E12-extra}). Past $t\approx0.30$ the warm-start drift
quantified in Section~\ref{sec:fsi:limits} dominates; this safe window
of ${\sim}60$ steps is enough lead-in for the contact experiment.

\begin{figure}[t]
\centering
\figmaybe[0.95\linewidth]{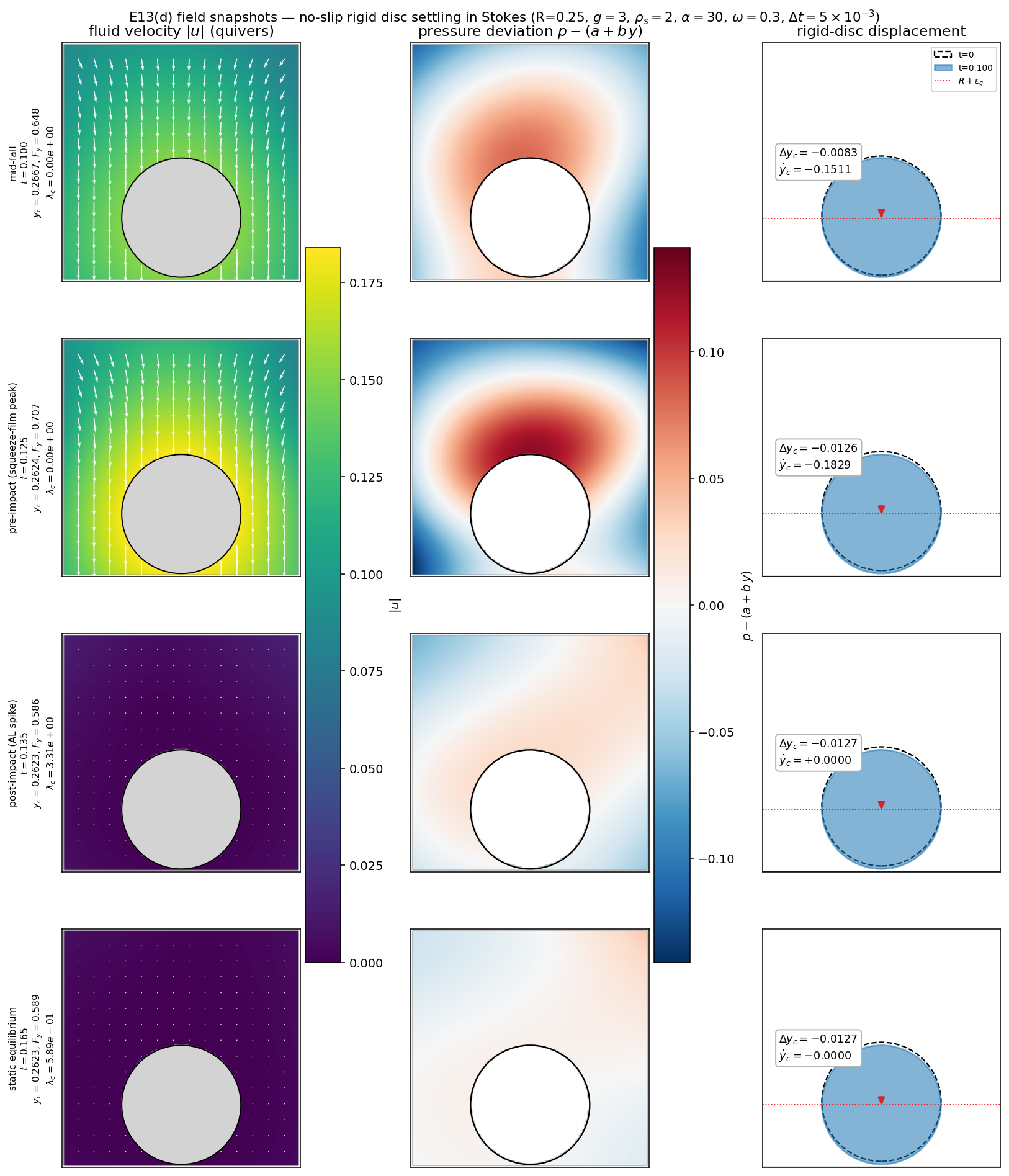}{E13(d) field snapshots: fluid velocity magnitude with quivers (left), pressure deviation (centre), and rigid-disc displacement (right) at four key time instants.}
\caption{E13(d) snapshots of the PINN--FEM coupled solution at four
instants spanning descent--impact--equilibrium: drag-balanced mid-fall
($t=0.100$), squeeze-film pre-impact ($t=0.125$), post-impact contact
spike ($t=0.135$, $\lambda_c=3.31$), and static equilibrium ($t=0.165$,
$\lambda_c=\Pi=0.589$). Pressure panels show $p-(a+by)$ with the linear
trend removed so the Stokes dipole and lift band are visible. The disc
crosses the contact event with the moving fluid topology handled by
collocation exclusion alone; the static-equilibrium reaction matches the
submerged weight to within $0.4\%$ under interface refinement.}
\label{fig:e13d-fields}
\end{figure}

\paragraph{Contact (E13d).} We shift to $y_c(0)=0.275$, $g=3$, so impact
occurs at $t^\star=\sqrt{2\cdot0.025/3}=0.129$, inside the safe window.
\emph{The PINN carries the disc through the contact event by collocation
exclusion alone}: as the gap closes, the strip $\{y<\varepsilon_g\}$
simply leaves the collocation set (Figure~\ref{fig:domains}c),
with no remeshing and no cut cells. The Alart--Curnier projector lands
the disc at $y_c=R+\varepsilon_g=0.26227$ to machine precision, and the
contact reaction converges to the submerged weight: at all three
interface resolutions $n_{\mathrm{iface}}\in\{16,32,64\}$ the
post-impact tail $\lambda_c(t\ge0.20)$ matches $\Pi=0.589$ to within
$0.4\%$ (Figure~\ref{fig:e13d-fields}, the static-equilibrium claim is
mesh-converged). Figure~\ref{fig:e13d-fields} shows the coupled solution
across the descent--impact--equilibrium cycle, the topology change
absorbed entirely by the sampler.

\subsection{Limitations of the current methodology}
\label{sec:fsi:limits}

The coupling has three limitations; the first two are quantified
in the appendix. \emph{(i) Warm-start drift.} The fluid PINN is warm-started
across time steps and each per-step RN solve is only approximately
converged, so the network slides off the Stokes manifold; a
training-independent interface-residual monitor grows by a factor
${\sim}80$ over a fifty-step trajectory if nothing is done. Periodically
restarting the network from scratch every three to five steps holds the
drift to a factor ${\sim}4$--$5$ at no extra inner-training cost and
seed-stably, but over much longer ($150$-step) horizons no tested remedy
\emph{eliminates} the drift; it only slows it
(Appendix~\ref{app:E21}). \emph{(ii) Pre-impact dynamics are
under-resolved.} The pre-impact squeeze-film $F_y$ signatures are not
converged in the PINN inner-solve budget (they are non-monotone, and at
the highest budget the disc never even enters the contact band), and a
matched FEM--FEM Stokes benchmark at the same operating point shows that
a fully resolved fluid never reaches the wall within $T=0.30$,
under-resolving the bulk Stokes drag by ${\sim}80\times$. We therefore
report the pre-impact dynamic signatures as PINN under-resolution
artefacts rather than physical findings (Appendix~\ref{app:contact-sweeps}).
The validated dynamic results are those for which the disc has already
touched the wall and the AL projector dominates: $\lambda_c\to\Pi$ and
$y_c\to R+\varepsilon_g$ to machine precision. Finally, steady Stokes
carries no fluid inertia and the AL contact is perfectly inelastic, so
the model admits no literal post-impact rebound. \emph{(iii) The
meshless advantage is robustness, not speed.} On a matched static
coupling (natural-orientation Poisson, $n_{\mathrm{iface}}=32$, the same
RN outer iteration, only the disc solver swapped), the PINN subdomain
costs ${\sim}4600\times$ more wall-clock per outer iteration than a
conforming P1 finite-element subdomain (a network training versus one
sparse solve, at our $10^3$-Adam-step inner budget). The case for the
PINN fluid is therefore not efficiency on a fixed domain but the
meshless treatment of moving geometry and topology change at contact,
where the body-fitted alternative pays in remeshing, ALE, or cut-cell
machinery instead.

\emph{Is the inner-solve floor removable?} Both limitations bottom out
at the same PINN inner-solve floor, and every experiment above trained
with Adam on a uniform collocation set under a soft incompressibility
penalty. To separate an \emph{optimisation} limit from a
\emph{representation} one we re-ran the two diagnostics under three
heavy gradient-descent levers (Appendix~\ref{app:E22}): a quasi-Newton
(L-BFGS) inner solve, adaptive collocation concentrated in the
lubrication layer, and improved conditioning (an exactly
divergence-free stream-function ansatz; gradient-norm loss balancing).
None lifts the floor at fixed width: the Stokes Steklov cap stays below
${\approx}50\%$ of the Taylor--Hood reference's leading $k\ge1$
eigenvalue (L-BFGS roughly triples the Adam value but never resolves
even the first mode in full; gradient-norm balancing is actively
harmful; the divergence-free ansatz is no better than
velocity--pressure), and the squeeze-film drag ratio remains
${\approx}0.03$ and \emph{deepens} as the gap narrows, with adaptive
sampling no better than uniform. Removing gradient descent overturns
the verdict on the cap: a random-feature inner solve (frozen
random-Fourier features with a linear head fitted by a single
least-squares solve, Appendix~\ref{app:E23}) flattens the per-mode
Steklov ratio, which then tracks a fully resolved finite-element solver
to within a few percent at every probed mode; the remaining
frequency-flat ${\approx}0.57$ offset is reproduced by the
finite-element solver itself when its traction is read through the same
pointwise functional, so it is a traction-reader convention, not a
representation deficit. We therefore reclassify the spectral cap as a
gradient-descent \emph{optimisation} pathology (spectral bias), not a
representation-capacity wall. The squeeze-film drag deficit, by
contrast, is not closed by any inner solve we tried and remains an
under-resolution limitation of the time-stepped application.

\section{Discussion and conclusion}
\label{sec:discussion}

PINN--FEM coupling is most clearly understood as a partitioned domain
decomposition method with a learned subdomain solver. Under that view
the classical Dirichlet--Neumann instability and its Robin--Neumann cure
transfer directly, modulo a discrete units conversion and the operational
specifics of training a network at each outer iteration. The
operator-spectral theory predicts the empirical PINN--FEM rates to
within $7\%$ on 1D and 2D Poisson; the Fourier Steklov estimator
recovers the resolvable interface eigenvalues an order of magnitude
cheaper and more accurately than nodal-basis or random-trace probes; and
RN's per-mode impedance matching wins decisively in the large-added-mass
regime that motivates cardiovascular FSI, where a tuned scalar
relaxation saturates. On the moving-geometry Stokes/rigid-disc
application, the PINN's collocation-based discretisation handles topology
change at contact without any meshing operation, and the
static-equilibrium contact reaction matches the submerged weight to
within $0.4\%$ under mesh refinement.

The binding limitations are the PINN inner-solve floor (a training, not a
coupling, problem) and the warm-start drift of the time-stepped fluid; a
matched FEM--FEM benchmark further shows the pre-impact dynamic signatures
to be under-resolution artefacts. The inner-solve floor is specifically an
\emph{optimisation} effect: no gradient-descent lever (second-order
optimisation, adaptive collocation, improved conditioning) lifts the
interface spectral cap at fixed width (Appendix~\ref{app:E22}), yet a
non-gradient random-feature least-squares inner solve flattens it and
recovers the interface operator to finite-element accuracy
(Appendix~\ref{app:E23}), so the cap is a gradient-descent spectral-bias
pathology rather than a representation-capacity wall. The natural next steps are
to give the fluid inertia (unsteady Stokes / Navier--Stokes) to obtain a
\emph{physical} post-impact rebound (validated against a matched
fully-resolved reference so it is not mistaken for an under-resolution
artefact), to promote the PINN to a \emph{deformed} interface for a
deformable solid, and to extend to 3D. A complementary route to the
per-iteration retraining cost quantified in
Section~\ref{sec:fsi:limits} is to replace the freshly-trained PINN
subdomain with a pretrained operator-learning surrogate
(DeepONet~\citep{LuJinPangZhangKarniadakis2021}, a Fourier neural
operator~\citep{LiKovachkiEtAl2021}, or a non-intrusive reduced-order
surrogate~\citep{HesthavenUbbiali2018}): an interface-conditioned operator
network would realise $\SN$ in a single forward pass per outer
iteration, trading per-step training for an offline cost, and whether
its spectral cap behaves like the PINN's is a question the Fourier
probe of Section~\ref{sec:theory:sn-estimator} could answer directly.
None of these is demonstrated
here; the defensible claims are the spectral-prediction agreement on
Poisson coupling, the added-mass RN win, the static-Stokes vector-RN
recipe, the free-fall drag match, and the static-equilibrium contact
result. Code, seeds, environment details, and the cached artefact behind
every figure are catalogued in Appendix~\ref{app:repro}.

\section*{Acknowledgments}

This work was performed under the auspices of the U.S. Department of Energy by Lawrence Livermore National
Laboratory under Contract DE-AC52-07NA27344 and was supported by the LLNL-LDRD Program under Project No. 24-ERD-022. LLNL-JRNL-2019742.

\bibliographystyle{plainnat}
\bibliography{refs}

\clearpage

\appendix

\section{Sign and discretisation conventions}
\label{app:conventions}

All fluxes are reported in the FEM-outward normal sign $\nF$, so the
discrete trace-to-flux maps on the two subdomains have opposite signs in
the interface equation. By construction the $n_{\mathrm{iface}}$ FEM
boundary nodes on $\SigmaI$ are reused verbatim as the PINN's interface
evaluation points (in 2D equispaced in $\theta$; in 1D the single point
$\{L\}$). Both 2D mesh builders (annulus and disc) accept the same
$n_{\mathrm{iface}}$ and produce identical interface DOF locations, so
the trace vector $\lambda\in\R^{n_{\mathrm{iface}}}$ is directly
interchangeable between solvers and no scattered-data interpolation
enters the trace-transfer step. This shared-grid choice is what makes the
interface mass matrix $\Miface$ of~\eqref{eq:units-fix} a single
well-defined object, which in turn makes the pointwise-flux unit
conversion clean.

\begin{remark}[Discrete units pitfall]
\label{rem:units}
The FEM side outputs the integrated reaction
\begin{equation}
\gint \;=\; (Ku-f)\big|_\SigmaI ,
\label{eq:gint}
\end{equation}
i.e.\ the restriction of the residual force vector to the interface
nodes, which carries units of force (a flux already integrated against
the FEM test functions), whereas the PINN's Neumann loss expects a
pointwise flux $\gpw$. The conversion is via the interface mass matrix,
\begin{equation}
\gpw = \Miface^{-1} \, \gint .
\label{eq:units-fix}
\end{equation}
Omitting~\eqref{eq:units-fix} scales the Neumann datum by $h_b$ and
inverts both the DN and RN headline conclusions in 2D; after the fix,
PINN--FEM tracks FEM--FEM to within $7\%$ on all rates we measured.
\end{remark}

\section{Boundary conditions on the natural and FSI-aligned configurations}
\label{app:bcs}

\begin{figure}[h]
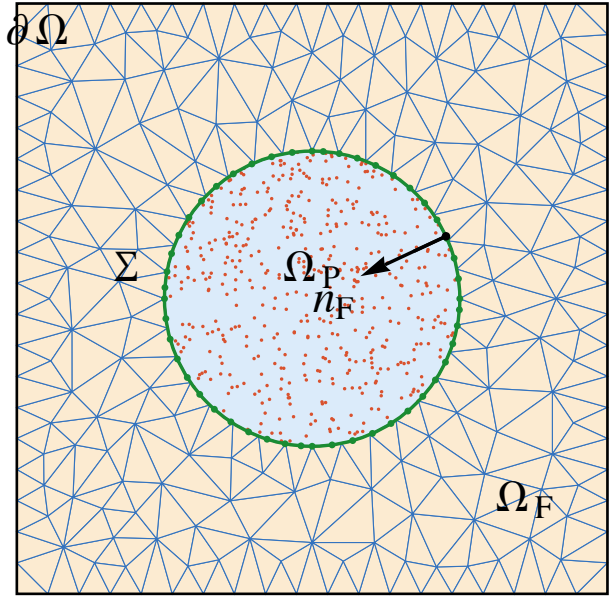

\centering
\figmaybe[0.7\linewidth]{figures/domain_natural.pdf}{Natural-orientation boundary conditions.}
\caption{Natural orientation (scalar Poisson). Dirichlet (or Robin) on
the FEM annulus side of $\SigmaI$, Neumann on the PINN disc side;
homogeneous Dirichlet on the outer square boundary.}
\label{fig:app-bcs:natural}
\end{figure}

\begin{figure}[h]
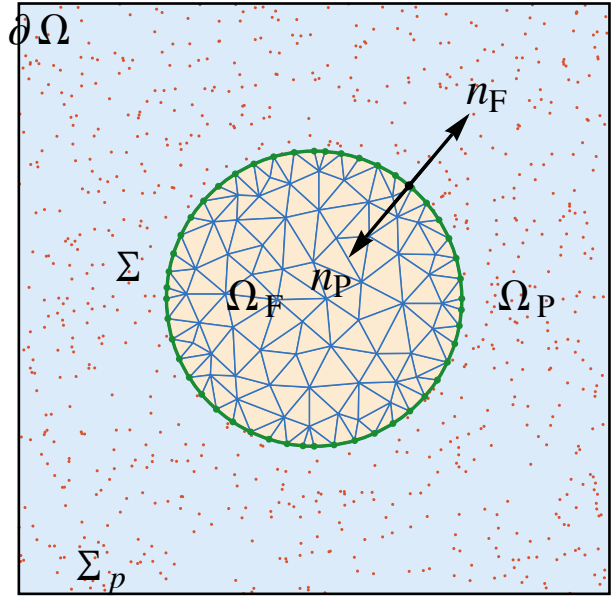

\centering
\figmaybe[0.7\linewidth]{figures/domain_config.pdf}{FSI-aligned-orientation boundary conditions.}
\caption{FSI-aligned orientation (Stokes--rigid-disc). The PINN fluid
annulus takes the velocity trace (Robin) on $\SigmaI$ and the FEM solid
disc takes the integrated traction (Neumann); no-slip on top/left/right
walls and Navier-slip on the bottom wall $\Sigma_p$.}
\label{fig:app-bcs:fsi}
\end{figure}

For convenience we collect, in a single visual reference, the boundary
conditions imposed on each piece of $\partial\OmegaF$,
$\partial\OmegaN$, and the shared interface $\SigmaI$ in the two role
assignments of Section~\ref{sec:problem}.
Figure~\ref{fig:app-bcs:natural} shows the \emph{natural} (scalar
Poisson) orientation used in E5/E7/E10, and
Figure~\ref{fig:app-bcs:fsi} the \emph{FSI-aligned}
(Stokes--rigid-disc) orientation of Section~\ref{sec:fsi}.

\section{Steklov--Poincar\'e derivation of the DN and RN
maps\texorpdfstring{~\eqref{eq:dn-rn-map}}{}}
\label{app:dn-rn-derivation}

This appendix derives the two trace-error recurrences of
Equation~\eqref{eq:dn-rn-map} step by step from the subdomain
trace-to-flux relations, making explicit why the bulk unknowns can be
eliminated (the point clarified in Section~\ref{sec:theory:interface}).

\paragraph{Subdomain responses and the interface equation.} For a given
interface trace $t$ on $\SigmaI$, each subdomain problem is well-posed
(Dirichlet data $t$ plus the subdomain forcing), so its interior solution
is the \emph{unique} PDE extension of $t$ and the interior degrees of
freedom are completely determined by $t$. Their only trace on the coupling
is the boundary flux, which by~\eqref{eq:steklov-def}--\eqref{eq:fluxes}
is affine in $t$,
\begin{equation}
\phi_{\tagF}(t) := \dnF\uF(t)\big|_\SigmaI = \SF\,t - \chiF , \qquad
\phi_{\tagN}(t) := \dnF\uN(t)\big|_\SigmaI = -\SN\,t + \chiN ,
\label{eq:app-flux-maps}
\end{equation}
the sign flip on the PINN side coming from $\nN=-\nF$. The exact coupled
solution has a single trace $\lambda^\star$ on both sides (trace
continuity) and matching fluxes there (flux continuity),
$\phi_{\tagF}(\lambda^\star)=\phi_{\tagN}(\lambda^\star)$, i.e.
\begin{equation}
\SF\lambda^\star-\chiF = -\SN\lambda^\star+\chiN
\;\Longleftrightarrow\;
(\SF+\SN)\,\lambda^\star = \chiF+\chiN ,
\label{eq:app-interface}
\end{equation}
which is~\eqref{eq:steklov-interface}. Because~\eqref{eq:app-flux-maps}
has already condensed each interior solve into an operator on $\SigmaI$,
every iteration below lives on the interface alone.

\paragraph{Dirichlet--Neumann.} At step $k$ the FEM side solves the
Dirichlet problem with trace $\lambda^k$ and exposes the flux
$\phi_{\tagF}(\lambda^k)=\SF\lambda^k-\chiF$. The PINN side then solves the
Neumann problem carrying that same flux (flux continuity), so its trace
$t$ satisfies $\phi_{\tagN}(t)=\phi_{\tagF}(\lambda^k)$, i.e.
$-\SN t+\chiN = \SF\lambda^k-\chiF$, giving
\begin{equation}
t = \SN^{-1}\bigl(\chiN+\chiF-\SF\lambda^k\bigr)
  = \SN^{-1}(\SF+\SN)\lambda^\star - \SN^{-1}\SF\,\lambda^k ,
\end{equation}
where the second equality used~\eqref{eq:app-interface}. The DN update is
$\lambda^{k+1}=t$, so subtracting $\lambda^\star=\SN^{-1}\SN\lambda^\star$
and writing $e^k=\lambda^k-\lambda^\star$,
\begin{equation}
e^{k+1} = \SN^{-1}\SF\,\lambda^\star - \SN^{-1}\SF\,\lambda^k
        = -\,\SN^{-1}\SF\,e^k =: \TDN\,e^k ,
\end{equation}
the first map in~\eqref{eq:dn-rn-map}.

\paragraph{Robin--Neumann.} The Neumann (PINN) side supplies the flux
response at the current trace, $g^k=\phi_{\tagN}(\lambda^k)=\chiN-\SN\lambda^k$.
The FEM side solves the Robin problem
$\dnF\uF+\alpha\uF=\alpha\lambda^k+g^k$ on $\SigmaI$; using
$\dnF\uF=\SF\lambda^{k+1}-\chiF$ for its new trace $\lambda^{k+1}$ this is
$(\SF+\alpha I)\lambda^{k+1}-\chiF=\alpha\lambda^k+g^k$, i.e.
\begin{equation}
(\alpha I+\SF)\,\lambda^{k+1}
  = \alpha\lambda^k + (\chiN-\SN\lambda^k) + \chiF
  = (\alpha I-\SN)\,\lambda^k + (\chiF+\chiN) .
\label{eq:app-rn-fixed}
\end{equation}
This is a preconditioned Richardson iteration for the interface
equation~\eqref{eq:app-interface}: its fixed point
$(\alpha I+\SF)\lambda^\star=(\alpha I-\SN)\lambda^\star+(\chiF+\chiN)$
collapses to $(\SF+\SN)\lambda^\star=\chiF+\chiN$, recovering the exact
solution. Subtracting the fixed point from~\eqref{eq:app-rn-fixed} gives
the error recurrence
\begin{equation}
e^{k+1} = (\alpha I+\SF)^{-1}(\alpha I-\SN)\,e^k .
\label{eq:app-rn-map}
\end{equation}
Under the working hypothesis of Section~\ref{sec:theory:theorems} that
$\SF$ and $\SN$ are simultaneously diagonalisable, both factors are
diagonal in the shared eigenbasis and therefore commute, so
$(\alpha I+\SF)^{-1}(\alpha I-\SN)=(\alpha I-\SN)(\alpha I+\SF)^{-1}=\TRN(\alpha)$,
the second map in~\eqref{eq:dn-rn-map}, with per-mode factor
$\mu_k(\alpha)=(\alpha-\sigma_{\tagN,k})/(\alpha+\sigma_{\tagF,k})$ used in
the proof of Theorem~\ref{thm:rn}. (When the eigenbases differ the two
orderings remain similar, hence isospectral, so the contraction rate is
unchanged; see Appendix~\ref{app:proofs}.) The under-relaxed update of
Algorithm~\ref{alg:dn-rn},
$\lambda^{k+1}\!\gets\!(1-\omega)\lambda^k+\omega\,\hat\lambda^{k+1}$ with
$\hat\lambda^{k+1}$ the right-hand side of~\eqref{eq:app-rn-fixed},
replaces $\TRN(\alpha)$ by $(1-\omega)I+\omega\,\TRN(\alpha)$; the Poisson
experiments use $\omega=1$ and the vector-Stokes ones $\omega\approx0.3$.

\section{Proofs of Theorems~\ref{thm:dn}--\ref{thm:pinn-floor}}
\label{app:proofs}

\subsection{Proof of Theorem~\ref{thm:dn}}

Under the simultaneous-diagonalisation hypothesis, $\SF$ and $\SN$ share
an eigenbasis $\{v_k\}$ with eigenvalues
$\sigma_{\tagF,k},\sigma_{\tagN,k}>0$. The map $\TDN=-\SN^{-1}\SF$ is
therefore diagonal in $\{v_k\}$ with eigenvalues
$-\sigma_{\tagF,k}/\sigma_{\tagN,k}$; the spectral radius is
$\max_k\sigma_{\tagF,k}/\sigma_{\tagN,k}$. Geometric divergence whenever
this maximum exceeds $1$ is the standard linear-iteration consequence.
\qed

\subsection{Proof of Theorem~\ref{thm:rn}}
\label{app:proofs:rn}

\emph{(i) Rate formula.} Diagonalising $\TRN(\alpha)$ in the shared
eigenbasis gives the per-mode factors
$\rho_k(\alpha)=\abs{\alpha-\sigma_{\tagN,k}}/(\alpha+\sigma_{\tagF,k})$
and $\rho(\TRN(\alpha))=\max_k\rho_k(\alpha)$.

\emph{(ii) Contractive band.} Fix $\alpha>0$ and a mode $k$. If
$\sigma_{\tagN,k}\le\alpha$ then
$\abs{\alpha-\sigma_{\tagN,k}}=\alpha-\sigma_{\tagN,k}<\alpha+\sigma_{\tagF,k}$
holds unconditionally (both eigenvalues are positive). If
$\sigma_{\tagN,k}>\alpha$ then
$\abs{\alpha-\sigma_{\tagN,k}}<\alpha+\sigma_{\tagF,k}
\iff \sigma_{\tagN,k}-\sigma_{\tagF,k}<2\alpha$. Hence
$\rho(\TRN(\alpha))<1$ if and only if
$\alpha>\alpha_{\min}=\tfrac12\max_k(\sigma_{\tagN,k}-\sigma_{\tagF,k})_+$,
and the band $(\alpha_{\min},\infty)$ is never empty (as
$\alpha\to\infty$, $\rho_k(\alpha)\to1^-$ for every $k$).

\emph{(iii) Equioscillation in the aligned case.} Suppose
$\sigma_{\tagF,k}=\sigma_{\tagN,k}=\sigma_k$ with
$\sigma_k\in[m,n]:=[\sigmaF^{\min},\sigmaN^{\max}]$, so
$\rho(\TRN(\alpha))\le\max_{\sigma\in[m,n]}h_\sigma(\alpha)$ with
$h_\sigma(\alpha)=\abs{\alpha-\sigma}/(\alpha+\sigma)$. For fixed
$\alpha$, $\sigma\mapsto h_\sigma(\alpha)$ is decreasing on
$\sigma\le\alpha$ and increasing on $\sigma\ge\alpha$, so the maximum
over $[m,n]$ is attained at an endpoint:
$\max\{h_m(\alpha),h_n(\alpha)\}$. As functions of $\alpha\in[m,n]$,
$h_m$ is increasing and $h_n$ decreasing, so the minimax is at the
crossing $h_m(\alpha)=h_n(\alpha)$, i.e.\
$(\alpha-m)(n+\alpha)=(n-\alpha)(\alpha+m)$, which reduces to
$\alpha^2=mn$: $\alphastar=\sqrt{mn}=\sqrt{\sigmaF^{\min}\sigmaN^{\max}}$,
with equioscillated value
$h_m(\alphastar)=h_n(\alphastar)=(\sqrt\kappa-1)/(\sqrt\kappa+1)$,
$\kappa=n/m$. \qed

\paragraph{Scope of $\alphastar$ for mismatched spectra.} When the two
spectra do not pair mode-by-mode, step (iii) no longer applies:
$\alphastar$ is the minimiser of the \emph{aligned-case envelope}, not of
the true $\max_k\rho_k(\alpha)$, and the iteration optimum can sit
elsewhere in the band. This is exactly what the numerical
$\alpha$-search shows on the 2D natural orientation (closed-form
$\alphastar=4.16$ vs.\ searched optimum $\alpha=10.31$,
Appendix~\ref{app:inner-pinn:e10}): the $2.5\times$ shift is the
practical price of the alignment hypothesis, mitigated by the bowl
around the optimum being broad.

\subsection{Proof of Theorem~\ref{thm:cotraining}}
\label{app:proofs:cotraining}

We give the full argument behind the energy-decrease statement. Work in
the shared eigenbasis $\{v_k\}$ of Theorems~\ref{thm:dn}--\ref{thm:rn},
write the FEM-side and PINN-side traces as
$\uF^k=\sum_k a_k^F v_k$, $\uN^k=\sum_k a_k^N v_k$, and let
$e^k=\lambda^k-\lambda^\star$ be the trace mismatch with components
$e_k^{(j)}$. We make the following standing assumptions, matching the
co-training implementation (Appendix~\ref{app:inner-pinn:e7}):
\begin{enumerate}[leftmargin=*,itemsep=1pt,topsep=1pt]
\item \textbf{Exact Robin FEM solve.} At each outer step the FEM Robin
sub-problem is solved exactly, so the FEM contribution to $E$ is the
exact Dirichlet energy $\tfrac12\norm{\nabla\uF^k}^2 =
\tfrac12\sum_k\sigma_{\tagF,k}(a_k^F)^2$ of the harmonic extension of the
current trace. The second equality is Green's first identity specialised
to a harmonic extension: $\uF^k$ solves $\Delta\uF^k=0$ on $\OmegaF$ with
the current trace on $\SigmaI$ and homogeneous data elsewhere, so
\begin{equation}
  \tfrac12\norm{\nabla\uF^k}_{L^2(\OmegaF)}^2
  = \tfrac12\!\int_{\partial\OmegaF}\!\uF^k\,\partial_n\uF^k
    - \tfrac12\!\int_{\OmegaF}\!\uF^k\,\Delta\uF^k
  = \tfrac12\!\int_{\SigmaI}\!\uF^k\,\partial_n\uF^k
  = \tfrac12\,\langle\,\uF^k|_{\SigmaI},\,\SF\,\uF^k|_{\SigmaI}\rangle ,
  \label{eq:dirichlet-steklov}
\end{equation}
where the bulk term vanishes because $\uF^k$ is harmonic and the boundary
integral collapses to $\SigmaI$ (the far face contributes nothing since
$\uF^k=0$ there, the side walls nothing since $\partial_n\uF^k=0$ there),
and $\partial_n\uF^k|_{\SigmaI}=\SF\,\uF^k|_{\SigmaI}$ is precisely the
Steklov (Dirichlet-to-Neumann) map. Expanding the trace in the shared
$\SF$-eigenbasis, $\uF^k|_{\SigmaI}=\sum_k a_k^F v_k$ with $\SF
v_k=\sigma_{\tagF,k}v_k$ and orthonormal $\{v_k\}$, diagonalises the
quadratic form to $\tfrac12\sum_k\sigma_{\tagF,k}(a_k^F)^2$. In words: the
bulk energy of the exact FEM extension is, mode by mode, the squared
trace amplitude weighted by that mode's Steklov eigenvalue, which is what
lets the FEM term of $E$ be written in interface coordinates alone.
\item \textbf{Gradient PINN step.} The PINN side performs one (or a few)
gradient-descent step(s) of size $\eta$ on
$\Ll_{\tagN}(\theta)=\tfrac12\sum_k\sigma_{\tagN,k}(a_k^N-\lambda_k)^2 +
\Ll_{\tagN}^{\mathrm{inh}}$, whose stationary point is the Neumann
response $a_k^N=\lambda_k-\sigma_{\tagN,k}^{-1}\chi_{\tagN,k}$; near that
point $\Ll_{\tagN}$ is quadratic with Hessian eigenvalues
$\sigma_{\tagN,k}>0$.
\item \textbf{Quadratic regime.} The Adam preconditioner is, to leading
order in $\eta$, a positive diagonal scaling; we absorb it into $\eta$
and treat the PINN update as preconditioned gradient descent on the
quadratic model, valid for $\eta$ below the inverse of the largest
Hessian eigenvalue.
\end{enumerate}
Define the per-mode energy contribution
$E_k = \tfrac12\sigma_{\tagF,k}(a_k^F)^2 + \tfrac12\sigma_{\tagN,k}(a_k^N-\lambda_k)^2 +
\tfrac{\alpha}{2}(e_k)^2$, so $E^k=\sum_k E_k + \text{(forcing)}$. The RN
update~\eqref{eq:dn-rn-map} acts mode-wise as
$e_k^{k+1}=\mu_k(\alpha)\,e_k^k$ with
$\mu_k(\alpha)=(\alpha-\sigma_{\tagN,k})/(\alpha+\sigma_{\tagF,k})$, so
the interface term of $E_k$ changes by
$\tfrac{\alpha}{2}(\mu_k^2-1)(e_k)^2$. The exact Robin solve makes the
FEM term track the trace, contributing a non-positive change to leading
order, and one gradient step on the convex quadratic $\Ll_{\tagN}$
decreases the PINN term by $\eta\,\sigma_{\tagN,k}^2(a_k^N-\lambda_k)^2 +
O(\eta^2)$. Collecting the mode-wise changes,
\begin{equation}
E^{k+1}-E^k \;=\; -\frac{\alpha}{2}\sum_k\bigl(1-\mu_k(\alpha)^2\bigr)(e_k)^2
\;+\; O(\eta)
\;\le\; -\frac{\alpha}{2}\bigl(1-\rho(\TRN(\alpha))^2\bigr)\norm{e^k}^2_{L^2(\SigmaI)} + O(\eta),
\label{eq:energy-decrease}
\end{equation}
where $\rho(\TRN(\alpha))=\max_k\abs{\mu_k(\alpha)}$ and the $O(\eta)$
term collects the second-order optimiser remainder and the bounded
inhomogeneous forcing $\chi_{\tagF}+\chi_{\tagN}$ (constant across
iterations, hence a fixed offset, not a growing term). For any $\alpha$
in the contractive band of Theorem~\ref{thm:rn} we have
$\rho(\TRN(\alpha))<1$, so $1-\rho(\TRN(\alpha))^2>0$ and the
leading term in~\eqref{eq:energy-decrease} is strictly negative whenever
$e^k\neq0$. Choosing the learning rate $\eta$ small enough that the
$O(\eta)$ remainder is dominated by this leading term (concretely,
$\eta \le c\,\alpha(1-\rho^2)\norm{e^k}^2/(\text{remainder constant})$
for a fixed $c<1$) gives $E^{k+1}\le E^k$ up to the bounded forcing
offset, which is the claim. Note that the admissible step size shrinks
with $\norm{e^k}^2$: for a \emph{fixed} $\eta$ the strict decrease
holds only while the trace mismatch dominates,
$\norm{e^k}^2\gtrsim O(\eta)$, after which the $O(\eta)$ remainder
takes over and the energy is merely bounded. This terminal
neighbourhood is precisely the ``bounded forcing term'' of the
statement---the optimisation analogue of the $O(\epsN)$ floor of
Theorem~\ref{thm:pinn-floor}. \qed

\subsection{Discrete consequences when eigenbases differ}

If $\SF$ and $\SN$ do not share an eigenbasis, the closed-form
$\alphastar$ is no longer the iteration optimum. The operator-norm bound
$\rho(\TRN(\alpha))\le\norm{(\alpha I-\SN)(\alpha I+\SF)^{-1}}$ remains
valid but is generally loose; in practice one can refine the closed-form
value by a numerical $\alpha$-search (E10,
Appendix~\ref{app:inner-pinn:e10}). On the 2D natural orientation the
searched optimum sits $2.5\times$ from the closed-form value
($\alphastar=4.16$ vs.\ $\alpha=10.31$), the mismatched-spectra effect
quantified in the scope paragraph of Appendix~\ref{app:proofs:rn}; the
bowl around the optimum is broad, and the empirical minimum is itself
protocol-dependent (under the cold-start benchmark of
Appendix~\ref{app:inner-pinn:e6} the closed-form value outperforms the
searched one), so the sweep-free default is a serviceable operating
point either way.

\subsection{Proof of Theorem~\ref{thm:pinn-floor}}
\label{app:proofs:pinn-floor}

Substituting the realised operator $\hSN=\SN+\EN$ into the RN map
of~\eqref{eq:dn-rn-map} and isolating the additive residual gives, by
linearity alone,
\begin{equation}
e^{k+1}=(\alpha I-\hSN)(\alpha I+\SF)^{-1}e^k+\tilde r^k
       =\underbrace{\TRN(\alpha)e^k-\EN(\alpha I+\SF)^{-1}e^k}_{=\,\That(\alpha)e^k}
        +\tilde r^k .
\label{eq:pinn-recurrence}
\end{equation}
Unrolling~\eqref{eq:pinn-recurrence} from $e^0$ telescopes to
$e^k=\That^k e^0+\sum_{j=0}^{k-1}\That^{\,k-1-j}\tilde r^j$. Taking norms
and using submultiplicativity, $\norm{\That^m}\le\rhohat^m$ with
$\rhohat=\norm{\That(\alpha)}$, and the uniform residual bound
$\sup_j\norm{\tilde r^j}\le C_\alpha\epsN$,
\begin{equation}
\norm{e^k}\le\rhohat^k\norm{e^0}+\Bigl(\sum_{i=0}^{k-1}\rhohat^i\Bigr)
\,C_\alpha\epsN
\le\rhohat^k\norm{e^0}+\frac{C_\alpha\epsN}{1-\rhohat},
\end{equation}
the geometric series converging because $\rhohat<1$;
letting $k\to\infty$ gives the $\limsup$ in~\eqref{eq:pinn-floor}. The
norm bound~\eqref{eq:rhohat} is the triangle inequality applied to
$\That=\TRN-\EN(\alpha I+\SF)^{-1}$ together with
$\norm{(\alpha I+\SF)^{-1}}=1/(\alpha+\sigmaF^{\min})$, which holds because
$\SF$ is SPD so $\alpha I+\SF$ has smallest eigenvalue
$\alpha+\sigmaF^{\min}$. Unlike Theorems~\ref{thm:dn}--\ref{thm:rn} no
simultaneous diagonalisation of $\SF$ and $\SN$ is used, so the bound
applies verbatim to the mismatched (square-boundary) geometry. \qed

\paragraph{Level 2: the floor is controlled by the training loss.}
The residual $\epsN$ is bound to a measured quantity by a standard
a-posteriori elliptic estimate. The PINN Neumann/Robin loss is
\begin{equation}
\Ll_{\tagN}(\theta)=\norm{-\Delta\uN-f}_{L^2(\OmegaN)}^2
+\norm{\dn\uN+g}_{L^2(\SigmaI)}^2 \;(+\text{ Dirichlet term}),
\end{equation}
i.e.\ the squared residuals of the strong-form problem solved by
$\uN^{\mathrm{exact}}$. Writing $w=\uN-\uN^{\mathrm{exact}}$, the
difference solves the same elliptic problem with data equal to those
residuals; coercivity of the elliptic operator (a G\aa rding/Lax--Milgram
estimate with constant $C_{\mathrm{stab}}$ governed by the inverse of the
smallest resolvable Steklov eigenvalue) gives
$\norm{w}_{H^1(\OmegaN)}\le C_{\mathrm{stab}}\,\Ll_{\tagN}(\theta)^{1/2}$,
and the trace inequality
$\norm{w|_{\SigmaI}}_{H^{1/2}(\SigmaI)}\le C_{\mathrm{tr}}\norm{w}_{H^1(\OmegaN)}$
yields
\begin{equation}
\epsN:=\norm{\lambda_{\tagN}^k-\lambda_{\tagN}^{\mathrm{exact}}}_{H^{1/2}(\SigmaI)}
\le C_{\mathrm{tr}}C_{\mathrm{stab}}\,\Ll_{\tagN}(\theta^k)^{1/2}.
\label{eq:loss-floor}
\end{equation}
Combining~\eqref{eq:loss-floor} with~\eqref{eq:pinn-floor}, the plateau
height obeys
$\limsup_k\norm{e^k}\le C_\alpha C_{\mathrm{tr}}C_{\mathrm{stab}}\,
\Ll_{\tagN}^{1/2}/(1-\rhohat)=O(\Ll_{\tagN}^{1/2})$: the floor is set by
how well the network is trained. This is the falsifiable form of the
theorem---halving $\Ll_{\tagN}^{1/2}$ should halve the floor---and is what
the training-budget sweep of Appendix~\ref{app:E25} confirms (log--log
slope $1.13$, $R^2=0.89$), with the caveat documented there that a
cold-restart variant of the iteration exposes a loss-independent floor
component. Bounding the optimisation part of
$\Ll_{\tagN}$ further by an NTK gradient-descent rate (so that $\epsN$ is
expressed through the network's neural-tangent spectrum and a step count)
is a separate, lazy-training-dependent argument; it is not needed for the
floor to appear and we leave it to future work.

\section{Validating the loss-controlled floor: a training-budget sweep (E25)}
\label{app:E25}

Remark~\ref{rem:loss-floor} makes the one directly falsifiable
prediction of Theorem~\ref{thm:pinn-floor}: the RN plateau height should
scale as $O(\Ll_{\tagN}^{1/2})$ in the \emph{achieved} training loss.
We test it on the 2D natural orientation of
Section~\ref{sec:exp:2d} ($n_{\mathrm{iface}}=32$, $h=0.06$, the
width-$64$, depth-$5$ network of Appendix~\ref{app:inner-pinn}), running
the RN iteration at the searched optimum $\alpha=10.31$ of
Appendix~\ref{app:inner-pinn:e10} so that the geometric transient of
Equation~\eqref{eq:pinn-floor} clears well inside the
$\Kouter=16$ budget and the plateau is exposed; the theorem applies for
any $\alpha$ in the contractive band, the choice moving only the
constant $C_\alpha/(1-\rhohat)$. The inner Adam budget per outer
iteration is swept over $\Tinner\in\{10,25,50,100,250,500,1000,2000,
4000\}$ with three seeds each; the \emph{floor} is the median relative
trace error $\norm{\lambda^k-\lambda^\star}/\norm{\lambda^\star}$ over
the last six outer iterations and $\Ll_{\tagN}$ the median achieved loss
over the same window. Crucially, the regression is of the measured floor
against the measured loss, not against the budget: under the
warm-started protocol the paper actually uses, the budget$\to$loss
mapping is non-monotone (cumulative training across outer iterations and
the Adam minibatch noise floor saturate the achievable loss above
$\Tinner\approx100$ at this learning rate), but the prediction of
Remark~\ref{rem:loss-floor} concerns the loss, however it was reached.

\begin{figure}[h]
\centering
\figmaybe[0.62\linewidth]{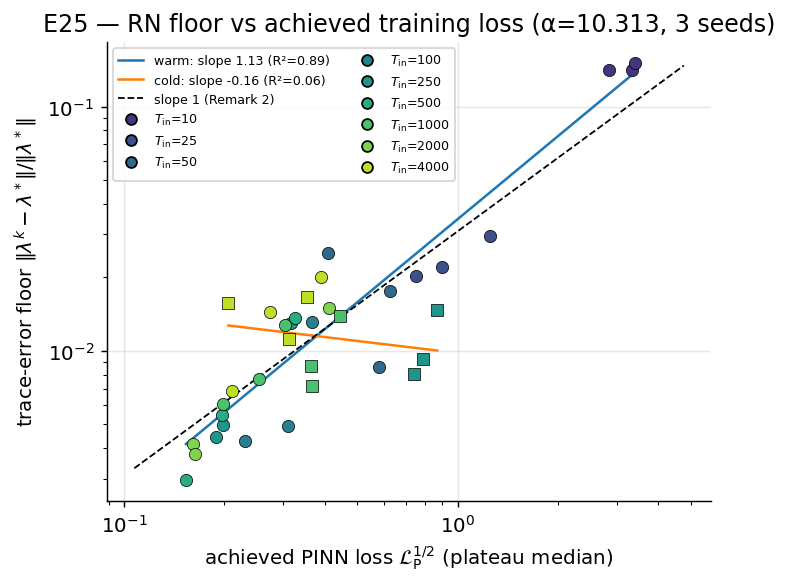}{E25 floor vs achieved loss, log--log, with slope-1 guide.}
\caption{E25: the RN trace-error floor against the achieved training
loss $\Ll_{\tagN}^{1/2}$ (log--log; circles: warm-started arm, colour =
inner budget; squares: cold-restart arm). Across a $400\times$ budget
range and $1.5$ decades of floor the warm-started runs follow the
slope-$1$ line of Remark~\ref{rem:loss-floor} (fitted slope $1.13$,
$R^2=0.89$, $27$ runs); the cold-restart arm saturates at a
loss-independent ${\approx}10^{-2}$ (fitted slope $-0.16$).}
\label{fig:e25-floor}
\end{figure}

\paragraph{The prediction holds in the warm-started protocol.} Across
$27$ runs spanning $1.5$ decades of floor
(Figure~\ref{fig:e25-floor}), the fitted log--log slope is $1.13$
against the predicted $1$, with $R^2=0.89$: halving the achieved
$\Ll_{\tagN}^{1/2}$ does, to within the seed scatter, halve the floor.
The seed-to-seed spread reinforces rather than dilutes the mechanism: at
a fixed budget, the seed with the larger achieved loss lands on the
larger floor, on the same line.

\paragraph{A caveat from the cold-restart arm.} A control arm
that retrains a \emph{fresh} network at every outer iteration (the
closest realisation of the theorem's per-step retraining hypothesis;
$\Tinner\in\{250,1000,4000\}$, three seeds) does \emph{not} show the
proportionality: its floor sits at ${\approx}10^{-2}$ essentially
independent of its (much larger and budget-sensitive) loss. This is
consistent with the structure of Theorem~\ref{thm:pinn-floor} rather
than a violation of it: the bound~\eqref{eq:pinn-floor} controls the
plateau by the \emph{sum} of the additive-residual term and the
operator-perturbation term $\EN(\alpha I+\SF)^{-1}$, and under cold
restarts each step's freshly re-randomised network bias re-enters
through the latter, which does not shrink with the per-step loss. The
proportional, loss-controlled reading of Remark~\ref{rem:loss-floor} is
therefore a property of the warm-started iteration the paper uses, in
which the network's systematic (cap) error is held fixed across outer
steps and the plateau is dominated by the training residual.

\section{Fourier Steklov probe: implementation details}
\label{app:sn-probe}

\begin{algorithm}[h]
\caption{\textsc{Fourier Steklov probe for a PINN subdomain}}
\label{alg:sn-probe}
\begin{algorithmic}[1]
\Require PINN $\NN$, interface $\SigmaI$ parameterised by the angle
$\theta\in[0,2\pi)$, modes $k=0,\ldots,K_{\max}$, budget $\Tinner$
\For{$k = 0$ to $K_{\max}$}
  \For{$\phi \in \{\cos, \sin\}$ (skip $\sin$ for $k=0$)}
    \State Train $\NN$ for $\Tinner$ steps to minimise the
    homogeneous PDE residual subject to soft Dirichlet
    $\uN|_\SigmaI=\phi(k\theta)$
    \State Read the radial flux $g_k^\phi(\theta)=\dnF\uN(\theta)$ on
    $\SigmaI$ by autograd
    \State Project $\hat\sigma_{\tagN,k}^\phi \gets \langle g_k^\phi,\phi(k\cdot)\rangle/\langle\phi(k\cdot),\phi(k\cdot)\rangle$
  \EndFor
\EndFor
\State \Return $\{\hat\sigmaN(k)\}_{k=0}^{K_{\max}}$
\end{algorithmic}
\end{algorithm}

\begin{figure}[h]
\centering
\begin{subfigure}[b]{0.49\linewidth}
  \centering
  \figmaybe{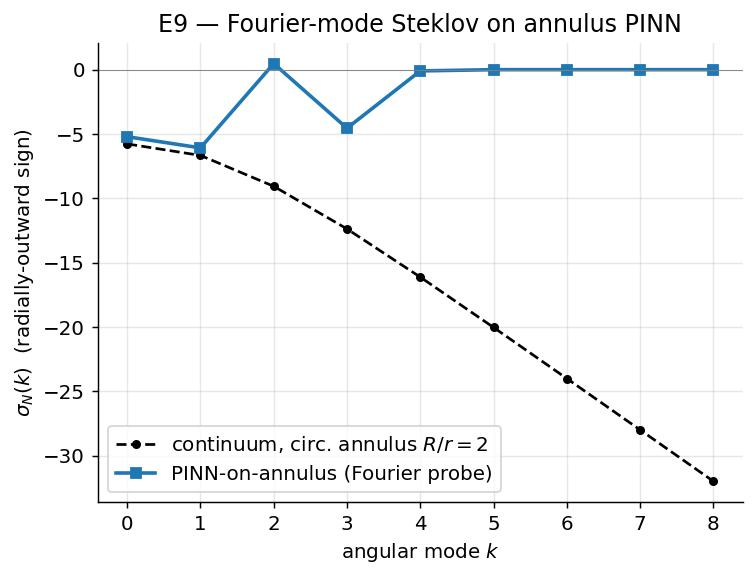}{Fourier-mode Steklov probe on the annulus PINN.}
  \caption{Fourier-mode estimate $\hat\sigmaN(k)$
  vs.\ continuum~\eqref{eq:sn-annulus}: $k=0,1$ matched within
  ${\sim}10\%$; $k\geq2$ collapses to ${\sim}0$ (the PINN resolution
  cap).}
\end{subfigure}\hfill
\begin{subfigure}[b]{0.49\linewidth}
  \centering
  \figmaybe{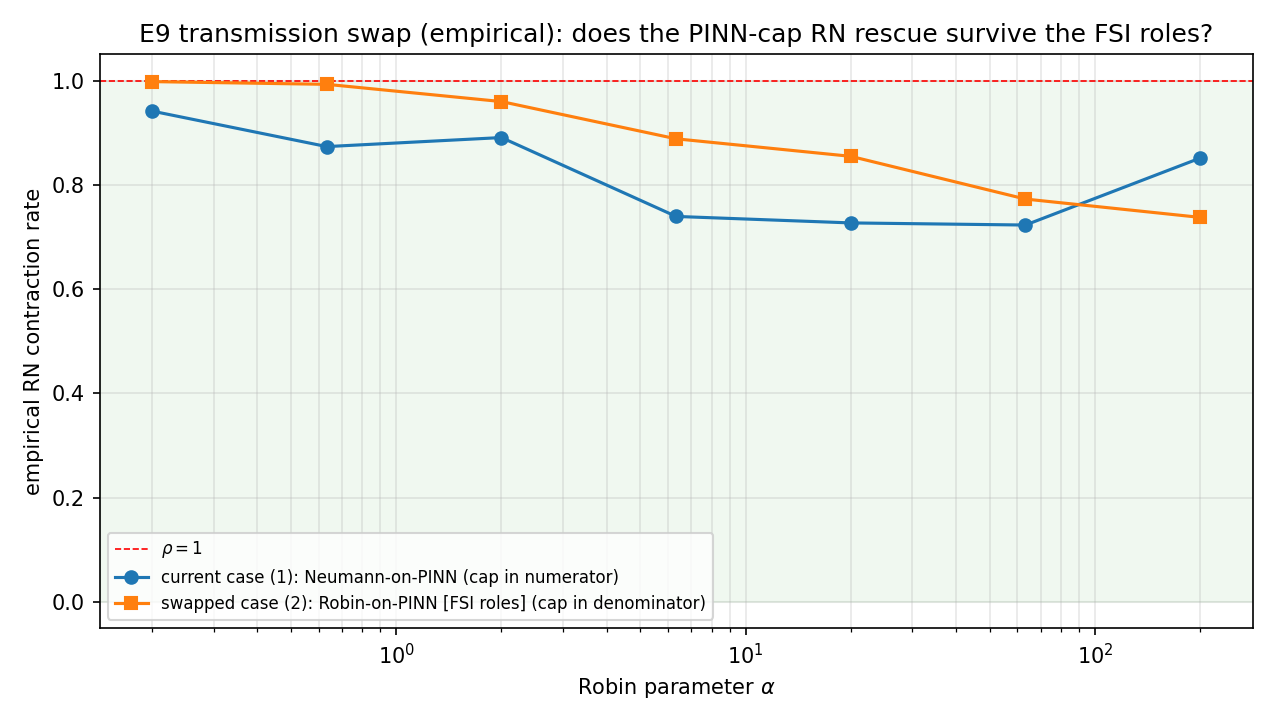}{Empirical RN contraction rate vs.\ Robin parameter for the two transmission assignments.}
  \caption{E9 transmission-role swap (single seed): empirical RN rate
  vs.\ $\alpha$ for case~(1) (Neumann-on-PINN) and the swapped case~(2)
  used by FSI (Robin-on-PINN). Both contractive across $[0.2,200]$ with
  near-identical best rates ($0.72$ vs.\ $0.74$).}
\end{subfigure}
\caption{FSI-aligned orientation: the annulus PINN's spectral cap (left)
and the role-invariance of the RN contraction under the FSI transmission
swap (right).}
\label{fig:e9-fourier-swap}
\end{figure}

The probe of Section~\ref{sec:theory:sn-estimator}, stated as
Algorithm~\ref{alg:sn-probe}, is implemented in two
variants, one for the disc PINN (FEM on the annulus) and one for the
annulus PINN (FEM on the disc; FSI-aligned). Each takes a user-supplied
trainer that returns a trained PINN given the boundary data and a flux
reader that returns the boundary flux $g(\theta)$, so the probe is
decoupled from any particular PINN architecture; the continuum
references are $\sigmaN(k)=k/R$ for the disc and
\eqref{eq:sn-annulus} for the annulus. For Section~\ref{sec:exp:2d} we
use $K_{\max}=8$, $\Tinner=2500$ Adam steps per probe, MLP width $=64$,
depth $=5$, soft Dirichlet weight $\mu_b=100$, soft PDE residual weight
$\mu_r=1$; wall time ${\sim}142$ s on a 4-thread Apple-silicon CPU.
Figure~\ref{fig:e9-fourier-swap} reports the resulting annulus-PINN
Steklov estimate $\hat\sigmaN(k)$ against the
continuum~\eqref{eq:sn-annulus} (left), exposing the spectral cap at
$k\ge2$, together with the empirical RN contraction rate under both
transmission-role assignments (right), the data behind the role-invariance
claim of Section~\ref{sec:exp:2d}.

\begin{figure}[h]
\centering
\figmaybe[0.62\linewidth]{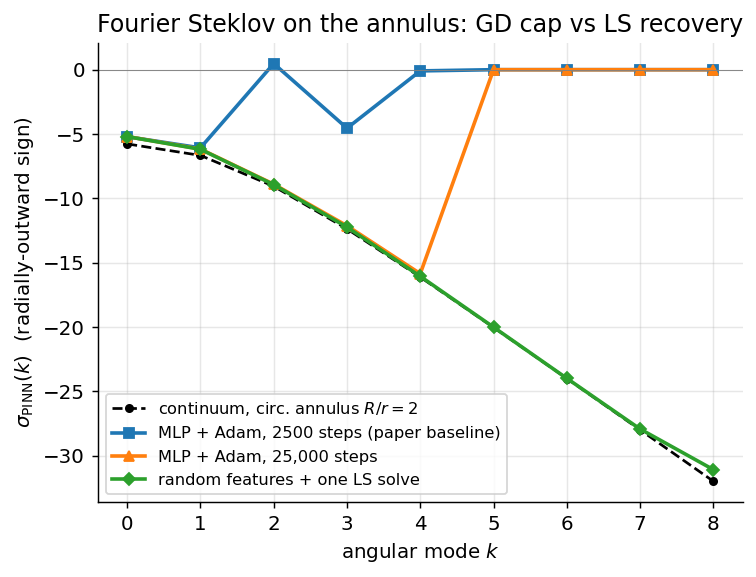}{Fourier
Steklov probe on the annulus: gradient-descent cap vs least-squares
recovery.}
\caption{The spectral cap of Figure~\ref{fig:e9-fourier-swap}(a) is an
optimisation effect: under the paper's training budget ($2500$ Adam
steps) the probe collapses at $k\geq2$; a tenfold budget recedes the
cap to $k=5$ but no further; a random-feature model fitted by a single
least-squares solve, through the same probe, recovers all modes
$k\leq8$ to the continuum reference.}
\label{fig:e9-cap-recovery}
\end{figure}

The collapse in Figure~\ref{fig:e9-fourier-swap}(a) is removable, and
its removal in this scalar geometry mirrors the Stokes-side findings of
Appendices~\ref{app:E22} and~\ref{app:E23} mode for mode.
Figure~\ref{fig:e9-cap-recovery} repeats the probe with two alternative
inner solvers at otherwise identical settings (same interface points,
projection, and sign conventions). Staying with gradient descent but
raising the budget tenfold ($\Tinner=25{,}000$ Adam steps per mode)
recedes the cap from $k=2$ to $k=5$: modes $k=2,3,4$ land within
$1$--$2\%$ of the continuum, after which the estimate collapses to
numerical zero, the scalar analogue of the budget-helps-only-to-a-plateau
behaviour of Remark~\ref{rem:spectral-cap}. Replacing the
gradient-trained network altogether with a random-feature model
($M=3200$ frozen random Fourier features of frequency scale $32$, double
precision) whose linear head is fitted by a single column-equilibrated
least-squares solve against the same homogeneous constraints (the PDE
residual at $4096$ annulus collocation points, the prescribed trace at
the $32$ interface points, and the homogeneous Dirichlet condition at
$1024$ outer-boundary points) recovers \emph{every} probed mode, with
relative error $0.1$--$2.8\%$ for $k\geq2$; the residual $7$--$10\%$
offsets at $k=0,1$ are shared by all arms and reflect the square outer
boundary against the circular-annulus reference, not the solver. As on
the Stokes probe, feature count and frequency scale must be raised
together: $M=1600$ at scale $16$ reaches only $k=4$, with a
rank-deficient least-squares system on the higher modes.

\section{Static-Stokes vector recipe: hyperparameters and ablations}
\label{app:E11-hyper}

\begin{figure}[h]
\centering
\figmaybe[0.8\linewidth]{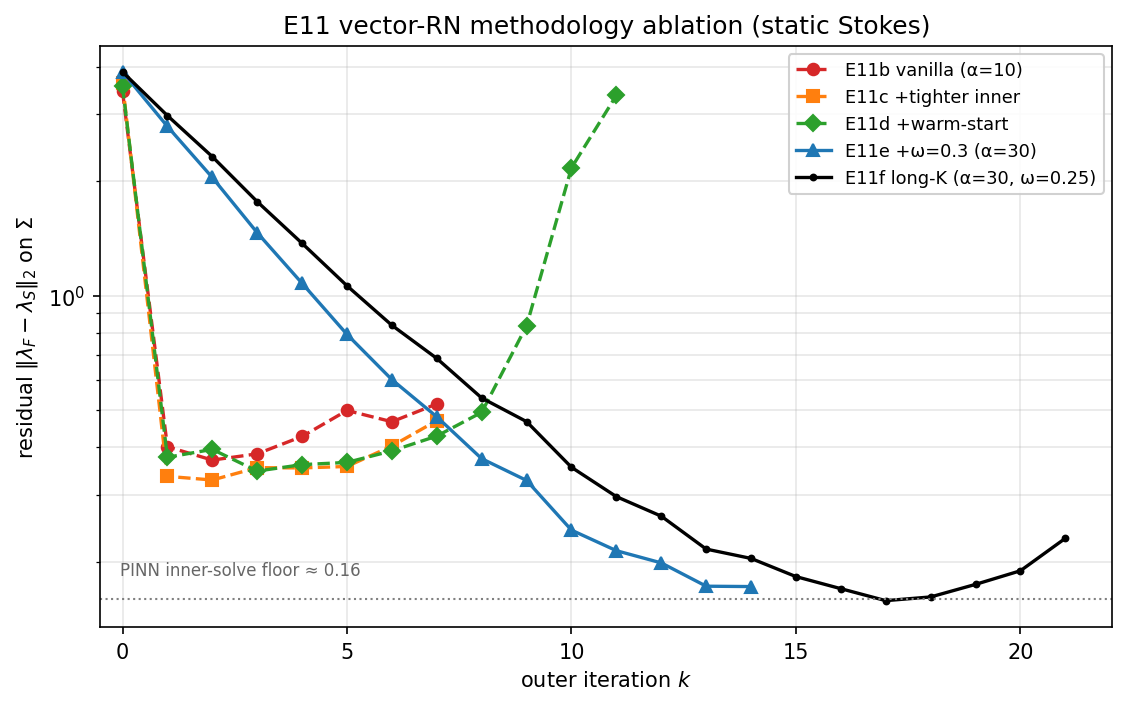}{Composite residual-vs-outer-iter plot for E11b/c/d/e/f.}
\caption{Methodology ablation on the static-Stokes vector RN iteration.
Each curve adds one ingredient of Section~\ref{sec:fsi:recipe}: (b)
vanilla cold-restart $\alpha=10$, (c) tightened inner solve, (d)
warm-start at $\alpha=10$ exposing the underlying RN-map divergence the
cold restart had masked, (e) warm-start $+$ $\omega=0.30$ at $\alpha=30$
giving the monotone contraction, (f) long-$K$ extension at $\omega=0.25$
showing the late-time warm-start drift that motivates the early-stop
ingredient. The PINN floor ${\approx}0.16$ is annotated.}
\label{fig:e11-ablation}
\end{figure}

The recipe of Section~\ref{sec:fsi:recipe} uses: PINN width $=32$, depth
$=4$, Tanh, no Fourier features (a Fourier-features ablation with
embedding frequencies $\{1,2,4,8\}$ at the same width/depth and
$\Tinner=2000$ does \emph{not} raise the spectral cap: the cap is
capacity-limited, not basis-limited, see
Figure~\ref{fig:e11-ablation}); $4000$ Adam
steps cold-start (once), $1500$ refine per outer iteration; loss weights
$w_{\mathrm{neu}}=200$, $w_{\mathrm{div}}=50$, $w_{\mathrm{wall}}=100$,
$w_{\mathrm{robin}}=200$, $w_{\mathrm{pres}}=5$; RN parameters
$\alpha=30$, $\omega=0.30$, $\Kouter=22$, early-stop growth factor
$1.05$. Wall time per outer iteration ${\sim}2.6$ s on CPU; the combined
recipe reaches residual floor $0.158$ at $k=17$ and relative trace error
$0.136$.

\begin{figure}[h]
\centering
\begin{subfigure}[b]{0.49\linewidth}
  \centering
  \figmaybe{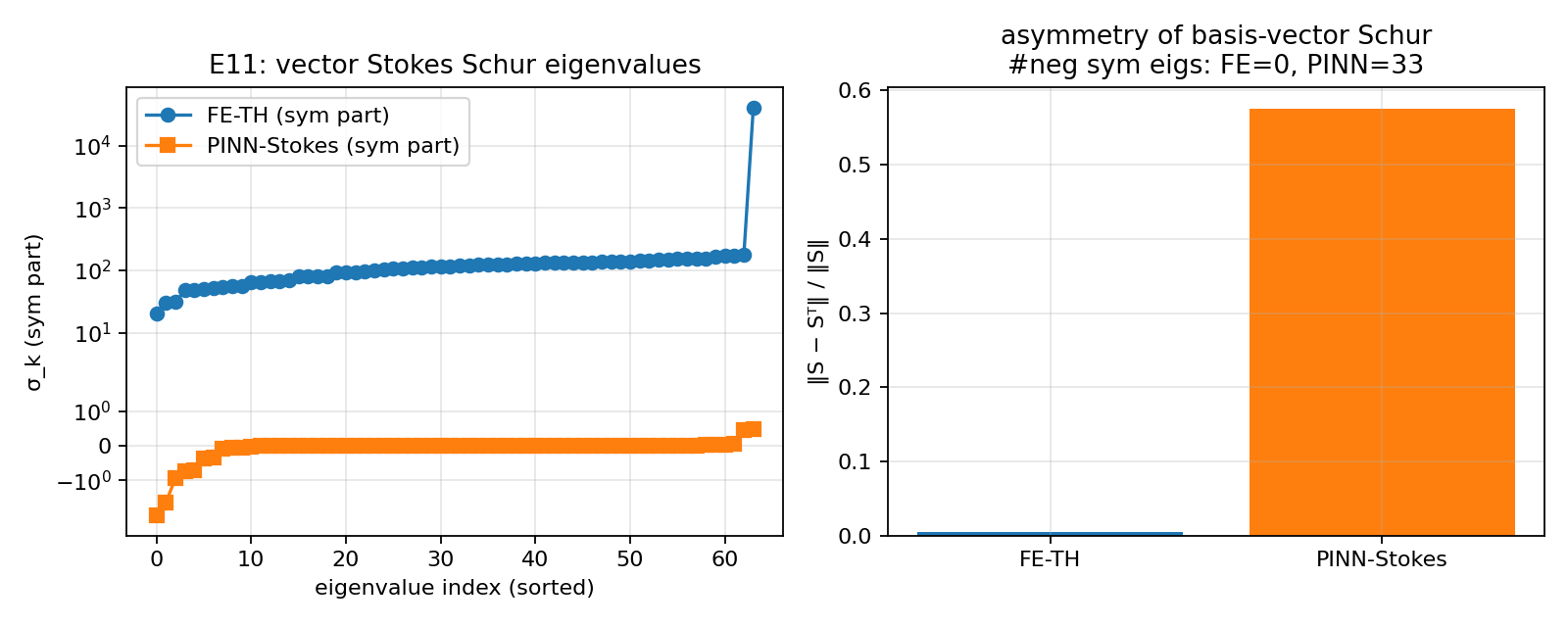}{E11 spectra of $\SF$ and $\SN$.}
  \caption{Eigenvalues of $\bSF$ (FE-TH) vs.\ $\bSN$ (PINN-Stokes),
  $M$-normalised: the PINN side is capped.}
\end{subfigure}\hfill
\begin{subfigure}[b]{0.49\linewidth}
  \centering
  \figmaybe{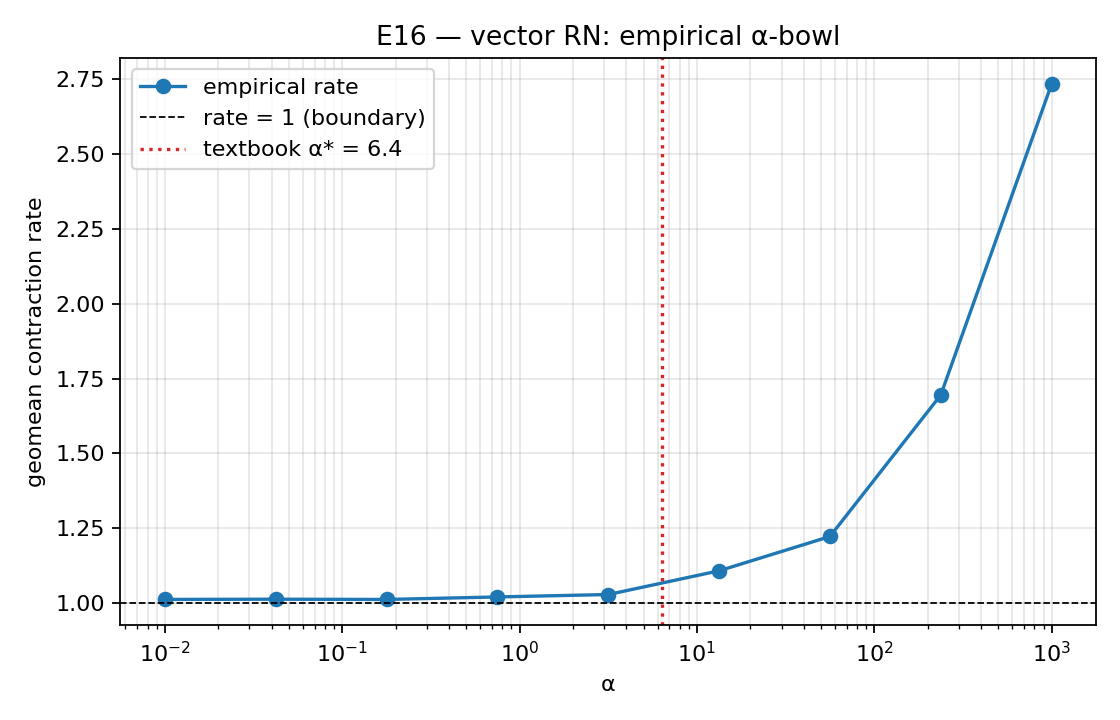}{E16: numerical $\alpha$-bowl for vector RN.}
  \caption{E16: vanilla vector RN gives no $\alpha$ with rate $<1$ across
  $[10^{-2},10^3]$.}
\end{subfigure}
\caption{Static-Stokes vector RN diagnostic (E11, E16).}
\label{fig:e11-spectra}
\end{figure}

The FEM vector Schur $\bSF$ ($64\times64$) is assembled by $64$
Taylor--Hood Stokes solves on the disc; the FEM spectrum spans
$[20,4\times10^4]$ while the PINN side is capped at
$\sigmaN^{\max}\approx0.5$ (basis probe) / $3.06$ at $k=2$ (Fourier
probe). The per-mode $2\times2$ block asymmetry drops from $\le135\%$
(basis-vector probe) to $\le90\%$ (Fourier probe), so part of the
original $58\%$ value was probe-method noise and part is a genuine
non-symmetry of the discrete Steklov. A numerical $\alpha$-search (E16)
on $[10^{-2},10^3]$ finds no contractive $\alpha$ with the vanilla
iteration, which is exactly the regime in which the four-ingredient
recipe restores transient contraction. Figure~\ref{fig:e11-spectra} plots
both spectra side by side (the capped PINN $\bSN$ against the FE
Taylor--Hood $\bSF$) and the flat E16 $\alpha$-bowl that has no point
below $1$, and Figure~\ref{fig:e11-ablation} traces the recipe ingredient
by ingredient.

\section{Large-added-mass two-slab configuration (Section~\ref{sec:exp:addedmass})}
\label{app:addedmass-config}

\begin{figure}[h]
\centering
\figmaybe[0.95\linewidth]{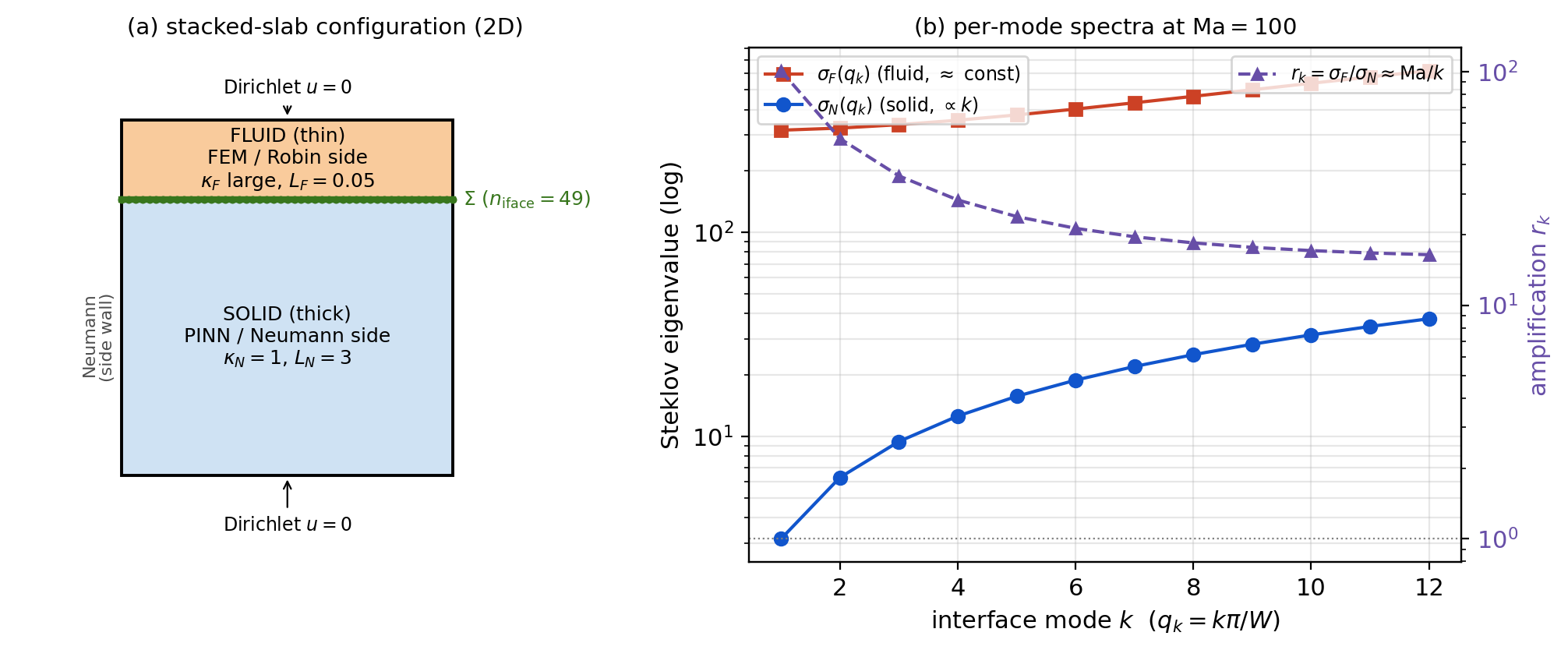}{Two-slab
added-mass configuration and its per-mode Steklov spectra.}
\caption{Configuration for the large-added-mass experiments of
Section~\ref{sec:exp:addedmass}. (a) The 2D stacked-slab geometry: a thin
fluid slab (FEM, Robin side) over a thick solid slab (PINN, Neumann side)
sharing the interface $\SigmaI$; the thick solid is drawn compressed for
visibility. (b) Closed-form flat-slab Steklov spectra at $\mathrm{Ma}=100$:
$\sigmaF(q_k)$ is nearly flat while $\sigmaN(q_k)$ grows with $k$, so the
DN amplification $r_k\approx\mathrm{Ma}/k$ (right axis) is large at the low
modes and decays through $1$, which is why a single scalar relaxation
cannot straddle the spectrum and RN's per-mode impedance matching wins.}
\label{fig:addedmass-config}
\end{figure}

The added-mass experiments of Section~\ref{sec:exp:addedmass} are genuine
2D coupled iterations on the stacked-slab geometry of
Figure~\ref{fig:addedmass-config}(a): a thin high-diffusivity ``fluid''
slab $\OmegaF=[0,W]\times[0,L_F]$ on top of a thick ``solid'' slab
$\OmegaN=[0,W]\times[-L_N,0]$, sharing the flat interface
$\SigmaI=\{y=0\}$ with $n_{\mathrm{iface}}=49$ nodes. The far faces carry
homogeneous Dirichlet data and the side walls carry natural Neumann
conditions, so the interface modes are $q_k=k\pi/W$. The role assignment
follows the Poisson convention of Section~\ref{sec:problem}
(Robin/Dirichlet-on-FEM, Neumann-on-PINN), with the FEM on the thin fluid
(the added-mass side, where the impedance belongs) and the PINN on the
thick solid. Both slabs are discretised by tensor-product $Q_1$ elements,
for which the flat-slab Steklov eigenvalue is the closed form
$\sigma(q)=\kappa\,q\coth(qL)$. The thin fluid gives a nearly
mode-independent $\sigmaF(q_k)\approx\kappa_F/L_F$ (the lumped added-mass
limit) and the thick solid a growing $\sigmaN(q_k)\approx\kappa_N q_k$, so
the per-mode amplification is $r_k=\sigmaF/\sigmaN\approx\mathrm{Ma}/k$
with the single dial $\mathrm{Ma}=\sigmaF(q_1)/\sigmaN(q_1)$
(Figure~\ref{fig:addedmass-config}b). The geometric parameters are fixed
at $W=1$, $L_F=0.05$, $L_N=3$, $\kappa_N=1$, and $\kappa_F$ is tuned per
target $\mathrm{Ma}$ from the first-mode ratio; the PINN solid sits on the
smooth, low-mode-dominated side, exactly the band a network resolves, so
its spectral cap truncates only the (already stable) high modes.

\paragraph{Governing equations.} Each slab carries a scalar steady
diffusion (Poisson--Laplace) problem,
\begin{equation}
  -\nabla\!\cdot\!\bigl(\kappa_i\,\nabla u_i\bigr)
  \;=\; -\,\kappa_i\,\Delta u_i \;=\; 0
  \quad\text{in } \Omega_i,\qquad i\in\{\tagF,\tagN\},
  \label{eq:addedmass-pde}
\end{equation}
coupled across the shared interface $\SigmaI=\{y=0\}$ by the standard
transmission conditions, continuity of the field and of the normal flux,
\begin{equation}
  u_\tagF = u_\tagN
  \quad\text{and}\quad
  \kappa_\tagF\,\partial_n u_\tagF
  = \kappa_\tagN\,\partial_n u_\tagN
  \qquad\text{on } \SigmaI,
  \label{eq:addedmass-transmission}
\end{equation}
with homogeneous Dirichlet data on the far faces $y=L_\tagF$ and
$y=-L_\tagN$ and homogeneous Neumann data on the side walls
$x\in\{0,W\}$.

\paragraph{Why these interface modes.} The lateral structure of the
problem is set entirely by the side walls. Separating variables in
\eqref{eq:addedmass-pde}, $u_i(x,y)=X(x)Y(y)$, gives
$X''/X=-Y''/Y=-q^2$, so the admissible lateral profiles solve
$X''+q^2X=0$ subject to the side-wall conditions. With homogeneous
\emph{Neumann} walls, $\partial_x u=0$ at $x=0$ and $x=W$, the surviving
solutions are the cosines $X(x)=\cos(q x)$ whose derivative
$-q\sin(qx)$ vanishes at \emph{both} ends: $\sin(q\cdot 0)=0$ holds
automatically, and $\sin(qW)=0$ forces $qW=k\pi$, i.e.
\begin{equation}
  q_k=\frac{k\pi}{W},\qquad k=0,1,2,\dots
  \label{eq:addedmass-qk}
\end{equation}
These $\{\cos(q_k x)\}$ are exactly the eigenfunctions of $-\partial_x^2$
on $[0,W]$ under Neumann boundary conditions (eigenvalues $q_k^2$); they
are mutually $L^2(0,W)$-orthogonal and complete, so they form an
\emph{interface eigenbasis} in which every trace on $\SigmaI$ expands
uniquely. The index $k$ counts lateral half-periods across the width
(wavelength $2W/k$): $k=0$ is the uniform mode and larger $k$ is finer
lateral structure. (Had the side walls been Dirichlet instead, the same
$q_k$ would appear with sines $\sin(q_k x)$; the Neumann choice is what
makes the constant $k=0$ mode admissible and the basis pure cosines.)
This is the flat-interface analogue of the angular Fourier modes
$e^{\mathrm{i}k\theta}$ used on the circular $\SigmaI$ elsewhere in the
paper.

\paragraph{Per-mode Steklov eigenvalue.} Because the modes decouple,
writing $u_i(x,y)=\sum_k a_{i,k}(y)\cos(q_k x)$ in
\eqref{eq:addedmass-pde} leaves each amplitude obeying
$a_{i,k}''=q_k^2\,a_{i,k}$, so the harmonic extension of a single
interface mode decays hyperbolically ($\cosh/\sinh$) across a slab and
the slab's Dirichlet-to-Neumann (Steklov) eigenvalue, mapping interface
trace to interface flux against the homogeneous far face, is exactly the
closed form $\sigma(q_k)=\kappa\,q_k\coth(q_k L)$ used above. The
operator-level analysis is thus an exact consequence of
\eqref{eq:addedmass-pde}--\eqref{eq:addedmass-transmission}; an independent
finite-element realisation of the same two-slab problem (included with the
reproducibility package) recovers these eigenvalues
to within the discretisation error.

Two studies are run. Study~A sweeps $\mathrm{Ma}\in[0.3,300]$ with
\emph{both} slabs by tensor-product $Q_1$ FEM (cheap; checks that the
operator-level crossover survives a genuine discrete 2D coupling) and
produces Figure~\ref{fig:addedmass}(a) and Table~\ref{tab:addedmass}.
Study~B fixes the representative $\mathrm{Ma}=100$ and replaces the solid
by a \emph{real trained PINN} (width/depth as in the FSI network, $400$
Adam steps per outer iteration, three seeds), running the five schemes
DN, DN${+}\omega^\ast$, DN${+}$Aitken, RN($\alphastar$), RN($\alpha$-search)
head-to-head against a monolithic-FEM reference; it produces
Figure~\ref{fig:addedmass}(b). In Study~B, RN($\alphastar$) converges to
the PINN floor ${\approx}5\times10^{-2}$ at the sweep-free geometric-mean
rate $0.45\pm0.004$ while DN and Aitken diverge or stall, confirming the
operator-level prediction in the presence of a trained subdomain solver.

\section{Alart--Curnier contact: closed-form one-step over-relaxation}
\label{app:contact}

The Alart--Curnier projection for normal contact between the disc and
the bottom wall $\Sigma_p$ reads
\begin{equation}
\lambda_c^+ = \max\bigl(0,\, \lambda_c - \omega_{\mathrm{AL}}\, \gamma_c\, \mathrm{gap}(y_c)\bigr),
\qquad
\mathrm{gap}(y_c) = y_c - R - \varepsilon_g,
\label{eq:al}
\end{equation}
with $\gamma_c>0$ the AL penalty stiffness and $\varepsilon_g$ the
relaxed gap. Substituting the backward-Euler position update
$y_c^{n+1}=y_c^{\mathrm{free}}+(\dt^2/m_s)\lambda_c^{n+1}$ into the active
branch gives a contraction factor
$\abs{1-\omega_{\mathrm{AL}}\gamma_c\dt^2/m_s}$; setting it to zero gives
$\omega_{\mathrm{AL}}^{\mathrm{opt}}=m_s/(\gamma_c\dt^2)$, under which the
projection converges in one step and lands the disc at
$y_c=R+\varepsilon_g$ exactly (machine precision), strictly improving on
the penalty-method $O(\gamma_c^{-1})$ bound. The $\max(\cdot,0)$ clamp
keeps $\lambda_c\ge0$, so over-relaxation beyond $1$ is safe.

Algorithm~\ref{alg:contact-fsi} assembles this contact projector, the
backward-Euler solid ODE, the vector-RN fluid--solid coupling of
Section~\ref{sec:fsi:model}, and the per-step collocation regeneration
into the full time-stepped contact-FSI solver referenced from
Section~\ref{sec:fsi:model}.

\begin{algorithm}[h]
\caption{\textsc{Time-stepped contact FSI} (backward-Euler outer step,
vector-RN fluid--solid inner coupling, Alart--Curnier AL contact). The
fluid is the PINN, the rigid-disc ODE the FEM side. As in
Algorithm~\ref{alg:dn-rn}, $(\lambda,g)$ denotes the RN
trace--traction pair; the unpaired symbol $g$ in the Require line and
the predictor step is gravity.}
\label{alg:contact-fsi}
\begin{algorithmic}[1]
\Require initial height/velocity $(y_c^0,\dot y_c^0)$; step $\dt$;
impedance $\alpha$, relaxation $\omega$, inner budget $\Kouter$, tol;
AL stiffness $\gamma_c$, gap $\varepsilon_g$, over-relaxation
$\omega_{\mathrm{AL}}=m_s/(\gamma_c\dt^2)$; disc radius $R$, mass $m_s$,
gravity $g$.
\For{$n = 0, 1, \ldots$}
  \State \textbf{Regenerate collocation.} Sample the current fluid domain
  $\Omegaf=\Omega_+\setminus\bigl(B((0,y_c^n),R)\cup\{y<\varepsilon_g\}\bigr)$;
  warm-start the PINN from step $n$ (carry weights and Adam moments).
  \State \textbf{Predictor.} Backward-Euler free (contact-off) height
  $y_c^{\mathrm{free}}\gets y_c^n+\dt\,\dot y_c^n - \dt^2 g$;
  initialise $\lambda_c\gets0$ and the interface data $(\lambda,g)$ from step $n$.
  \For{$k = 0, 1, \ldots, \Kouter-1$}
    \State \textbf{Fluid (PINN) solve.} Train the $(u,p)$ network on
    $\Omegaf$ enforcing steady Stokes, $\nabla\!\cdot u=0$, no-slip on the
    rigid walls, Navier-slip on $\Sigma_p$, and the Robin condition
    $\sigma_f n+\alpha u=\alpha\lambda+g$ on $\SigmaI$.
    \State \textbf{Read traction.} $g_F\gets(\sigma_f n)|_\SigmaI$ via the
    pointwise transfer~\eqref{eq:units-fix};
    $F_y\gets\int_\SigmaI(\sigma_f n)\cdot e_y\,ds$.
    \State \textbf{Solid (FEM) + contact.} Backward-Euler update
    $y_c\gets y_c^{\mathrm{free}}+(\dt^2/m_s)\lambda_c$ with the
    Alart--Curnier projection
    $\lambda_c\gets\max\!\bigl(0,\lambda_c-\omega_{\mathrm{AL}}\gamma_c\,(y_c-R-\varepsilon_g)\bigr)$;
    set $\dot y_c\gets(y_c-y_c^n)/\dt$ and the contact reaction
    $F_c\gets\lambda_c$.
    \State \textbf{Kinematic trace.} $\lambda_N\gets\dot y_c\,e_y$ on
    $\SigmaI$ (the velocity the disc imposes on the fluid).
    \State \textbf{vector-RN update.}
    $(\lambda,g)\gets(1-\omega)(\lambda,g)+\omega(\lambda_N,g_F)$.
    \If{residual increases for the first time \textbf{or}
    $\norm{\lambda_N-\lambda}<\mathrm{tol}$} \textbf{break} \EndIf
  \EndFor
  \State \textbf{Accept step.} $(y_c^{n+1},\dot y_c^{n+1},\lambda_c^{n+1})\gets(y_c,\dot y_c,\lambda_c)$.
\EndFor
\end{algorithmic}
\end{algorithm}

\section{Coupling baselines: RN versus the relaxation ladder}
\label{app:baselines}

\begin{figure}[h]
\centering
\begin{subfigure}[b]{0.49\linewidth}
  \centering
  \figmaybe{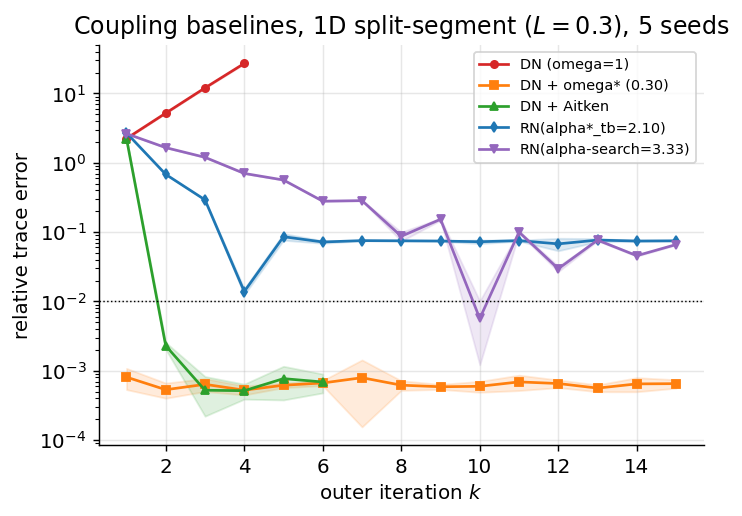}{1D relaxation ladder.}
  \caption{1D split-segment (single mode): optimal constant relaxation
  annihilates the one mode in $1$--$2$ steps; RN sits at the PINN floor.}
\end{subfigure}\hfill
\begin{subfigure}[b]{0.49\linewidth}
  \centering
  \figmaybe{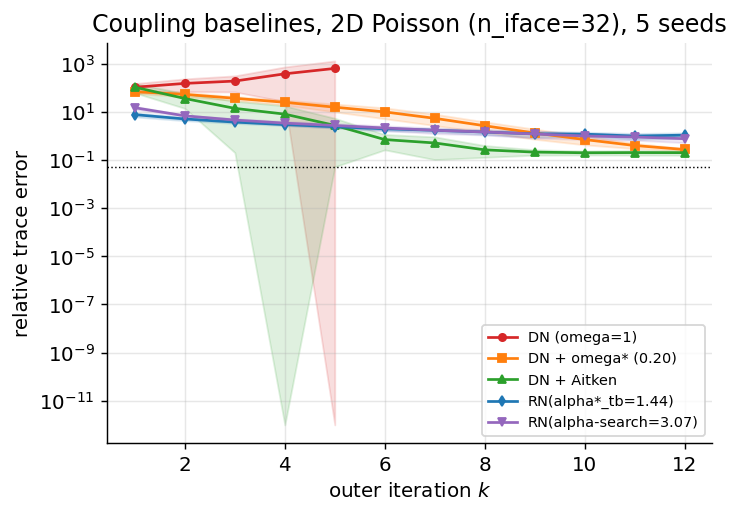}{2D relaxation ladder.}
  \caption{2D disc-in-square (multi-mode): DN diverges; all cures
  contract but tuned relaxation beats RN with non-overlapping bands.}
\end{subfigure}
\caption{RN versus the classical relaxation ladder at low added mass.
Relative trace error per outer iteration; band is $\pm$std over five
seeds.}
\label{fig:baselines}
\end{figure}

This appendix gives the full low-added-mass relaxation-ladder benchmark
summarised in Section~\ref{sec:exp:addedmass}. We benchmark RN against
plain DN ($\omega=1$), DN with the best constant under-relaxation
$\omega^\ast$, DN with Aitken/Irons--Tuck dynamic relaxation, and RN at
both the closed-form $\alphastar$ and a numerically searched optimum,
over five seeds at the operating points of
Sections~\ref{sec:exp:poisson}--\ref{sec:exp:2d}. (The quasi-Newton
interface accelerators of partitioned FSI, notably
IQN-ILS~\citep{DegrooteBatheVierendeels2009}, are the natural next rung
of this ladder; we leave them as the stronger history-based baseline
for a future revision.) In 1D the interface carries a
single mode, so the optimal constant $\omega^\ast=1/(1-g)$ annihilates
the error in essentially one step and RN merely settles at the PINN
floor. In the multi-mode 2D disc-in-square a single scalar relaxation
cannot be optimal for every mode, yet on this \emph{low-added-mass}
problem it nonetheless \emph{beats} RN (final error $0.20$--$0.27$ vs.\
RN's $0.76$--$1.1$, non-overlapping five-seed bars), and RN's
$\alpha$-search optimum is itself seed-unstable
(Table~\ref{tab:baselines}, Figure~\ref{fig:baselines}). RN does not win
on iteration count here; its advantage is being \emph{sweep-free}: a
single closed-form impedance, no $\omega$-sweep and no per-iteration
history. The premium for that convenience is visible in the table: rate
$0.84\pm0.02$ at the closed-form $\alphastar$ vs.\ $0.76\pm0.03$ after
an $\alpha$-sweep. As Section~\ref{sec:exp:addedmass} shows, this
ranking inverts
in the large-added-mass regime that motivates FSI.

\begin{table}[h]
\centering
\caption{Coupling baselines: RN versus the relaxation ladder
(mean$\,\pm\,$std over five seeds). ``final err'' is the relative trace
error after the outer budget ($\Kouter=15$ in 1D, $12$ in 2D); ``2D
rate'' is the geometric-mean per-step contraction. Plain DN diverges on
every seed; every cure bounds the iteration, but in the multi-mode 2D
case tuned scalar relaxation beats RN with non-overlapping bars. The 1D
std on the relaxation rows is $<\!10^{-4}$ and is omitted.}
\label{tab:baselines}
\small
\begin{tabular}{@{}lrrrl@{}}
\toprule
scheme & 1D final err & 2D final err & 2D rate & tuning required \\
\midrule
DN ($\omega = 1$)            & diverges & diverges & $1.5{\pm}0.2$ & none \\
DN $+\,\omega^*$             & $6.5{\times}10^{-4}$ & $0.27{\pm}0.03$ & $0.61{\pm}0.01$ & 1-D $\omega$ sweep \\
DN $+$ Aitken                & $6.1{\times}10^{-4}$ & $0.20{\pm}0.04$ & $0.57{\pm}0.01$ & per-iter.\ history \\
RN($\alphastar$, closed form) & $7.5{\times}10^{-2}$ & $1.1{\pm}0.3$ & $0.84{\pm}0.02$ & closed-form $\alpha$ \\
RN($\alpha$-search)          & $6.6{\times}10^{-2}$ & $0.76{\pm}0.24$ & $0.76{\pm}0.03$ & $\alpha$ sweep \\
\bottomrule
\end{tabular}
\end{table}

\section{What cures the warm-start drift (E21)}
\label{app:E21}

We turn the cure for the E12 warm-start drift into the independent
variable. We re-run the E12 free-fall and, after every time step,
evaluate an independent Stokes momentum-plus-continuity residual on a
fresh collocation set that never enters training; its growth factor
$G:=\operatorname{median}(\text{tail})/\operatorname{median}(\text{settled
floor})$ is the headline drift metric. A single knob selects the
strategy: periodic \emph{cold} restart, periodic \emph{Adam-reset},
\emph{strong} (larger inner budget), \emph{reproject} (keep training
until the monitor recovers), and \emph{hybrid} (residual-triggered cold
restart).

\begin{figure}[h]
\centering
\figmaybe[0.95\linewidth]{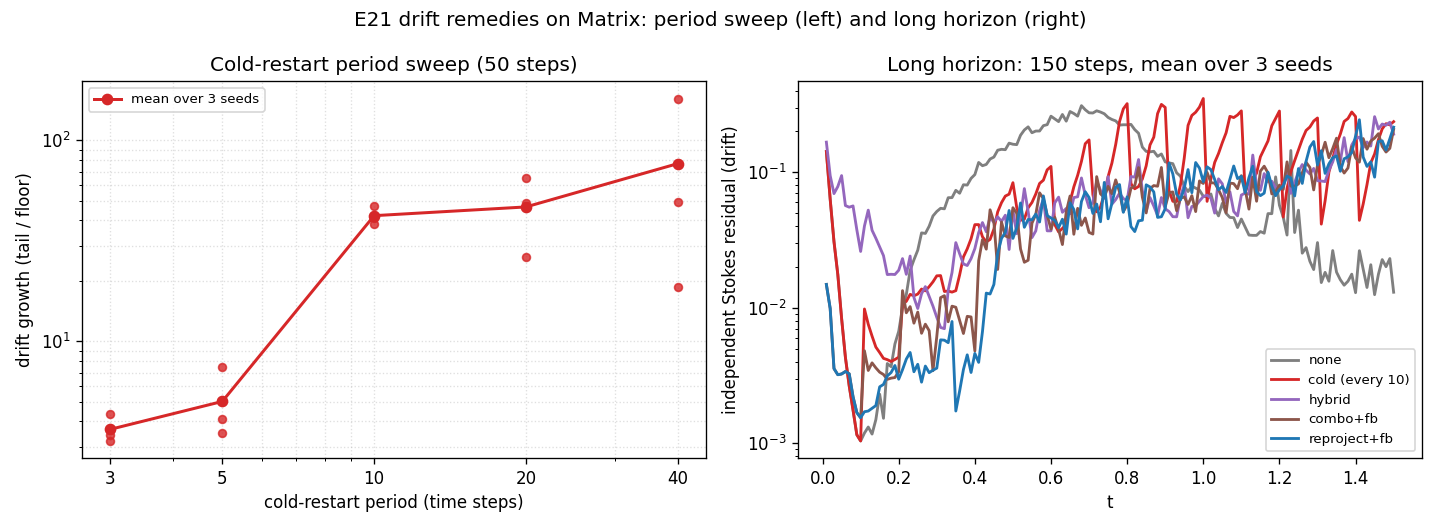}{Left: cold-restart period sweep, drift growth vs.\ restart period. Right: 150-step long-horizon Stokes-residual monitor vs.\ time.}
\caption{E21 drift remedies. \emph{Left:} cold-restart period
sweep (full $50$-step free-fall, $3$ seeds as dots, mean as line):
drift growth falls monotonically, with no U-shaped sweet spot. \emph{Right:}
long horizon ($150$ steps), geometric-mean off-manifold Stokes residual
over $3$ seeds; every remedy drifts back to $O(0.1)$ by $t=1.5$.}
\label{fig:e21}
\end{figure}

The headline cure is the cold-restart period (Table~\ref{tab:e21},
Figure~\ref{fig:e21} left): drift growth falls \emph{monotonically} as
restarts become more frequent, from $G\approx77$ at every $40$ steps
(${\approx}$ do-nothing) down to $G\approx3.7$ at every $3$, with no
U-shaped sweet spot; the cost is wall time, not residual. Restarting
every $3$--$5$ steps holds the drift to $G\approx3.7$--$5.0$ and is
seed-stable. Alternatives are weaker or less dependable: \emph{strong}
reaches $G\approx37$ at $2\times$ the wall, \emph{Adam-reset} barely
helps ($G\approx120$--$230$), and \emph{reproject} is sharply
\emph{bimodal} across five seeds (only $1/5$ below ${\sim}4$). Pushed to
$150$ steps (Figure~\ref{fig:e21} right) the picture sobers: every
remedy lets the absolute residual creep back to $O(0.1)$, the
do-nothing baseline is bimodal across seeds, and the cheap cold/period
cure that dominates at $50$ steps now mostly trades the drift's rate for
cost. The sober verdict is that no tested remedy \emph{eliminates}
$150$-step drift; frequent cold restart is the cheapest way to slow it.
A single-seed Fourier-feature probe suppressed the absolute drift
${\sim}20$--$30\times$ at $50$ steps, but one combination diverged, so we
record it only as a lead for a multi-seed follow-up.

\begin{table}[h]
\centering
\caption{Curing the warm-start drift: cold-restart period sweep and
alternatives on the E12 free-fall, scored by the training-independent
off-manifold Stokes-residual growth factor $G$ over $50$ steps (lower is
better; mean over three seeds, \emph{reproject} over five).}
\label{tab:e21}
\small
\begin{tabular}{@{}llrlrr@{}}
\toprule
 & remedy / period & $G$ (mean) & seed spread & extra Adam & wall \\
\midrule
\multirow{5}{*}{cold}
 & every $3$  & $\mathbf{3.7}$ & $3.2$--$4.3$  & $0$ & $960$\,s \\
 & every $5$  & $\mathbf{5.0}$ & $3.5$--$7.5$  & $0$ & $740$\,s \\
 & every $10$ & $42$           & $38$--$47$    & $0$ & $570$\,s \\
 & every $20$ & $47$           & $26$--$65$    & $0$ & $560$\,s \\
 & every $40$ & $77$           & $19$--$161$   & $0$ & $570$\,s \\
\midrule
strong & $n_{\mathrm{ref}}{=}1000$ & $37$ & $29$--$44$ & $0$ & $1260$\,s \\
reproject & ($5$ seeds) & $85$ & $3.6$--$171$ & $6$k--$60$k & $430$\,s \\
Adam-reset & every $10$--$20$ & $120$--$230$ & --- & $0$ & $160$--$470$\,s \\
\bottomrule
\end{tabular}
\end{table}

\section{Probing the inner-solve floor: optimisation, sampling, and architecture (E22)}
\label{app:E22}

Both limitations of Section~\ref{sec:fsi:limits} bottom out at the PINN
inner-solve floor, and every experiment in the body trained with Adam, on
a uniform collocation set, under a soft incompressibility penalty. The
Fourier-feature ablation (Figure~\ref{fig:e11-ablation}) already showed
that injecting high frequencies through the input basis does not move the
cap. Here we test the three remaining heavy levers---a second-order inner
solve, adaptive collocation, and improved conditioning---against two
diagnostics: the Stokes Steklov cap (the resolvable interface spectrum,
measured as in Appendix~\ref{app:sn-probe} but against a Taylor--Hood
Stokes reference) and the pre-impact squeeze-film drag ratio
$F_y^{\mathrm{PINN}}/F_y^{\mathrm{FE}}$ for a disc descending toward the
wall. All three drivers are additive and share the cap diagnostic; none
touches the coupling source tree.

\paragraph{The cap does not move (Table~\ref{tab:e22-cap}).}
We re-measure the leading $k\ge1$ Steklov eigenvalue under
optimiser~$\times$~width and under architecture~$\times$~loss-balancing,
reporting the operator ratio $\hat\sigma_{\max}^{k\ge1}(\text{PINN})/\sigma_{\max}^{k\ge1}(\text{FE})$
against a Taylor--Hood reference whose leading eigenvalue is
${\approx}120$. A quasi-Newton (L-BFGS) polish roughly triples the
Adam-only ratio---from $0.03$--$0.18$ to $0.24$--$0.50$, with the single
best value $0.50$ at width $128$---yet the ratio never reaches unity and
in every configuration stays below the $0.5$ threshold, so the PINN does
not resolve even the first interface mode in full.\footnote{A damped
Gauss--Newton / energy-natural-gradient variant was also implemented but
did not complete within the wall-clock budget, as it forms and solves
$J^\top J$ at each iteration; the reported second-order arm is L-BFGS.}
Gradient-norm loss balancing is actively harmful, collapsing the $k\ge1$
response to near zero ($\text{ratio}\approx0.005$), and the exactly
divergence-free stream-function ansatz is \emph{worse} than the plain
velocity--pressure formulation at matched width. No gradient-descent
lever---optimiser, sampler, or conditioning---moves the cap at fixed
width. This is strong evidence the cap is an \emph{optimisation} floor
rather than a basis or sampling artefact; Appendix~\ref{app:E23} confirms
it directly by lifting the cap with a non-gradient (least-squares) inner
solve.

\begin{table}[h]
\centering
\caption{E22 Stokes Steklov cap. Operator ratio (PINN over Taylor--Hood,
leading $k\ge1$ eigenvalue ${\approx}120$); higher is better, $1.0$ would
match the finite-element reference. \emph{Top:} second-order inner solve
$\times$ width (velocity--pressure, fixed loss weights). \emph{Bottom:}
architecture $\times$ loss balancing at width $64$. In every run the
resolvable cap stays at the first mode (no configuration crosses the
$0.5$ threshold).}
\label{tab:e22-cap}
\small
\begin{tabular}{@{}llrrr@{}}
\toprule
\multicolumn{2}{@{}l}{second-order $\times$ width} & $w{=}32$ & $w{=}64$ & $w{=}128$ \\
\midrule
\multicolumn{2}{@{}l}{Adam}        & $0.03$ & $0.16$ & $0.18$ \\
\multicolumn{2}{@{}l}{L-BFGS}      & ---    & $0.25$ & $0.25$ \\
\multicolumn{2}{@{}l}{Adam $+$ L-BFGS} & $0.24$ & $0.40$ & $\mathbf{0.50}$ \\
\midrule
ansatz & balancing & \multicolumn{3}{c}{ratio ($w{=}64$)} \\
\midrule
\multirow{2}{*}{velocity--pressure}
 & fixed     & \multicolumn{3}{c}{$0.18$ (Adam), $0.46$ (L-BFGS)} \\
 & grad-norm & \multicolumn{3}{c}{$0.004$} \\
\multirow{2}{*}{stream-function}
 & fixed     & \multicolumn{3}{c}{$0.10$ (Adam), $0.19$ (L-BFGS)} \\
 & grad-norm & \multicolumn{3}{c}{$0.007$} \\
\bottomrule
\end{tabular}
\end{table}

\paragraph{Scaling gradient descent up does not lift the cap either
(E24).} A GPU scale-up of the same diagnostic (Adam on four H100
accelerators, $8192$ collocation points) separates the two remaining
gradient-descent axes, budget and size. Budget helps, but only to a
plateau: at width $128$, depth $5$, doubling the budget from $10^{5}$
to $2\times10^{5}$ steps moves the leading $k\ge1$ ratio only from
$0.50$ to $0.53$ (receding the resolvable cap from $k=1$ to $k=4$),
still barely half the finite-element reference. Size actively hurts: at
comparable budgets the ratio falls to $0.48$ at width $256$, $0.28$ at
depth $10$, and $0.10$ at width $512$/depth $8$, a monotone degradation
with network size that is the opposite of a representation-capacity
wall and consistent with the spectral-bias reading. The wall cost of
these runs ($9$--$24$\,h each) also marks the practical end of the
budget axis.

\paragraph{The drag deficit is not closed and deepens with the gap
(Table~\ref{tab:e22-drag}).} Concentrating collocation in the lubrication
layer---residual-weighted (RAD), residual-refined (RAR), or explicit
gap-weighting---does not recover the bulk Stokes drag. The best ratio is
${\approx}0.03$, obtained at the widest gap, and it \emph{worsens} as the
gap narrows: the Taylor--Hood reference drag sharpens from $-237$ at gap
$0.05$ to $-1017$ at gap $0.015$ while the PINN stays near $-8$, so the
ratio falls to ${\approx}0.008$. The residual-driven samplers do worse
than uniform---RAR even returns the wrong sign at the tightest gap---and
an L-BFGS polish at the hardest gap changes the ratio only in the fourth
decimal. Packing the layer cannot substitute for inner-solve capacity.

\begin{table}[h]
\centering
\caption{E22 squeeze-film drag ratio $F_y^{\mathrm{PINN}}/F_y^{\mathrm{FE}}$
(closer to $1$ is better) for a disc descending toward the wall, by
sampler and wall gap. The Taylor--Hood reference drag is $-237$, $-417$,
$-1017$ at gaps $0.05$, $0.03$, $0.015$. An L-BFGS polish at gap $0.015$
(not shown) leaves every entry unchanged to three decimals.}
\label{tab:e22-drag}
\small
\begin{tabular}{@{}lrrr@{}}
\toprule
sampler & gap $0.05$ & gap $0.03$ & gap $0.015$ \\
\midrule
uniform       & $0.032$ & $0.014$ & $0.006$ \\
gap-weighted  & $\mathbf{0.033}$ & $0.012$ & $0.008$ \\
RAD           & $0.003$ & $0.004$ & $0.001$ \\
RAR           & $0.003$ & $0.001$ & $-0.000$ \\
\bottomrule
\end{tabular}
\end{table}

\section{The cap is an optimisation effect: a non-gradient inner solve (E23)}
\label{app:E23}

Appendix~\ref{app:E22} shows that no \emph{gradient-descent} lever moves
the Steklov cap. Every one of those runs shares a backbone---a global
network trained by gradient descent---whose neural-tangent-kernel spectral
bias is exactly the mechanism that stalls high-frequency content
\citep{WangYuPerdikaris2022,RahamanEtAl2019}. To test whether the cap is
that optimisation pathology or a genuine representational limit, we remove
gradient descent. We replace the inner Stokes solve with a
\emph{random-feature} model: a frozen random-Fourier hidden layer
$\phi_j(x)=\sin(\sigma\,W_j\!\cdot\!x+b_j)$ with $(W_j,b_j)$ sampled once
and held fixed, and a trainable \emph{linear} head $[u,v,p]=\phi(x)\,\Theta$.
Because the velocity, pressure, momentum and divergence operators are all
linear in $\Theta$, the homogeneous-Stokes Steklov problem becomes a single
weighted-collocation least-squares solve for $\Theta$---the random-feature
method / extreme-learning-machine PDE solver
\citep{ChenChiEYang2022,DongLi2021}---with no gradient descent and hence no
spectral bias. We sweep the feature count $M$ (capacity) and the frequency
scale $\sigma$, and re-measure the same Fourier Steklov cap.

\paragraph{Removing gradient descent lifts the cap and flattens it
(Table~\ref{tab:e23-ratio}).} The least-squares solve raises the leading
$k\ge1$ operator ratio to ${\approx}0.60$, above the best gradient-descent
value ($0.50$, Table~\ref{tab:e22-cap}), once the frequency scale matches
the feature budget (a scale too high for the number of features leaves the
system underdetermined and the solve fails---the $0.11$ entry). The ratio
then \emph{saturates} at $0.597$ across $M\in\{800,2000,3200\}$ and
$\sigma\in\{8,16,32\}$, identical to three decimals, with the mean
least-squares residual driven down to ${\approx}3\%$ at the largest
$M,\sigma$ (so the plateau is a converged solve, not a conditioning
artefact). The random-Fourier basis ($\sin$) outperforms a tanh
extreme-learning-machine ($0.60$ vs $0.34$). Decisively, the
\emph{per-mode} ratio is now \emph{flat} in mode number
(Table~\ref{tab:e23-fefe}, last column), rather than decaying as it does
under gradient descent: the high-frequency collapse---the cap proper---is
gone.

\begin{table}[h]
\centering
\caption{E23 random-feature least-squares solve: leading $k\ge1$ Steklov
operator ratio (PINN over Taylor--Hood reference) versus feature count $M$
and frequency scale $\sigma$ (random-Fourier $\sin$ features, plus one
tanh control). Compare the best gradient-descent value $0.50$
(Table~\ref{tab:e22-cap}). The ratio saturates at $0.597$ once $\sigma$
matches $M$; ``--'' marks the underdetermined regime (scale too high for
the feature budget, least-squares residual ${\approx}0.95$).}
\label{tab:e23-ratio}
\small
\begin{tabular}{@{}lrrrr@{}}
\toprule
 & $\sigma{=}4$ & $\sigma{=}8$ & $\sigma{=}16$ & $\sigma{=}32$ \\
\midrule
$M{=}800$  & $0.37$ & $0.59$ & $0.60$ & --\,$(0.11)$ \\
$M{=}2000$ & $0.40$ & $0.59$ & $\mathbf{0.60}$ & $\mathbf{0.60}$ \\
$M{=}3200$ & ---    & $0.59$ & $\mathbf{0.60}$ & $\mathbf{0.60}$ \\
\midrule
tanh ($M{=}2000$) & & & $0.34$ & \\
\bottomrule
\end{tabular}
\end{table}

\paragraph{The residual ${\approx}0.57$ offset is a traction-reader
convention (Table~\ref{tab:e23-fefe}).} A uniform, frequency-flat factor
does not look like spectral bias; it looks like a fixed mismatch between
two ways of reading the interface traction. The reference operator
$\SF$ (the Taylor--Hood Schur complement) reads the \emph{mass-consistent}
integrated flux $M_\Sigma^{-1}\!\int_\Sigma(\sigma n)\,\varphi$, whereas
the random-feature and PINN probes read the \emph{analytic pointwise}
stress $\sigma\cdot n=(-pI+2\mu_f\varepsilon(u))\cdot n$ at the $\Sigma$
vertices. To isolate this we drive the \emph{same} Fourier probe with an
genuine finite-element Stokes solver and read its traction both ways. The
mass-consistent reader reproduces the reference exactly (ratio $1.000$ at
every mode, confirming the probe pipeline is faithful); the
analytic-pointwise reader---the same functional the network uses---lands
at $0.57$ of the reference, with the \emph{same} per-mode profile as the
random-feature solve, which it tracks to within a few percent throughout
(Table~\ref{tab:e23-fefe}). The $0.57$ offset is therefore a property of
the pointwise-traction convention shared by any solver, not a deficit of
the network: read consistently, the random-feature interface operator
agrees with the finite-element one. The cap is an optimisation
(spectral-bias) effect, removable by a non-gradient inner solve; it is not
a representation-capacity wall.

\begin{table}[h]
\centering
\caption{FE--FE convention check. Per-mode $\max|\mathrm{eig}|$ ratios to
the reference $\SF$ (the Taylor--Hood Schur complement projected to mode $k$),
$\cos$ polarisation. \emph{FE-consistent}: finite-element solver read with
the mass-consistent flux (reproduces the reference). \emph{FE-pointwise}:
\emph{same} finite-element solver read with the analytic pointwise
$\sigma\cdot n$ the network uses. \emph{RFM}: the random-feature
least-squares solve ($M{=}2000$, $\sigma{=}32$) read the same pointwise
way. The pointwise reader sits at ${\approx}0.57$ for the finite-element
solver itself, and the random-feature solve tracks it to a few percent.}
\label{tab:e23-fefe}
\small
\begin{tabular}{@{}rrrr@{}}
\toprule
$k$ & FE-consistent & FE-pointwise & RFM \\
\midrule
$1$ & $1.000$ & $0.602$ & $0.600$ \\
$2$ & $1.000$ & $0.585$ & $0.561$ \\
$3$ & $1.000$ & $0.556$ & $0.522$ \\
$4$ & $1.000$ & $0.545$ & $0.495$ \\
$5$ & $1.000$ & $0.545$ & $0.487$ \\
$6$ & $1.000$ & $0.575$ & $0.507$ \\
\bottomrule
\end{tabular}
\end{table}

\section{Contact-window dynamics: convergence failures and the FEM--FEM benchmark}
\label{app:contact-sweeps}

This appendix documents in full the pre-impact dynamic signatures that
Section~\ref{sec:fsi:limits} reports as PINN under-resolution artefacts.

\begin{figure}[h]
\centering
\figmaybe[0.7\linewidth]{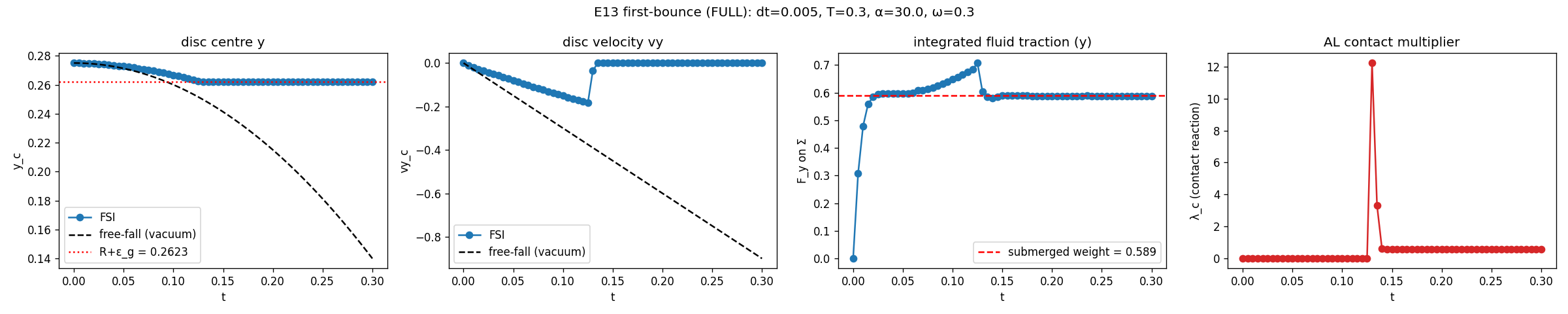}{E13(d) no-slip baseline: trajectories of $y_c$, $\dot y_c$, $F_y$, $\lambda_c$.}
\caption{E13(d), no-slip baseline. Impact at $t=0.130$ (one $\dt$ from
the vacuum prediction $0.129$); the AL projector lands the disc exactly
at $R+\varepsilon_g$; $\lambda_c$ converges to the submerged weight
$\Pi=0.589$. The $+20\%$ pre-impact $F_y$ overshoot appears only at
$(n_{\mathrm{iface}},n_{\mathrm{steps}}^{\mathrm{refine}})=(32,300)$ and
is not preserved under refinement
(Figures~\ref{fig:p0-e13d-niface}, \ref{fig:p0-e13d-refine}).}
\label{fig:e13d}
\end{figure}

\begin{figure}[h]
\centering
\figmaybe[0.95\linewidth]{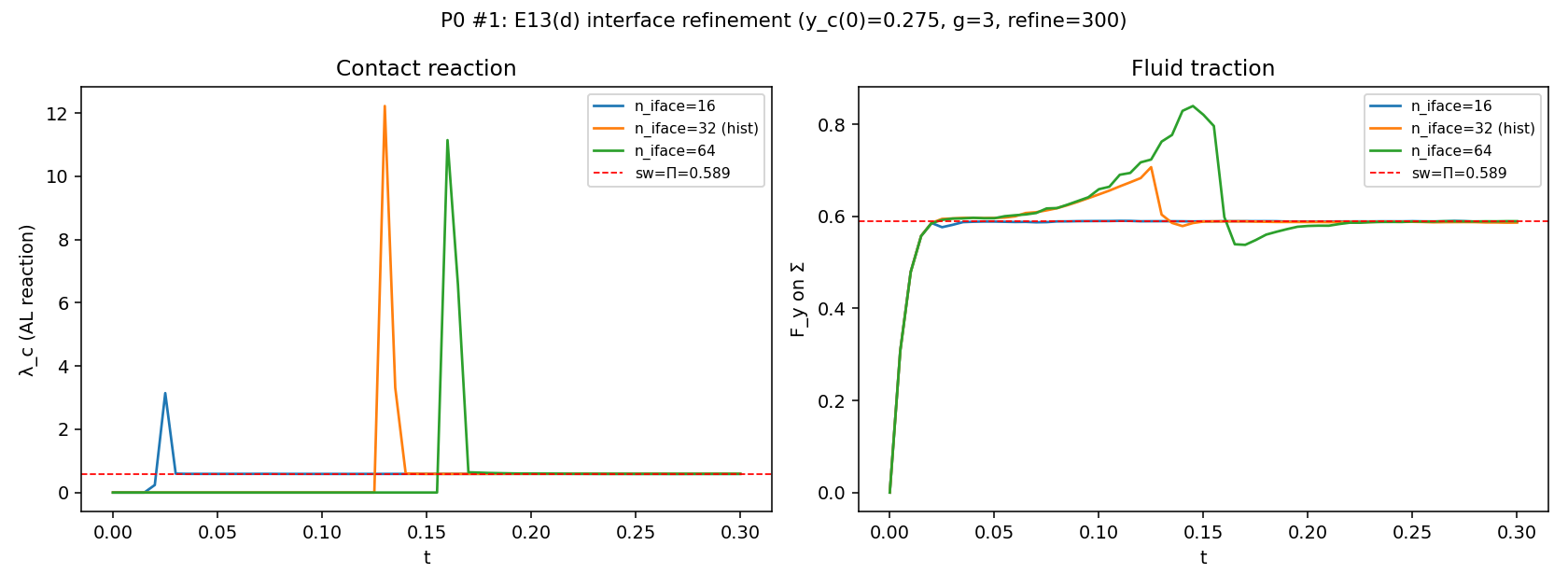}{$\lambda_c(t)$ and $F_y(t)$ for E13(d) at $n_{\mathrm{iface}}\in\{16,32,64\}$.}
\caption{Static-equilibrium robustness under interface refinement (E13d,
$\mathrm{refine}=300$). Left: the post-impact tail ($t\ge0.20$) settles
at $\lambda_c\approx\Pi=0.589$ for all three resolutions, within $0.4\%$.
Right: $F_y$ pre-impact peak grows with $n_{\mathrm{iface}}$
($0.59\to0.71\to0.83$) while the post-equilibrium tail again settles at
$\Pi$. ($n_{\mathrm{iface}}=16$ is partly degenerate because
$\varepsilon_g=h_{\mathrm{iface}}/4$ scales with the mesh.)}
\label{fig:p0-e13d-niface}
\end{figure}

\begin{figure}[h]
\centering
\figmaybe[0.78\linewidth]{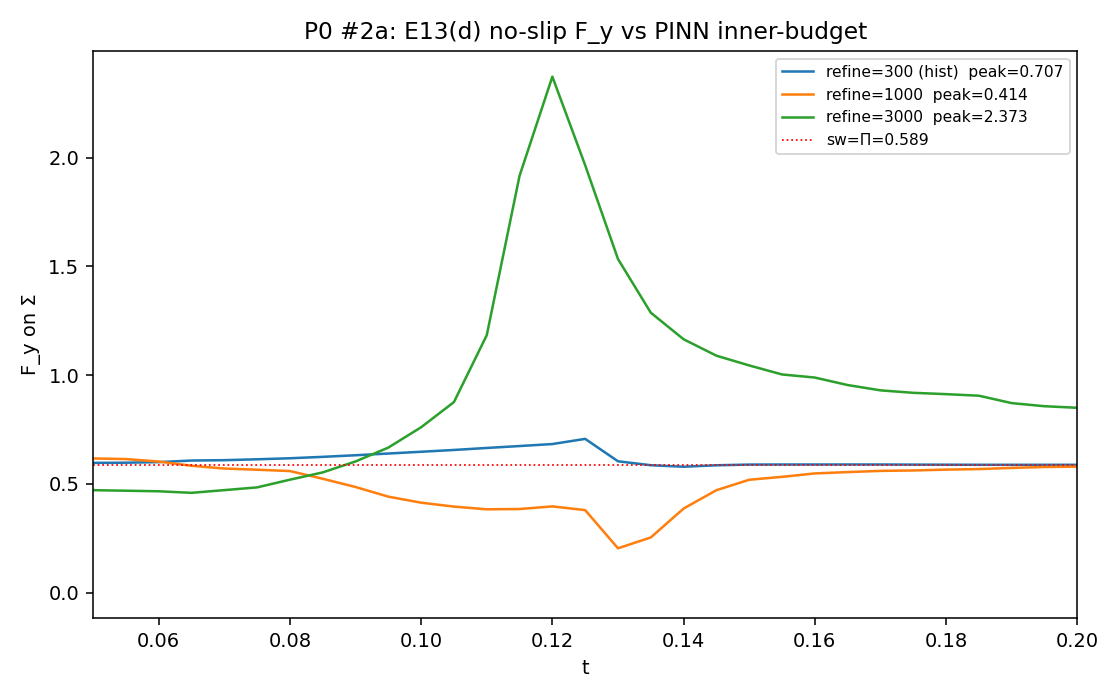}{$F_y(t)$ at $n_{\mathrm{steps}}^{\mathrm{refine}}\in\{300,1000,3000\}$.}
\caption{Pre-impact $F_y$ peak is not budget-converged on E13(d): a mild
overshoot at $\mathrm{refine}=300$ ($0.71$), a \emph{dip} below $\Pi$ at
$1000$ (peak $0.41$ mid-descent), and a sharp spike at $3000$ ($2.37$).}
\label{fig:p0-e13d-refine}
\end{figure}

\paragraph{E13(d) pre-impact peak is not budget-converged.}
Figure~\ref{fig:e13d} shows the full no-slip E13(d) trajectory ($y_c$,
$\dot y_c$, $F_y$, and the contact reaction $\lambda_c$). At
$(n_{\mathrm{iface}},n_{\mathrm{steps}}^{\mathrm{refine}})=(32,300)$,
$F_y$ shows a $+20\%$ pre-impact overshoot, but this is a
budget-specific artefact: at $\mathrm{refine}=1000$ the squeeze-film
phase \emph{dips} below $\Pi$ (peak $0.41$ mid-descent), and at
$\mathrm{refine}=3000$ the pre-impact peak jumps to $2.37$
(${\sim}3.3\times$). The static-equilibrium claim $\lambda_c\to\Pi$, by
contrast, survives both refinement axes cleanly
(Figure~\ref{fig:p0-e13d-niface}, \ref{fig:p0-e13d-refine}).

\paragraph{FEM--FEM benchmark.} Replacing the PINN-Stokes trainer by a
remeshed-per-$\Delta t$ Taylor--Hood Stokes solver at the identical
operating point, the
FE disc \emph{never reaches} the wall in $T=0.30$: it settles into a
near-static hover at $\dot y_c=-0.0023$ with $F_y=1.178=m_sg$ exactly.
The inferred drag $c\approx510$ vs.\ the PINN's $\approx6$ at impact is an
${\sim}80\times$ discrepancy; the wall accounts for ${\sim}2.5\times$
(centered-disc baseline $c\approx200$), and mesh refinement to $h=0.025$
pushes the FE drag \emph{up} ($c\to755$), so the ground-truth Stokes drag
is at least an order of magnitude larger than $c\approx6$. The
pre-impact $F_y$ peak is therefore a PINN under-resolution of bulk Stokes
drag, not a physical squeeze-film signature.

\begin{figure}[h]
\centering
\figmaybe[0.95\linewidth]{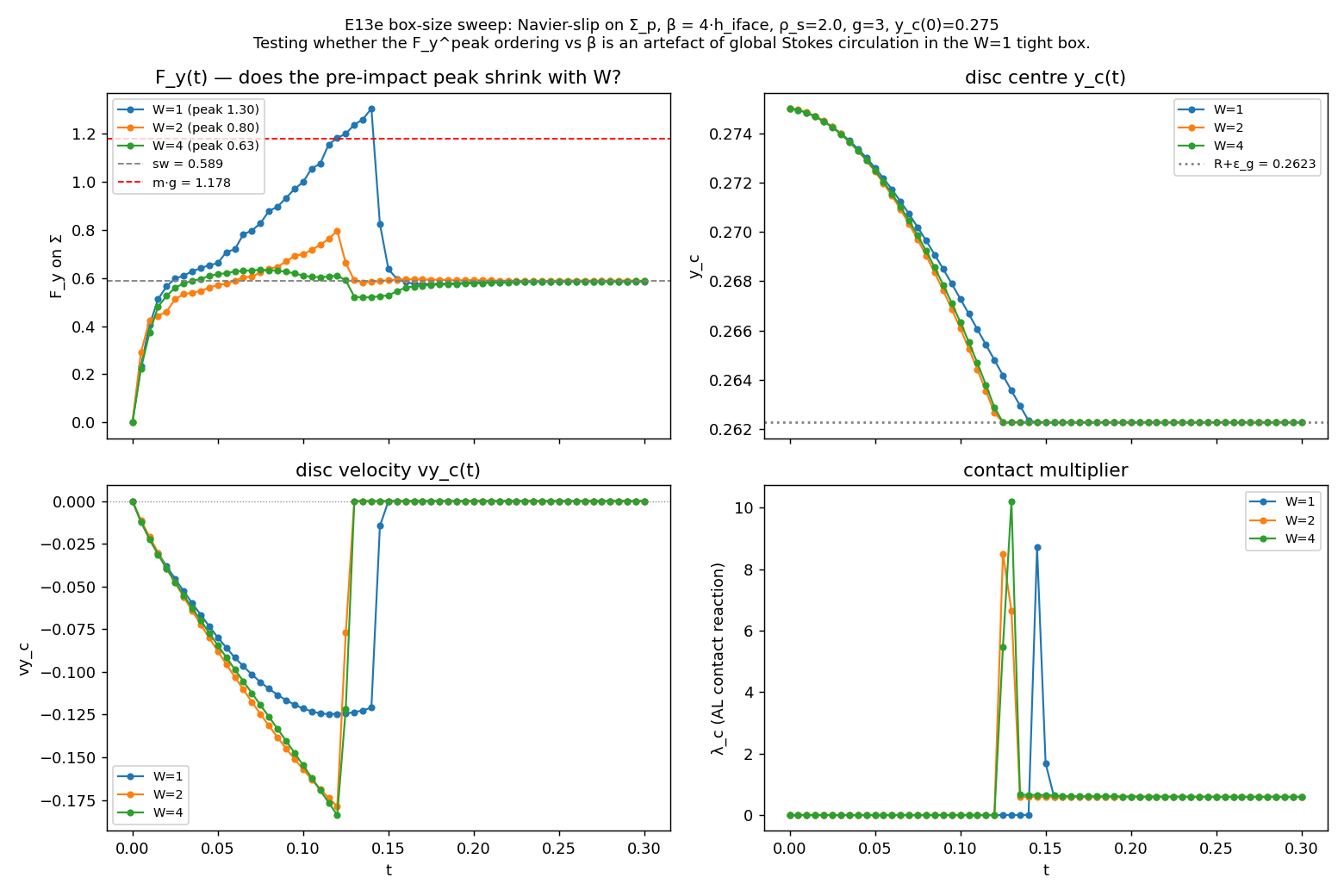}{E13(e) box-size sweep: $F_y$, $y_c$, $\dot y_c$, $\lambda_c$ vs.\ $t$ overlaid for $W\in\{1,2,4\}$.}
\caption{E13(e) box-size sweep at $\beta=4h$, $\mathrm{refine}=300$. The
pre-impact $F_y$ peak drops from $1.30$ ($W=1$) to $0.80$ ($W=2$) to
$0.63$ ($W=4$), i.e.\ $+121\%\to+35\%\to+8\%$ above $\Pi=0.589$. At $W=4$
the slip-side overshoot is essentially gone, isolating it as primarily a
global Stokes-circulation artefact of the confined box.}
\label{fig:e13e-box-sweep}
\end{figure}

\paragraph{E13(e) Navier-slip sweep.} With
$\beta\in\{h/4,h,4h\}$, $F_y^{\mathrm{peak}}$ rises from $+20\%$
(no-slip) to $+121\%$ ($\beta=4h$) and $\dot y_c$ at impact drops from
$-0.183$ to $-0.124$; all three still fail the $>1$\,mm rebound gate. The
peak is again non-monotone in $n_{\mathrm{steps}}^{\mathrm{refine}}$
($1.30\to1.70\to1.22$) and the regime flips at higher budgets, so the
qualitative ordering ($F_y^{\mathrm{peak}}$ increases with $\beta$)
survives but the specific $+121\%$ figure does not. A box-size sweep
($W\in\{1,2,4\}$, $\beta=4h$) shows the peak collapses monotonically
($+121\%\to+35\%\to+8\%$): the $\beta$-monotonicity at $W=1$ is largely a
global-Stokes-circulation artefact of the tight box
(Figure~\ref{fig:e13e-box-sweep}).

\begin{figure}[h]
\centering
\begin{subfigure}[b]{0.62\linewidth}
  \centering
  \figmaybe{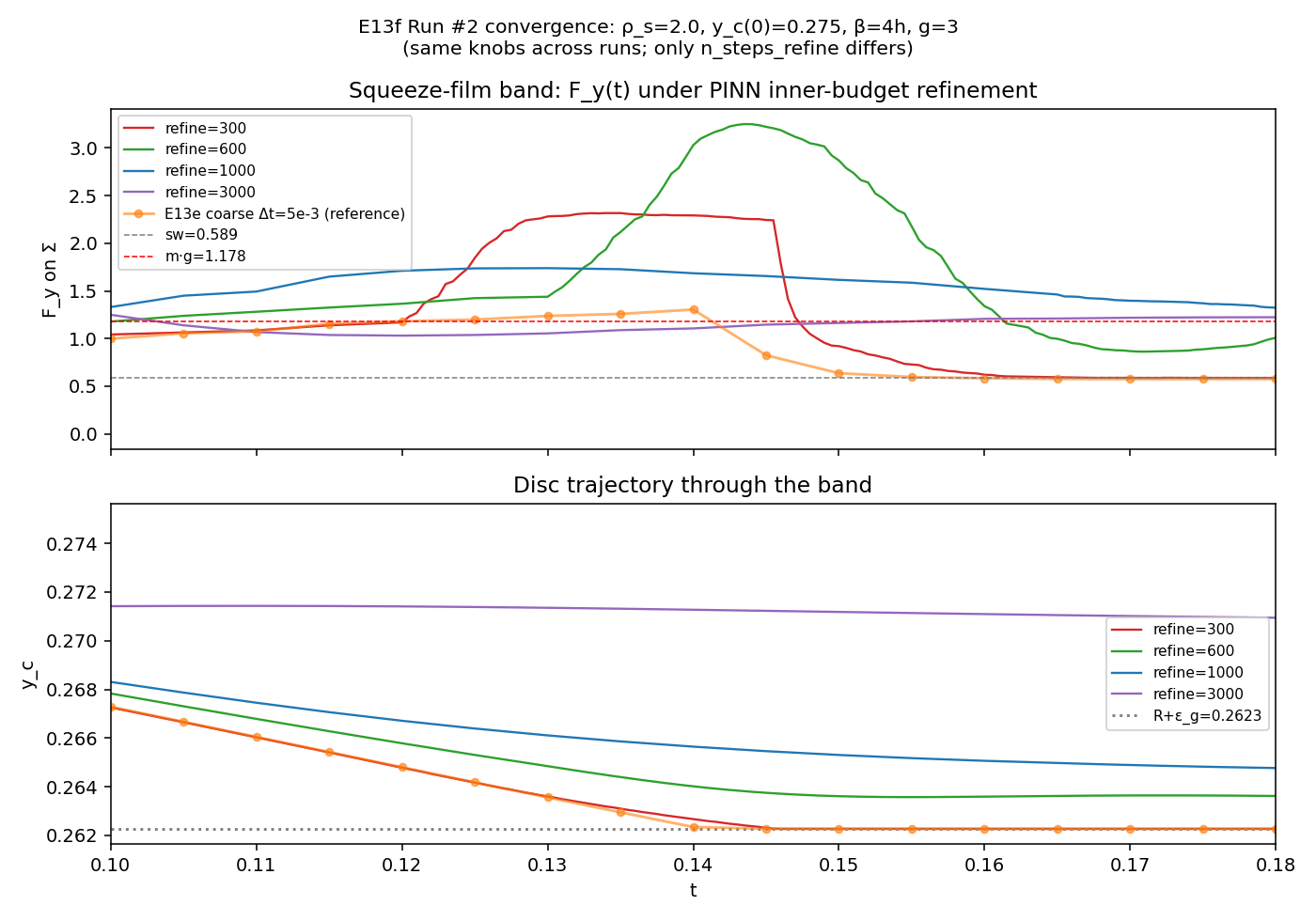}{$F_y(t)$ and $y_c(t)$ overlay across refine budgets.}
  \caption{$F_y(t)$ (top) and disc trajectory (bottom) for each
  refinement budget, coarse-$\dt$ E13(e) reference in orange.}
\end{subfigure}\hfill
\begin{subfigure}[b]{0.34\linewidth}
  \centering
  \figmaybe{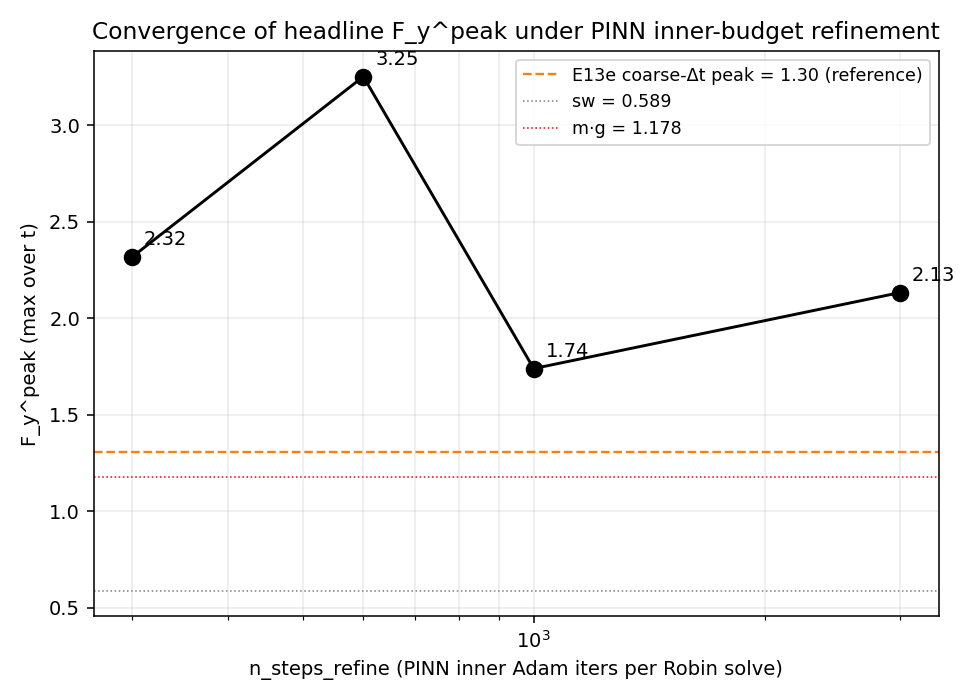}{$F_y^{\mathrm{peak}}$ vs.\ refine budget.}
  \caption{Peak $F_y$ vs.\ inner budget (log-$x$): $2.32\to3.25\to1.74\to
  2.13$, non-monotone, not on an asymptote.}
\end{subfigure}
\caption{Convergence failure of the contact-window signature. At
$\mathrm{refine}=3000$ the disc never enters the adaptive-$\dt$ region;
the squeeze-film picture vanishes. We treat the contact-window signature
as a PINN inner-solve artefact.}
\label{fig:e13f}
\end{figure}

\paragraph{E13(f) adaptive $\dt$ and convergence failure.} With adaptive
time stepping near contact ($\dt=5\times10^{-4}$ when the gap drops below
$0.02$), a four-point inner-budget sweep
$n_{\mathrm{steps}}^{\mathrm{refine}}\in\{300,600,1000,3000\}$ shows the
peak $F_y$ is non-monotone ($2.32\to3.25\to1.74\to2.13$) and at the
highest budget the disc never enters the squeeze-film band at all
(Table~\ref{tab:e13f}, Figure~\ref{fig:e13f}). The contact-window
signature is a PINN inner-solve artefact under this combination of
parameters.

\begin{table}[h]
\centering
\caption{E13(f) Run~\#2 ($\rho_s=2.0$, $y_c(0)=0.275$, $\beta=4h$,
$g=3$) under PINN inner-budget refinement. The peak $F_y$ is non-monotone
and the trajectory regime itself flips: contact occurs only at
$\mathrm{refine}=300$.}
\label{tab:e13f}
\small
\begin{tabular}{@{}cccccp{4.0cm}@{}}
\toprule
$n_{\mathrm{steps}}^{\mathrm{refine}}$ & $F_y^{\mathrm{peak}}$ & $t^\star$ & $\min y_c$ & contact & regime \\
\midrule
$300$  & $2.32$ & $0.134$ & $0.262$ & yes & contact (Run~\#2 original) \\
$600$  & $3.25$ & $0.144$ & $0.263$ & no  & near-contact creep \\
$1000$ & $1.74$ & $0.130$ & $0.264$ & no  & creep band \\
$3000$ & $2.13$ & $0.070$ & $0.270$ & no  & no contact, free-fall arrest \\
\bottomrule
\end{tabular}
\end{table}

\section{1D split-segment Poisson coupling in detail}
\label{app:1d}

\begin{figure}[h]
\centering
\begin{subfigure}[b]{0.49\linewidth}
  \centering
  \figmaybe{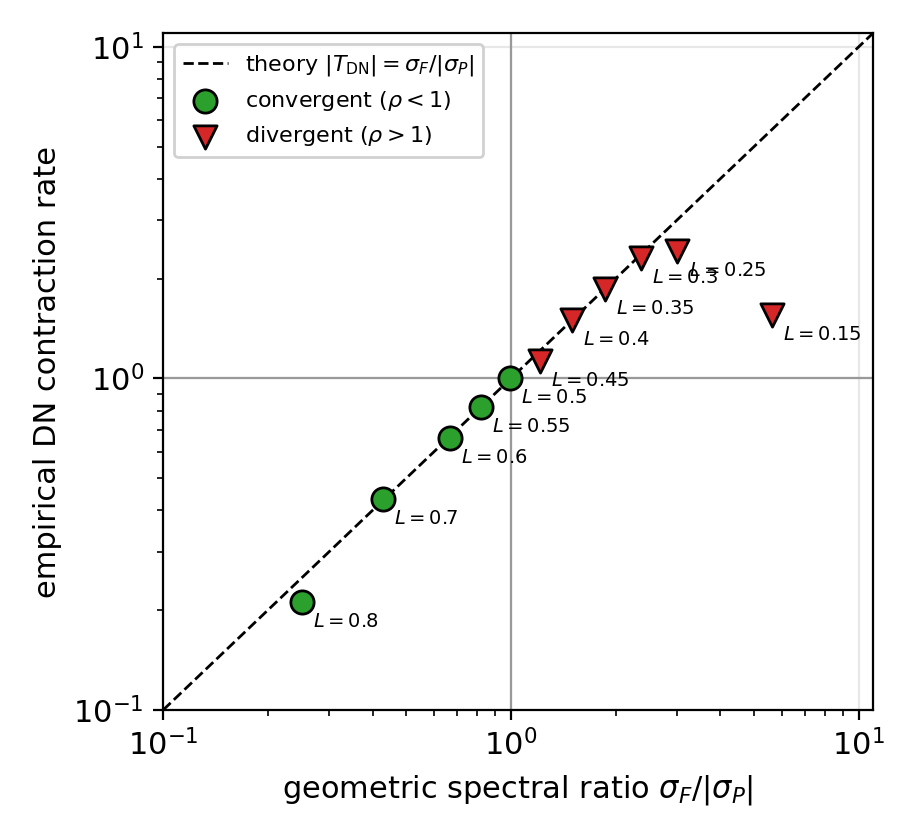}{DN rate vs.\ theory along $L$-sweep.}
  \caption{DN rate matches $\sigmaF/\abs{\sigmaN}=(1-L)/L$ within $2\%$ on
  $7/11$ points; the divergence sign is recovered on all $11$.}
  \label{fig:1d:a}
\end{subfigure}\hfill
\begin{subfigure}[b]{0.49\linewidth}
  \centering
  \figmaybe{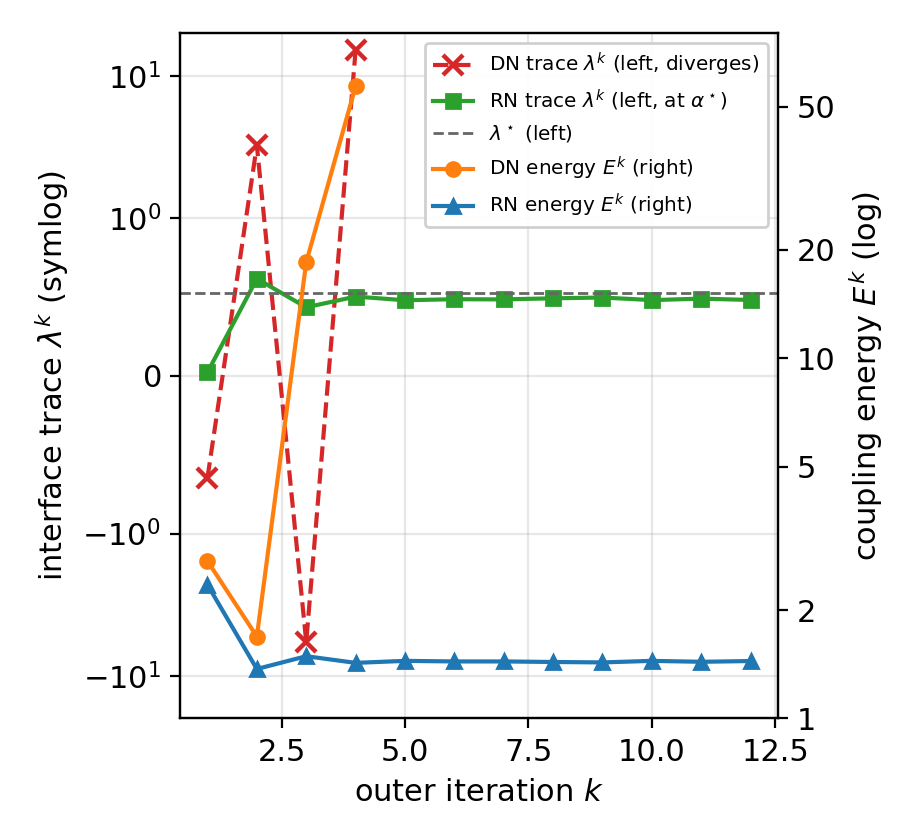}{1D trace and energy under DN vs.\ RN.}
  \caption{At $L=0.3$ (left axis, symlog): the RN interface trace
  $\lambda^k$ settles onto $\lambda^\star$, while the DN trace oscillates
  and diverges away from it. The coupling energy $E^k$ (right axis, log)
  grows ${\sim}30\times$ under DN and descends to a plateau under RN, as
  Theorem~\ref{thm:cotraining} predicts.}
  \label{fig:1d:b}
\end{subfigure}
\caption{1D Poisson PINN--FEM coupling: theory match (a) and the RN cure
with its energy signature (b). Notation follows the text:
$\sigmaF,\sigmaN$ are the FEM/PINN Steklov eigenvalues and $\lambda^k$
the interface trace at outer iteration $k$.}
\label{fig:1d}
\end{figure}

This appendix expands the 1D summary of
Section~\ref{sec:exp:poisson}. We split $\Omega=(0,1)$ at $x=L$, give
$(0,L)$ to a P1 FEM solver and $(L,1)$ to an MLP PINN with hard
right-Dirichlet ansatz $\uN(x)=(1-x)\NN_\theta(x)$. With manufactured
$u^\star(x)=\sin(2\pi x)e^{-2x}$ and an $L$-sweep densified near the
spectral crossover $L=0.5$, the empirical DN rate matches the theory
line $\sigmaF/\sigmaN=(1-L)/L$ to within $2\%$ on $7$ of $11$ sweep
points (within $7\%$ on $8$), recovering the divergence sign at all
$11$ (Figure~\ref{fig:1d:a}). At $L=0.3$ ($\alphastar\approx2.18$),
DN diverges within four iterations while RN stays bounded and contracts
to the inner-PINN noise floor. The accompanying coupling energy makes
the mechanism of Theorem~\ref{thm:cotraining} visible
(Figure~\ref{fig:1d:b}): the DN energy grows ${\sim}30\times$ over four
iterations, whereas the RN energy drops ${\sim}40\%$ and then plateaus.

\section{Inner-PINN (natural-orientation) results in detail}
\label{app:inner-pinn}

Section~\ref{sec:exp:2d} reports the headline natural-orientation numbers
(FEM on the annulus, PINN on the inner disc); this appendix collects the
full diagnostics. The orientation isolates the PINN as the \emph{bounded}
subdomain (the disc), where the continuum Steklov picture
$\sigmaN(k)=k/R$ holds exactly. Table~\ref{tab:e5} summarises E5.

\begin{table}[h]
\centering
\caption{2D Poisson (natural orientation, E5): post unit-conversion
PINN--FEM iteration rates match the FEM--FEM baseline (E8); RN at
$\alphastar$ contracts as Theorem~\ref{thm:rn} predicts. The two RN
columns are reported at \emph{different} impedances, each at its own
optimum, since the PINN-side spectrum differs from the all-FE
spectrum.}
\label{tab:e5}
\small
\begin{tabular}{@{}lcc@{}}
\toprule
Metric & PINN--FEM (E5) & FEM--FEM ref.\ (E8) \\
\midrule
DN empirical rate & $\mathbf{1.442}$ (geometric) & $1.452$ \\
RN rate at own $\alpha^\star$ & $\mathbf{0.78}$ ($\alpha = 4.16$) & $0.92$ ($\alpha = 0.61$) \\
RN field error, $k\!=\!1 \to 10$ & $0.025 \to 3{\times}10^{-3}$ & $5.4{\times}10^{-3} \to 2.1{\times}10^{-3}$ \\
Numerical $\alpha$-bowl min (E10) & $\alpha=10.31$, rate $0.668$ & $\alpha = 0.61$ \\
\bottomrule
\end{tabular}
\end{table}

The Steklov spectra are $\sigmaF\in[0.257,3.768]$ (annulus FEM) and
$\abs{\sigmaN}\in[0,67.4]$ (disc PINN; the basis-vector probe gives a
noisier, larger $\abs{\sigmaN}^{\max}$ than the Fourier probe). The
textbook $\alphastar=\sqrt{\sigmaF^{\min}\abs{\sigmaN}^{\max}}=4.16$; the
cleaned (baseline-subtracted) $S_N$ has $\abs{\sigmaN}^{\max}\approx8.05$,
giving $\alphastar=1.44$. A numerical $\alpha$-search (E10) recovers the
iteration optimum at $\alpha=10.31$ with empirical rate $0.668$; the
same sweep measures rate $0.752$ at its $\alpha=1.41$ grid point, so
the cleaned-spectrum default $1.44$ is also a serviceable operating
point, and the choice of probe moves the operating point within the
broad bowl rather than in or out of the contractive band.

\subsection{E5: full DN-vs-RN diagnosis}
\label{app:inner-pinn:e5}

\begin{figure}[h]
\centering
\begin{subfigure}[b]{0.49\linewidth}
  \centering
  \figmaybe{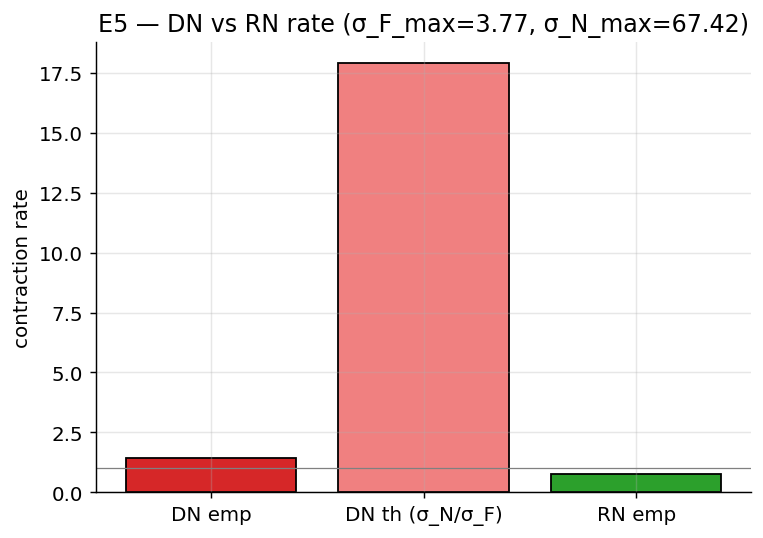}{E5 rate summary.}
  \caption{Empirical DN rate $1.442$ matches the FEM--FEM analogue (E8)
  rate $1.452$ to within $0.7\%$.}
\end{subfigure}\hfill
\begin{subfigure}[b]{0.49\linewidth}
  \centering
  \figmaybe{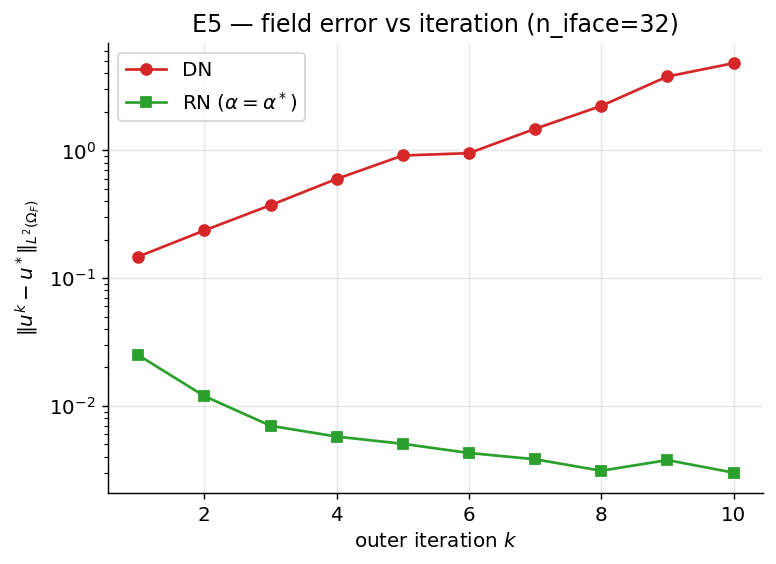}{E5 field error: RN at $\alpha^\star=4.16$.}
  \caption{RN at $\alphastar=4.16$ contracts the field error $0.025\to
  3{\times}10^{-3}$ in $10$ iterations; DN diverges geometrically.}
\end{subfigure}
\caption{Natural orientation (E5): DN vs.\ RN trace-error rate and
field-error contrast.}
\label{fig:app-e5}
\end{figure}

Figure~\ref{fig:app-e5} contrasts the two schemes on the natural
orientation: the empirical DN trace-error rate matching the FEM--FEM
analogue (left) and the RN field-error contraction at $\alphastar=4.16$
against the geometric DN divergence (right).

\subsection{E10: numerical \texorpdfstring{$\alpha$}{alpha}-search}
\label{app:inner-pinn:e10}

\begin{figure}[h]
\centering
\figmaybe[0.6\linewidth]{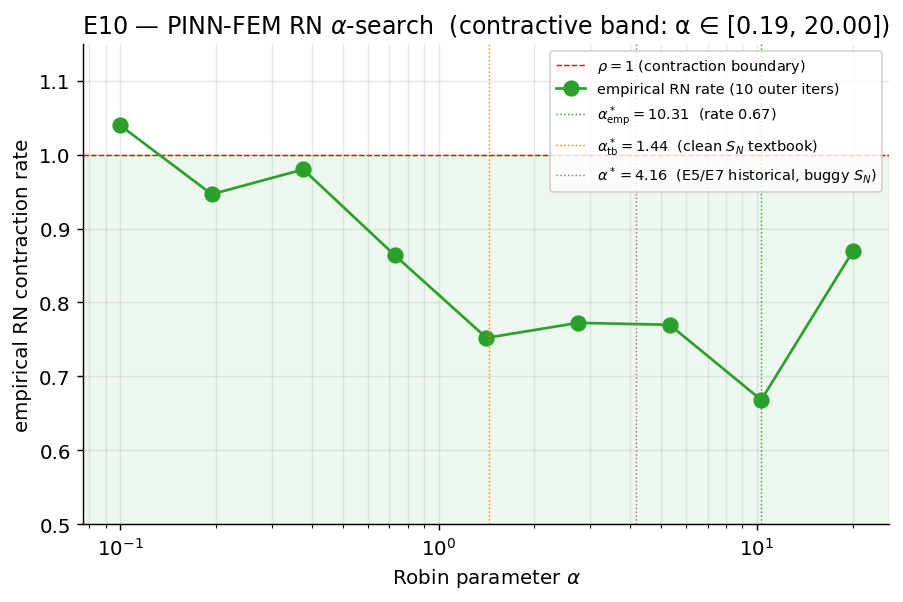}{E10 $\alpha$-bowl.}
\caption{E10 numerical $\alpha$-search on the inner-PINN orientation
($9$-point sweep over $[0.1,20]$, $\Kouter=10$, $\Tinner=1500$); minimum
at $\alpha=10.31$, rate $0.668$. The bowl is broad and the optimum is
shifted ${\sim}2.5\times$ from the closed-form $\alphastar=4.16$,
consistent with the operator-norm bound of Appendix~\ref{app:proofs}
being loose when the discrete eigenbases do not coincide.}
\label{fig:app-e10}
\end{figure}

Figure~\ref{fig:app-e10} shows the numerical $\alpha$-bowl for this
orientation; the broad minimum sits at $\alpha=10.31$, shifted from the
closed-form $\alphastar=4.16$ because the discrete eigenbases do not
exactly coincide.

\subsection{E7: co-training variant}
\label{app:inner-pinn:e7}

\begin{figure}[h]
\centering
\begin{subfigure}[b]{0.49\linewidth}
  \centering
  \figmaybe{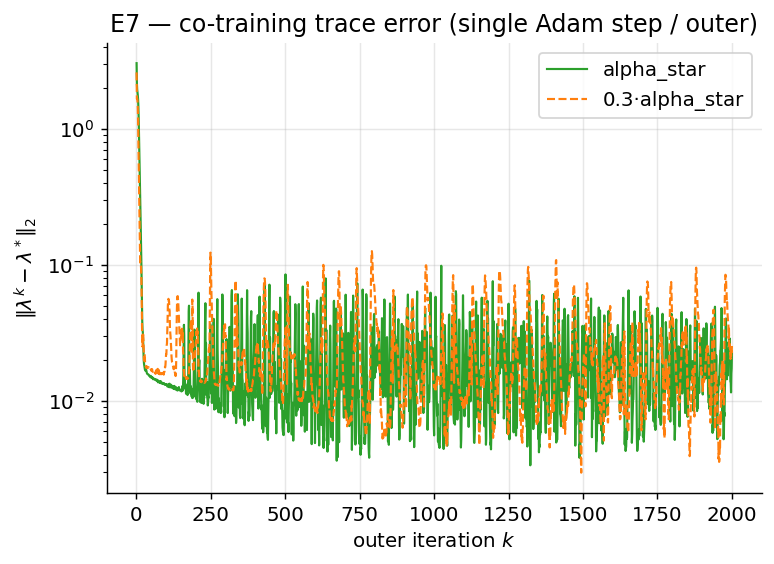}{E7 co-training trace error.}
  \caption{Trace error $3.06\to2.23{\times}10^{-2}$ in $\Kouter=2000$
  co-training steps at $\alpha=\alphastar$.}
\end{subfigure}\hfill
\begin{subfigure}[b]{0.49\linewidth}
  \centering
  \figmaybe{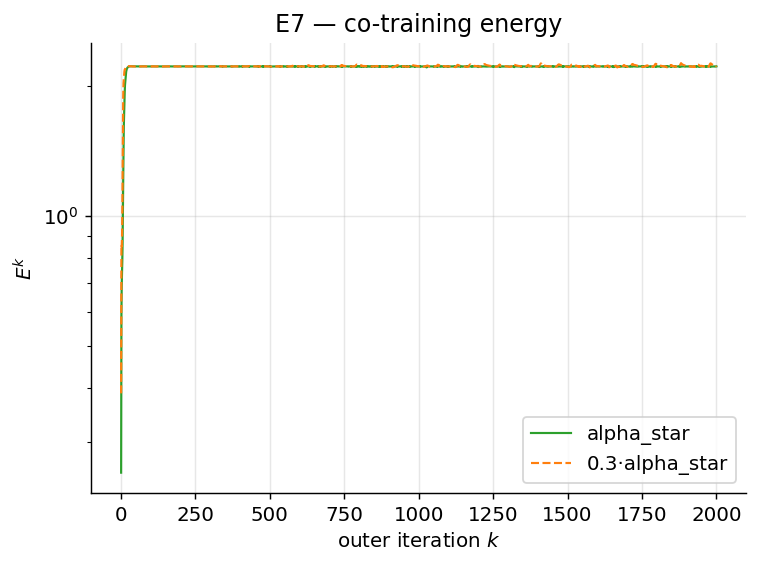}{E7 co-training composite energy.}
  \caption{Composite energy stays bounded around $2.22$ and does not
  grow, in agreement with Theorem~\ref{thm:cotraining}.}
\end{subfigure}
\caption{Inner-PINN co-training (E7) at $\alpha\in\{\alphastar,0.3
\alphastar\}$, $5$ PINN-side Adam steps per outer iteration, the
regime where the leading $-\frac{\alpha}{2}(1-\rho^2)\norm{e^k}^2$ term
of Appendix~\ref{app:proofs:cotraining} dominates the $O(\eta)$
remainder.}
\label{fig:app-e7}
\end{figure}

Figure~\ref{fig:app-e7} reports the co-training variant of
Theorem~\ref{thm:cotraining}, in which the PINN takes a few gradient
steps per outer iteration rather than solving to convergence: the trace
error contracts (left) while the composite coupling energy stays bounded
(right), as the theorem predicts.

\subsection{E26: head-to-head with the Schwarz PINN--FOM coupling of
Snyder et al.}
\label{app:inner-pinn:e6}

The closest prior PINN--FEM coupling is the Schwarz alternating method
of \citet{SnyderTezaurWentland2023}: a multiplicative \emph{overlapping}
Schwarz iteration exchanging \emph{Dirichlet} data at the two Schwarz
boundaries, the PINN's transmission condition imposed weakly through the
loss (their WDBC variant, the one they use for PINN--FOM coupling, where
they report the enforcement variant matters little), and the subdomain
network retrained at every Schwarz iteration. We implement the scheme
faithfully on the natural-orientation testbed and run it head-to-head
against DN and RN under fully matched conditions: the identical global
problem and manufactured solution, identical network (width $64$, depth
$5$), identical inner budget ($1500$ Adam steps per outer iteration,
identical collocation count and weak-BC weight, so the per-outer-iteration
cost is the same by construction, $8$--$10$\,s for every scheme),
the same perturbed initial trace $\lambda^\star+0.3\,\mathcal N(0,I)$ at
each scheme's own interface, cold network initialisation warm-started
across outer iterations, and no manufactured-flux warm start anywhere
(RN starts from $g^0=0$). The FE subdomain is the square minus a disc of
radius $r_1$ and the PINN disc has radius $r_2$, so the Schwarz
boundaries are $\Sigma_1$ (the FE inner boundary, which receives the
PINN trace) and $\Sigma_2$ (the PINN disc edge, which receives the FE
values), with overlap $\delta=r_2-r_1\in\{0.05,0.10\}$; DN/RN use the
sharp split at $r=0.25$. Their published experiments are 1D
advection--diffusion boundary layers, so this benchmark transplants
their \emph{algorithm} to our geometry rather than reproducing their
results. Three seeds; the error is the relative trace error at the
scheme's own interface against the manufactured solution.

\begin{table}[h]
\centering
\caption{E26 head-to-head at matched budget (median over three seeds;
``--'' = tolerance not reached within $\Kouter=25$, with the per-seed
success count in parentheses where it is not $3/3$). Schwarz converges
fast with a generous overlap but is seed-fragile at the small one and
non-monotone (its error drifts back up from the minimum under continued
warm-started retraining); RN at the closed-form $\alphastar$ is monotone
on every seed (final $=$ min) and reaches the lowest floor.}
\label{tab:e26-schwarz}
\small
\begin{tabular}{@{}lrrrr@{}}
\toprule
scheme & iters $\to2{\times}10^{-2}$ & iters $\to10^{-2}$ & min err & err at $k{=}25$ \\
\midrule
Schwarz WDBC, $\delta{=}0.05$ & $11$ \,($2/3$) & -- \,($1/3$) & $1.3{\times}10^{-2}$ & $3.7{\times}10^{-2}$ \\
Schwarz WDBC, $\delta{=}0.10$ & $\mathbf{2}$ & $\mathbf{3}$ & $2.4{\times}10^{-3}$ & $7.1{\times}10^{-3}$ \\
RN ($\alphastar{=}4.16$, closed form) & $13$ & $17$ & $\mathbf{1.8{\times}10^{-3}}$ & $\mathbf{1.8{\times}10^{-3}}$ \\
RN ($\alpha{=}10.31$, E10 search) & $18$ & $22$ \,($2/3$) & $7.6{\times}10^{-3}$ & $7.8{\times}10^{-3}$ \\
DN & \multicolumn{4}{c}{diverges on all seeds} \\
\bottomrule
\end{tabular}
\end{table}

\begin{figure}[h]
\centering
\figmaybe[0.75\linewidth]{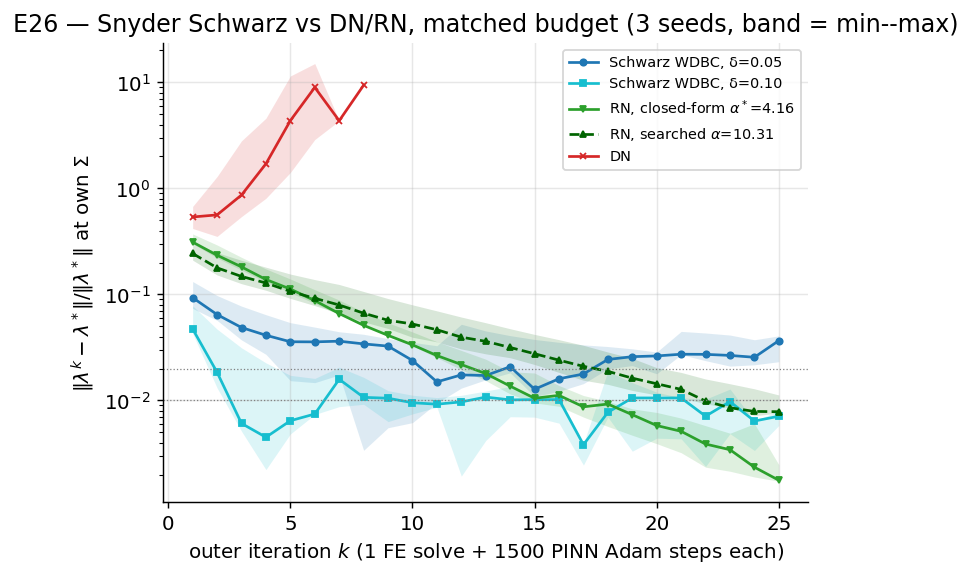}{E26: Snyder Schwarz vs.\ DN/RN convergence trajectories, three seeds.}
\caption{E26 head-to-head at matched budget: relative trace error per
outer iteration (median over three seeds, band $=$ min--max; dotted
lines at the two tolerances of Table~\ref{tab:e26-schwarz}). The
$\delta=0.10$ Schwarz iteration plunges within $2$--$3$ iterations but
plateaus above the RN floor and oscillates; the $\delta=0.05$ variant
stalls at ${\sim}3\times10^{-2}$; RN at the closed-form $\alphastar$
descends monotonically through both tolerances to the lowest floor; DN
diverges.}
\label{fig:app-e6}
\end{figure}

The result (Table~\ref{tab:e26-schwarz}, Figure~\ref{fig:app-e6}) is
informative in both directions. \emph{(i) With a generous overlap the
Schwarz iteration is fast.} At $\delta=0.10$ (40\% of the disc radius)
it reaches $2\times10^{-2}$ in a median of two outer iterations and
$10^{-2}$ in three, against thirteen and seventeen for RN at
$\alphastar$; classical Schwarz theory, with its contraction rate
improving with overlap, predicts exactly this advantage, and the counts
are consistent with the ${\sim}10$-iteration PINN--FOM convergence
reported by \citet{SnyderTezaurWentland2023}. \emph{(ii) The speed is
bought with overlap and robustness.} At $\delta=0.05$ the same scheme is
seed-fragile (one of three seeds reaches $10^{-2}$; one never reaches
$2\times10^{-2}$) and non-monotone: on most seeds the error drifts back
up from its minimum under continued warm-started retraining (median
$1.3\times10^{-2}$ at the minimum vs.\ $3.7\times10^{-2}$ at $k=25$),
the warm-start mechanism quantified in Appendix~\ref{app:E21}. RN
contracts monotonically on every seed, reaches the lowest floor of any
scheme ($1.7$--$2.5\times10^{-3}$), and requires no overlap. The latter
is a requirement, not a preference, in the FSI application of
Section~\ref{sec:fsi}: the subdomains meet at a sharp physical
fluid--solid interface that cannot overlap, and the Dirichlet--Dirichlet
Schwarz iteration converges only for a non-empty overlap (as
\citet{SnyderTezaurWentland2023} note), so the transmission-condition
(DN/RN) family is the one available there. \emph{(iii) A protocol note
in favour of the sweep-free default.} Under this cold-start benchmark
the closed-form $\alphastar=4.16$ outperforms the E10-searched
$\alpha=10.31$ (whose sweep used the warm-started protocol of
Appendix~\ref{app:inner-pinn:e10}): the $\alpha$-bowl is broad and its
empirical minimum protocol-dependent, which strengthens the case for
the closed-form default.

\section{FSI solver figures and FEM--FEM reference meshes}
\label{app:fsi-figs}

\begin{figure}[h]
\centering
\figmaybe[0.95\linewidth]{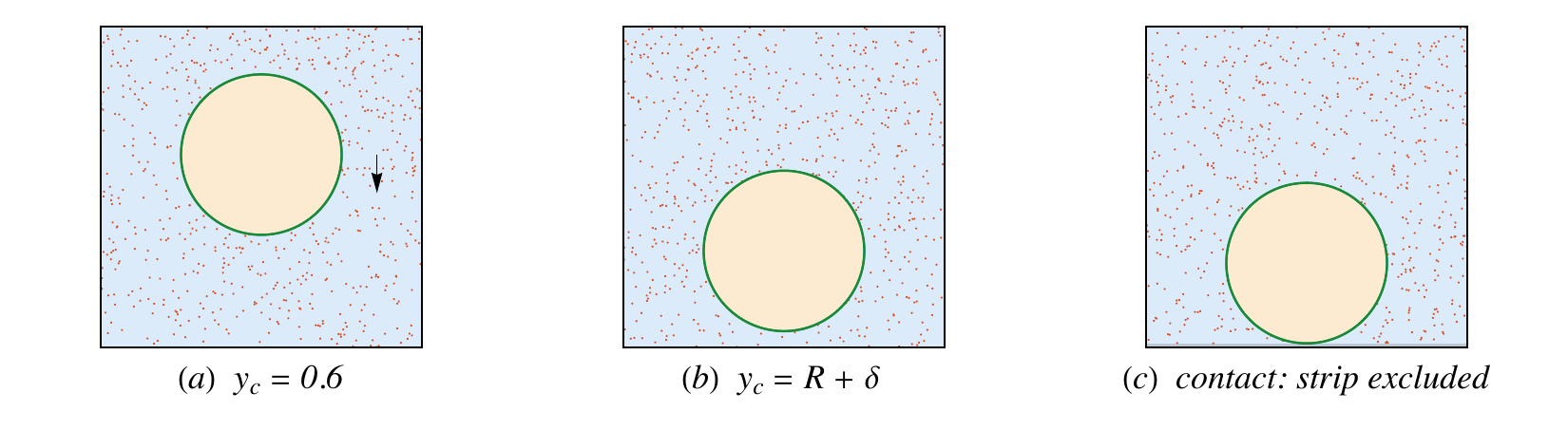}{Three
configurations of the moving fluid topology near contact.}
\caption{Three configurations of the moving fluid domain
$\Omegaf(t)$. (a) free fall, $y_c=0.6$: collocation points fill
$\Omega_+\setminus B$. (b) approach, $y_c=R+\delta$ with $\delta\sim
\varepsilon_g$. (c) contact, $y_c=R+\varepsilon_g$: the strip
$\{y<\varepsilon_g\}$ is excluded from the collocation set (shaded). The
PINN never sees a degenerate gap because the sampler simply stops
drawing points there; no cut-cell or remeshing operation is required.}
\label{fig:contact-topology}
\end{figure}

\begin{figure}[h]
\centering
\figmaybe[0.7\linewidth]{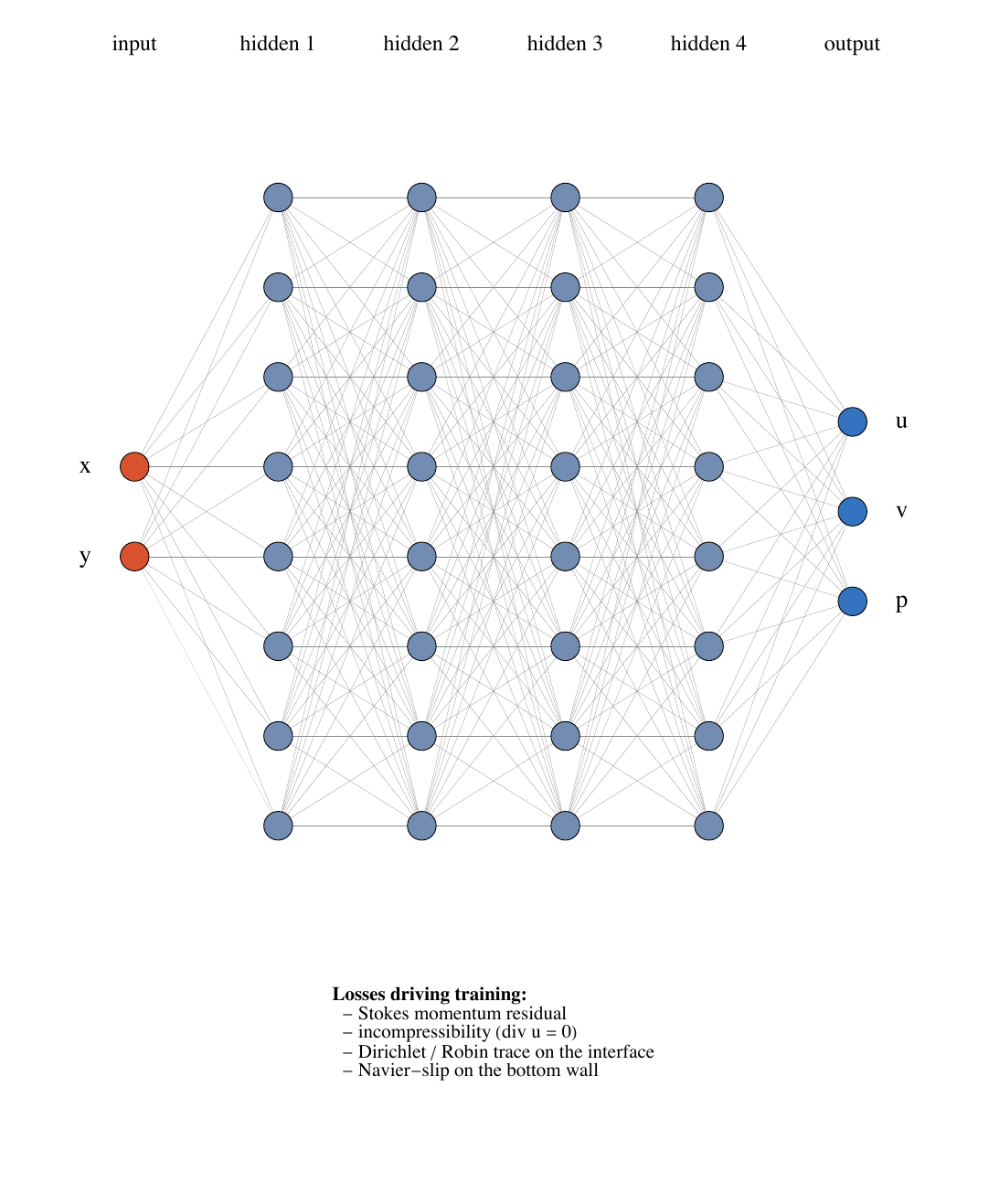}{PINN-Stokes MLP architecture.}
\caption{Architecture of the PINN-Stokes network of
Section~\ref{sec:fsi}: a fully-connected MLP with inputs $(x,y)$, four
hidden Tanh layers of width $32$, and outputs $(u,v,p)$. The losses combine
the Stokes residual, the divergence-free constraint, the no-slip
Dirichlet walls, the Robin condition on $\SigmaI$, and the Navier-slip
condition on $\Sigma_p$.}
\label{fig:nn-arch}
\end{figure}

\begin{figure}[h]
\centering
\figmaybe[0.7\linewidth]{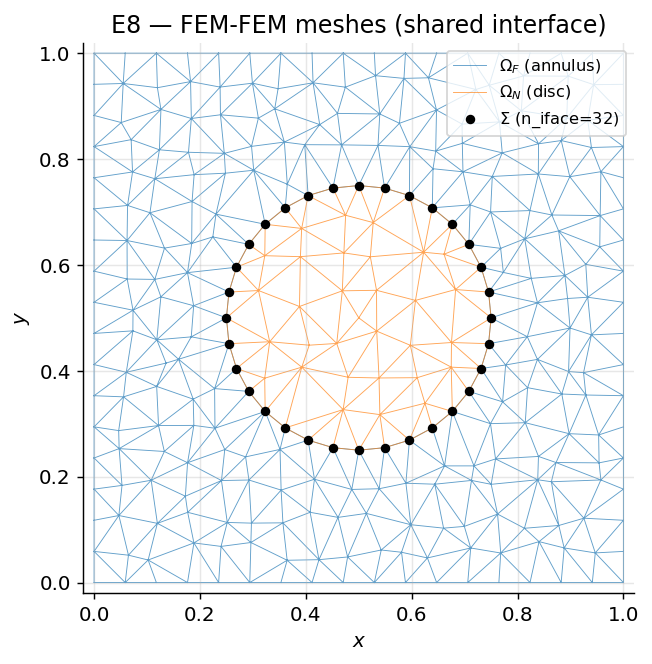}{FEM-FEM reference meshes (E8).}
\caption{FEM--FEM reference meshes used in E8 (Poisson, natural
orientation) and, with the FSI-aligned role assignment, in E9-FEMFEM.
Annulus and disc FE meshes share the same $n_{\mathrm{iface}}=32$
interface vertices on $\SigmaI$. The disc FE mesh is the one used by the
rigid-solid side of the FSI experiments; the annulus FE mesh serves only
as the FEM--FEM Steklov reference and is replaced by the PINN's mesh-free
collocation in every live FSI run.}
\label{fig:app-femfem-meshes}
\end{figure}

This appendix collects the structural figures for the FSI application of
Section~\ref{sec:fsi}. Figure~\ref{fig:contact-topology} details the
three configurations of the moving fluid topology near contact,
Figure~\ref{fig:nn-arch} gives the PINN-Stokes
network architecture (the block-level view of the partitioned solver
that implements Algorithm~\ref{alg:contact-fsi} is
Figure~\ref{fig:fsi-coupling} in the introduction),
and Figure~\ref{fig:app-femfem-meshes} the FE meshes used both for the E8
Poisson FEM--FEM reference and, with the FSI-aligned role assignment, for
the E9-FEMFEM Steklov benchmark.

\section{Additional E12 free-fall diagnostics}
\label{app:E12-extra}

\begin{figure}[h]
\centering
\figmaybe[0.75\linewidth]{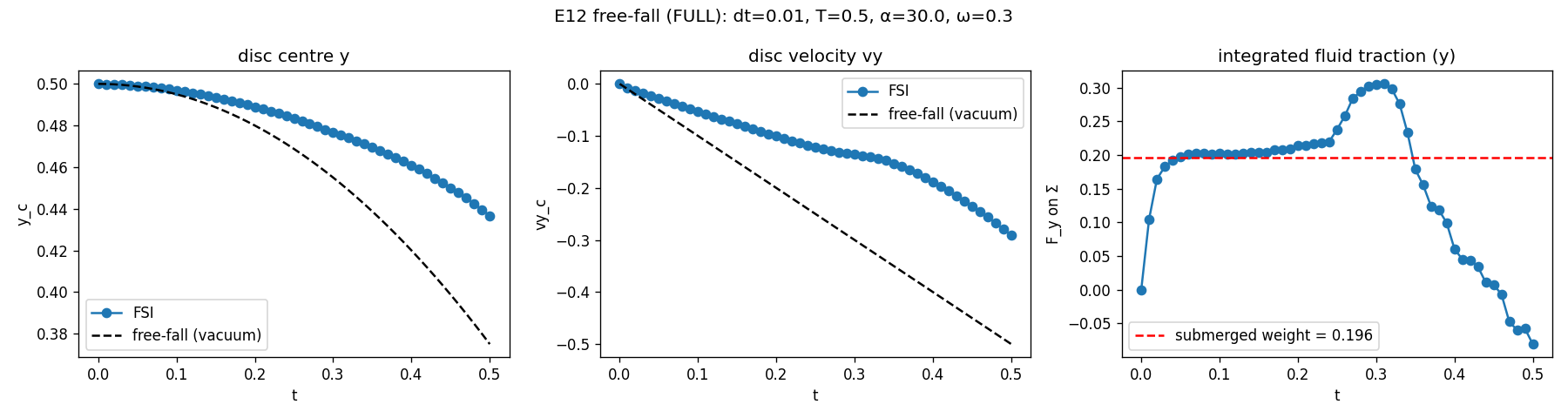}{E12 free-fall trajectory: $y_c(t)$, $\dot y_c(t)$, $F_y(t)$.}
\caption{E12 free-fall: $F_y$ matches the buoyancy level $\rho_f\pi R^2
g=0.196$ through $t\approx0.25$ (drag-balanced fall), then warm-start
drift breaks the iteration. The Newton-balance fit
$m\ddot y_c+\Pi=c|\dot y_c|$ on $t\in[0.04,0.25]$ gives effective drag
$c\approx0.32$, consistent with a confined sub-terminal disc.}
\label{fig:e12}
\end{figure}

\begin{figure}[h]
\centering
\figmaybe[0.95\linewidth]{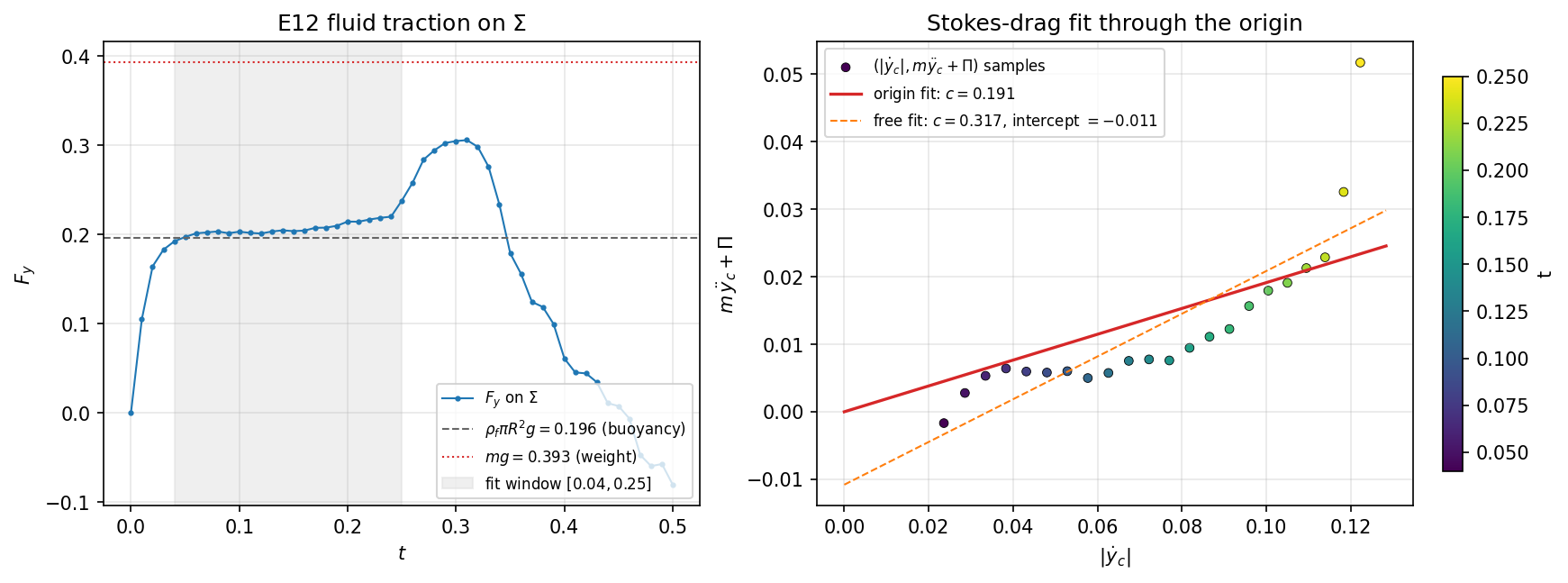}{E12 Stokes-drag fit.}
\caption{Quantitative drag-balance diagnostic. Left: $F_y$ plateaus at
the buoyancy level $0.196$ in the fit window $[0.04,0.25]$. Right: Newton
balance $m\ddot y_c+\Pi$ vs.\ $|\dot y_c|$; origin and free fits give
$c\in[0.19,0.32]$ with near-zero intercept. The implied terminal
velocity $\Pi/c\approx1.0$ is well beyond the observed maximum $|\dot
y_c|=0.12$, so the disc is sub-terminal throughout E12.}
\label{fig:e12-drag}
\end{figure}

These two figures expand the free-fall (E12) result of
Section~\ref{sec:fsi:contact}. Figure~\ref{fig:e12} shows the trajectory
and the buoyancy-level $F_y$ plateau before warm-start drift sets in, and
Figure~\ref{fig:e12-drag} the quantitative Newton drag-balance fit that
places the disc in the confined sub-terminal regime.

\section{Reproducibility}
\label{app:repro}

The complete implementation is a single pip-installable Python package
(Python 3.14, PyTorch 2.11, scikit-fem 12.0, numpy 2.4); the code,
together with the cached artefact behind every figure, will be released
publicly alongside the published version of this paper. All experiments
are deterministic given the seeds set at the top
of each driver. The full 1D + 2D Poisson sweep reproduces in ${\sim}30$
minutes on a 4-thread CPU; the full FSI sweep (E11--E13f) in ${\sim}90$
minutes. The heavier HPC batteries (the E21 drift-remedy, E22
inner-solve-floor, E23 random-feature, and E24 gradient-descent
scale-up sweeps) were run on a single
compute node with four NVIDIA H100 accelerators and an InfiniBand
interconnect, under a TOSS~4 image and the Slurm scheduler; the headline
1D/2D Poisson, static-Stokes, free-fall, and static-equilibrium contact
results, the E25 training-budget floor sweep, and the E26 Schwarz
head-to-head all reproduce on the laptop CPU.

\section{Glossary}
\label{app:glossary}

\begin{table}[h]
\centering
\small
\begin{tabular}{@{}ll@{}}
\toprule
Symbol & Meaning \\
\midrule
$\OmegaF$, $\OmegaN$ & FEM subdomain, PINN subdomain \\
$\SigmaI$ & shared interface $\overline{\OmegaF}\cap\overline{\OmegaN}$ \\
$\SF$, $\SN$ & Steklov--Poincar\'e operators on the two subdomains \\
$\sigmaF, \sigmaN$ & eigenvalues of $\SF, \SN$ \\
$\alphastar$ & spectrally-derived Robin impedance, $\sqrt{\sigmaF^{\min}\,\sigmaN^{\max}}$ \\
$\TDN, \TRN$ & DN and RN trace-error iteration maps \\
$\hSN$, $\EN$ & realised (perturbed) PINN Steklov operator, $\hSN=\SN+\EN$ \\
$\epsN$ & per-step PINN approximation-plus-optimisation error (Theorem~\ref{thm:pinn-floor}) \\
$\That$, $\rhohat$ & realised RN iteration operator and its norm \\
$\Ll_{\tagN}$ & achieved PINN training loss (controls the floor, Remark~\ref{rem:loss-floor}) \\
$\Tinner$, $\Kouter$ & inner training budget per outer iteration, outer-iteration budget \\
$\Miface$ & interface mass matrix (pointwise flux conversion) \\
$\gint$, $\gpw$ & integrated FE reaction and its pointwise conversion $\gpw=\Miface^{-1}\gint$ \\
$\mathrm{Ma}$ & added-mass number, $\sigmaF(q_1)/\sigmaN(q_1)$ \\
AL & augmented Lagrangian (contact projector) \\
$\varepsilon_g$ & relaxed gap parameter in AL \\
$\gamma_c$ & AL penalty stiffness \\
$\omega$ & RN under-relaxation \\
$\omega_{\mathrm{AL}}$ & AL over-relaxation \\
$\beta$ & Navier-slip length on $\Sigma_p$ \\
$\Pi$ & submerged weight, $(\rho_s-\rho_f)\pi R^2 g$ \\
\bottomrule
\end{tabular}
\end{table}

\end{document}